\documentclass[draft,seceqn,secthm]{elsart}

\usepackage{amssymb}\usepackage{amsmath}
\usepackage{hyperref}
\usepackage{mathrsfs}

\def \dist{\operatorname {dist}}
\def \supp{\operatorname {supp}}
\def \ord{\operatorname {ord}}
\begin{document}

\begin{frontmatter}



\title{Weighted Hardy-Sobolev Spaces and\\
Complex Scaling of Differential Equations\\ with Operator
Coefficients}


\author{Victor Kalvin\thanksref{AKA}}
\thanks[AKA]{This work was funded by grant number
108898 awarded by the Academy of Finland.} \ead{vkalvin@gmail.com}
\address{Department of Mathematical Information Technology, University of Jyv\"{a}skyl\"{a}, P.O. Box 35 (Agora),
 FIN-40014 University of Jyv\"{a}skyl\"{a}, Finland}

\date{}

\begin{abstract}
In this paper we study weighted Hardy-Sobolev spaces of vector
valued
 functions analytic on double-napped cones of the complex plane.
We introduce these spaces as a tool for complex scaling of linear
ordinary differential equations with dilation analytic unbounded
operator coefficients. As examples we consider boundary value
problems in cylindrical domains and domains with quasicylindrical
ends.
\end{abstract}

\begin{keyword}
complex scaling\sep quasicylindrical ends\sep cylindrical ends\sep
cusps\sep dilation analytic\sep Hardy classes\sep Hardy-Sobolev
spaces\sep operator coefficients
\end{keyword}
\end{frontmatter}


 \tableofcontents
\section{Introduction}
\subsection{Motivation and structure of the paper}
Complex scaling method has a long history in mathematical physics,
going back to the original work of Combes~\cite{combes}. The idea
was developed by Aguilar,  Balslev,  Combes, and  Simon in
\cite{aguilar combes,balslev combes,Simon}. Van Winter
independently studied the complex scaling method invoking
 a technique of spaces of analytic functions~\cite{Van
Winter I,Van Winter}, however this approach did
 not get ensuing development. Further substantial steps
 were maid by Simon~\cite{simon}, Hunziker~\cite{Hunziker},
 G\`{e}rard~\cite{Gerard},
  a microlocal approach to the complex
scaling was devised by  Helffer and Sj\"{o}strand~\cite{helffer
sjostrand};  for a historical account see e.g. \cite{Hislop
Sigal}.  More recently, a  version of the microlocal
approach~\cite{helffer sjostrand} appeared in
the paper~\cite{J. Wunsch and M. Zworski} by Wunsch and Zworski. 
Mazzeo and Vasy \cite{R+1} implemented the complex scaling for studying the Laplacian on symmetric spaces.

The complex scaling method is widely used in mathematical physics
as an effective tool for meromorphic continuation of the
resolvent. Beyond that point, the method
 is finding increasing use in the construction of artificial
perfectly matched absorbing layers for reliable numerical analysis
of problems of mathematical physics, see
e.g.~\cite{Bonnet1,ref3,ref4,ref5,PREPRINT}; here the original
idea is due to
 B\'{e}renger~\cite{berenger}. In spite of attracting a lot
of attention, the method  of complex scaling itself is  poorly
explored for boundary value problems in domains with
quasicylindrical ends (quasicylindrical ends are unbounded domains
that can be transformed in a neighbourhood of infinity to a
half-cylinder by  suitable diffeomorphisms). To the author's
knowledge, in this setting there are only results obtained by
applying the modal expansion method in cylindrical ends. The modal
expansions
   immediately justify the complex scaling,
 see e.g.~\cite{arsenev},
  but unfortunately the class of
 problems, where modal
 expansions  can be employed, is quite  narrow. On the other hand,
 a  theory of elliptic boundary value problems
 in domains with quasicylindrical ends or with cusps  was developed
by Feigin \cite{Feigin}, Bagirov and Feigin~\cite{BagirovFeigin},
Maz'ya and Plamenevskii \cite{MP1,MP2}, see also Kozlov et
al~\cite{R17,ref7} and references therein. Nevertheless analytic
properties of solutions that may be used for the complex scaling
are beyond the scope of this theory. We first need 
analyticity of solutions with respect to a local coordinate in a
 neighbourhood of infinity. In the case of a
cylindrical end, this coordinate is the axial coordinate of the
cylinder, for a quasicylindrical end it is a curvilinear
coordinate along the quasicylinder.  It is natural to assume that
 the coefficients of operators, the right hand
side of the problem, and the boundary are analytic with respect to
only that local coordinate. In particular, we will be interested in
behavior of analytic continuations of solutions as the coordinate
tends to infinity in a complex sector.


In order to study analytic properties of solutions near infinity, we implement a special localization procedure in every quasicylindrical end. This procedure preserves the analyticity of solutions in the neighborhoods of localization. The localized solutions satisfy some subsidiary problems in infinite cylinders, our task reduces to studying of complex scaling for the subsidiary problems.  In many instances one can
treat the subsidiary boundary value problems in infinite cylinders
as linear ordinary differential equations with unbounded operator
coefficients. We develop the complex scaling for these equations.

In many approaches to the complex scaling, original and  scaled operators are studied at first separately, and then the relations between them are clarified with the help of a sufficiently large set of analytic functions. Our approach is principally different: we consider the original and scaled operators as one operator acting in  spaces of analytic functions.
Namely, we study differential operators with unbounded operator coefficients in Hardy-Sobolev spaces of
analytic functions. Here we essentially rely on methods of the mentioned above theory of elliptic boundary value problems. In this paper we are mainly concerned with the Hardy-Sobolev spaces. Exploring the complex scaling, we  restrict ourselves to a
relatively simple case, assuming that the values of parameters are limited so that the operators are Fredholm. A more general case will be considered
elsewhere.



The structure of the paper is as follows. We start with a detailed outline of the paper in Section~\ref{ss2}, where we formulate our results in a simplified form and illustrate them by some examples. The main text consists of three parts
(Sections~\ref{sec1}--\ref{sec3}) plus Appendix.
  Section~\ref{sec1} is devoted to a close examination of weighted Hardy spaces of
functions analytic on double-napped cones of the complex plane. In
these spaces we consider Fredholm polynomial operator pencils, Section~\ref{OP}.
In Section~\ref{sec2}, by means of the Fourier-Laplace transformation, we pass from the Hardy spaces to the
Hardy-Sobolev spaces and study properties of the latter. Finally, results
of Sections~\ref{sec1} and \ref{sec2} are applied in
Section~\ref{sec3}, where we investigate linear ordinary
differential equations with unbounded operator coefficients. Every
section is equipped with a short outline. Some technical proofs of Section~\ref{sec1} are located in Appendix.




\subsection{Outline and examples of applications}\label{ss2}
Let $\varphi\in(0,\pi]$ and let $e^{i\varphi}\Bbb R$ stand for the
line $\{\lambda\in\Bbb C:\Im\lambda= \Re\lambda\tan \varphi \}$ in
the complex plane $\Bbb C$. By  $\mathcal
K^\varphi=\{\lambda\in\Bbb C: \lambda=e^{i\psi}\xi, \xi\gtrless 0,
\psi\in(0,\varphi)\}$ we denote the open double-napped cone. We
define the Hardy class $\mathscr H(\mathcal K^\varphi;X)$ of
  analytic functions $\mathcal K^\varphi\ni\lambda\mapsto \mathcal
F(\lambda)\in X$ with values in a Hilbert space $X$ as a class of
all the functions $\mathcal F$ satisfying the estimate $\|\mathcal
F; L_2(e^{i\psi}\Bbb R;X)\|\leq Const(\mathcal F)$ for
$\psi\in(0,\varphi)$. Here $L_2(e^{i\psi}\Bbb R;X)$ is the space
of all square summable $X-$valued functions on $e^{i\psi}\Bbb R$.
The functions $\mathcal F\in \mathscr H(\mathcal K^\varphi;X)$ can
be extended to almost all points of the boundary $\partial
\mathcal K^\varphi=e^{i\varphi}\Bbb R\cup \Bbb R$ by the
non-tangential limits (i.e. by the limits of $\mathcal F(\lambda)$
 in $X$ as $\lambda$ goes to $\mu\in\partial\mathcal K^\varphi$ in $\mathcal K^\varphi$ by a
non-tangential to $\partial \mathcal K^\varphi$ path). The
extended functions satisfy the inclusions $\mathcal F\in L_2(\Bbb
R; X)$ and $\mathcal F\in L_2(e^{i\varphi}\Bbb R; X)$; in the case
$\varphi=\pi$ we distinguish the banks $\lim_{\psi\to 0+}
e^{i\psi}\Bbb R$ and $\lim_{\psi\to \varphi-} e^{i\psi}\Bbb R$ of
the  boundary $\partial \mathcal K^\pi$. The Hardy class $\mathscr
H(\mathcal K^\varphi;X)$ endowed with the norm
$$
\|\mathcal F; \mathscr
H(\mathcal K^\varphi;X)\|=\|\mathcal F; L_2(\Bbb R;
X)\|+\|\mathcal F; L_2(e^{i\varphi}\Bbb R; X)\|
$$
 is a Banach space. For any $\psi\in[0,\varphi]$ we can identify the elements $\mathcal F\in\mathscr H(\mathcal
 K^\varphi;X)$ with their traces $\mathcal F\!\!\upharpoonright_{e^{i\psi}\Bbb
 R}\in L_2(e^{i\psi}\Bbb R; X)$, then $\mathscr H(\mathcal
 K^\varphi;X)$
is dense  in $L_2(e^{i\psi}\Bbb R; X)$. As is known, the Fourier
transformation
 yields an isometric isomorphism of the space
$L_2(\Bbb R; X)$ onto itself. It turns out that the Fourier
transformation also yields an isometric isomorphism between
$\mathscr H(\mathcal K^\varphi;X)$ and a  Hardy class $\mathsf
H(K^\varphi; X)$ in the dual double-napped cone
$K^\varphi=\{z\in\Bbb C: z=e^{-i\psi}t, t\lessgtr 0,
\psi\in(0,\varphi)\}$; here the elements of the Hardy classes are
identified with their traces on $\Bbb R$, the class $\mathsf
H(K^\varphi; X)$ is defined in the same way as $\mathscr
H(\mathcal K^\varphi;X)$. Note that spaces similar to $\mathscr
H(\mathcal K^\varphi;X)$ and $\mathsf H(K^\varphi; X)$ were
considered in context of the complex scaling method by Van Winter
\cite{Van Winter I,Van Winter}.

We introduce the weighted Hardy class $\mathscr H^\ell(\mathcal
K^\varphi;X)$, $\ell\in\Bbb R$, that consists of all functions
 analytic in the cone $\mathcal K^\varphi$ and
satisfying the uniform in $\psi\in(0,\varphi)$ estimate
\begin{equation}\label{intro1}
\|(1+|\cdot|^2)^{\ell/2}\mathcal F(\cdot);L_2(e^{i\psi}\Bbb R;
X)\|\leq Const(\mathcal F).
\end{equation}
 Let $\mathscr
W^\ell(e^{i\psi}\Bbb R; X)$ denote the weighted $L_2$-space with
the norm equal to the left hand side of the
estimate~\eqref{intro1}. The class  $\mathscr H^\ell(\mathcal
K^\varphi;X)$ supplied with the norm
$$
\|\mathcal F;\mathscr H^\ell(\mathcal K^\varphi;X)\|=\|\mathcal
F;\mathscr W^\ell(\Bbb R; X)\|+\|\mathcal F;\mathscr
W^\ell(e^{i\varphi}\Bbb R; X)\|
$$
 is a Banach space. We extend the functions $\mathcal F\in \mathscr H^\ell(\mathcal K^\varphi;X)$ to almost all points of the boundary $\partial
\mathcal K^\varphi$ by the non-tangential limits. Here again we
can identify the elements $\mathcal F\in\mathscr H^\ell(\mathcal
K^\varphi;X)$ with their traces $\mathcal
F\!\!\upharpoonright_{e^{i\psi}\Bbb R}$ for any
$\psi\in[0,\varphi]$, then $\mathscr H^\ell(\mathcal K^\varphi;X)$
is dense in $\mathscr W^\ell(e^{i\psi}\Bbb R; X)$ . As is known,
the Fourier transformation  maps  $\mathscr W^\ell(\Bbb R; X)$ to
the Sobolev space $\mathsf W^\ell(\Bbb R;X)$ of tempered
distributions with values in $X$; see e.g. \cite{R4}. Clearly the
Fourier transformation also maps $\mathscr H^\ell(\mathcal
K^\varphi;X)$ to some subset $\mathsf H^\ell(K^\varphi; X)$ of
$\mathsf W^\ell(\Bbb R;X)$, $\mathsf H^\ell(K^\varphi; X)$ is
dense in $\mathsf W^\ell(\Bbb R;X)$. It turns out that $\mathsf
H^\ell(K^\varphi; X)$ is the Hardy-Sobolev space of order $\ell$.
In the main text we suppose that $\ell$ is a real number, but in
this introductory part  we restrict ourselves to the case of
nonnegative integer $\ell$. In this case the Hardy-Sobolev space
$\mathsf H^\ell(K^\varphi; X)$ can be introduced as a space of all
functions
 $F\in \mathsf
H(K^\varphi; X)$ such that the complex derivatives $D^j_z F$,
$j=1, \dots,\ell$, are also in the space $\mathsf H(K^\varphi;
X)$. We define the norm in $\mathsf H^\ell(K^\varphi; X)$ by the
equality
$$
\|F;\mathsf H^\ell(K^\varphi; X)\|=\sum_{j=0}^\ell\|D_z^j
F;\mathsf H(K^\varphi; X)\|.
$$
For any $\psi\in[0,\varphi]$ the space $\mathsf H^\ell(K^\varphi;
X)$ can be viewed as a space of functions on the line
$e^{-i\psi}\Bbb R$, then $\mathsf H^\ell(K^\varphi; X)$ is dense
in the Sobolev space $\mathsf W^\ell(e^{-i\psi}\Bbb R; X)$. Note
that for $\ell\geq 1$ the functions from $\mathsf
H^\ell(K^\varphi; X)$ are in the class
$C^{\ell-1}(\overline{K^\varphi};X)$ of $\ell-1$ times
continuously differentiable functions in the closure
$\overline{K^\varphi}$. For a fixed $\psi\in[0,\varphi]$ the
Fourier  transformation
$$
F(z)=\frac 1 {\sqrt{2\pi}}\int_{e^{i\psi}\Bbb R} e^{iz\lambda}
\mathcal F(\lambda)\,d\lambda, \quad z\in e^{-i\psi}\Bbb R,
$$
defined first on the smooth compactly supported functions
$e^{i\psi}\Bbb R\ni\lambda\mapsto\mathcal F(\lambda)\in X$, can be
extended to an isomorphism between  $\mathscr W^\ell(e^{i\psi}\Bbb
R;X)$ and $\mathsf W^\ell(e^{-i\psi}\Bbb R; X)$, and also  between
$\mathscr H^\ell(\mathcal K^\varphi; X)$ and $\mathsf
H^\ell(K^\varphi; X)$.

Boundary value problems in infinite cylinders and differential
equations with operator coefficients are traditionally studied  in
a scale of Sobolev spaces with exponential weights;
see~\cite{R17,ref7} and references therein. Therefore we introduce
the exponential weight function $e_\zeta:z\mapsto\exp(-i\zeta z)$,
where $\zeta$ is a complex weight number. We also want to shift
the vertex of the cone $K^\varphi$ to a point $w$ of the complex
plane, let $K^\varphi_w=\{z\in \Bbb C: z-w\in K^\varphi\}$. By
$\mathsf H^\ell_\zeta(K^\varphi_w; X)$ we denote the weighted
Hardy-Sobolev space  that consists of all functions $F$ such that
$\bigl(e_\zeta F\bigr )(\cdot-w)\in \mathsf H^\ell(K^\varphi;X)$.
As the norm of $F$ in $\mathsf H^\ell_\zeta(K_w^\varphi; X)$ we
take the value $\|\bigl(e_\zeta F\bigr )(\cdot-w); \mathsf
H^\ell(K^\varphi; X)\|$.  Similarly we define the weighted Sobolev
space $\mathsf W_\zeta^\ell(e^{-i\psi}\Bbb R+w;X)$ that consists
of all functions $F$ such that $\bigl(e_\zeta F\bigr
)(\cdot-w)\in\mathsf W^\ell(e^{-i\psi}\Bbb R;X)$; here
$e^{-i\psi}\Bbb R+w$ denotes the line $\{z\in\Bbb C: z-w\in
e^{-i\psi}\Bbb R\}$. Let us stress that the behavior of the weight
$e_\zeta$ of the space $\mathsf W_\zeta^\ell(e^{-i\psi}\Bbb
R+w;X)$ depends on the angle $\psi$. For a fixed
$\psi\in[0,\varphi]$ the inverse Fourier-Laplace transformation
\begin{equation}\label{flt}
 \mathcal
F(\lambda)=\frac {1} {2\pi}\int_{e^{-i\psi}\Bbb R+w}e^{-i\lambda
z}F(z)\,dz,\quad \lambda\in e^{i\psi}\Bbb R +\zeta,
\end{equation}
  can be extended to an
isomorphism between $\mathsf W_\zeta^\ell(e^{-i\psi}\Bbb R+w;X)$
and $\mathscr W_w^\ell(e^{i\psi}\Bbb R+\zeta;X)$,  where $\mathscr
W_w^\ell(e^{i\psi}\Bbb R+\zeta;X)$ is the space of functions
$e^{i\psi}\Bbb R+\zeta\ni \lambda\mapsto \mathcal F(\lambda)\in X$
with the finite norm
$$
\|\mathcal F;\mathscr W_w^\ell(e^{i\psi}\Bbb
R+\zeta;X)\|=\|\exp\{iw(\cdot-\zeta)\}\mathcal F(\cdot-\zeta);
\mathscr W^\ell (e^{i\psi}\Bbb R;X)\|.
$$
Let $\mathcal K^\varphi_\zeta=\{\lambda\in\Bbb C: \lambda-\zeta\in
\mathcal K^\varphi\}$ be the cone $\mathcal K^\varphi$ shifted by
the weight number $\zeta$. By $\mathscr H^\ell_w(\mathcal
K^\varphi_\zeta; X)$ we denote the weighted Hardy space that
consists of all functions such that $\exp\{iw\cdot\}\mathcal
F(\cdot-\zeta)\in\mathscr H^\ell(\mathcal K^\varphi;X)$. The norm
in $\mathscr H^\ell_w(\mathcal K^\varphi_\zeta; X)$ is given by
$$
\|\mathcal F;\mathscr H^\ell_w(\mathcal K^\varphi_\zeta;
X)\|=\|\exp\{iw(\cdot-\zeta)\}\mathcal F(\cdot-\zeta); \mathscr
H^\ell(\mathcal K^\varphi;X)\|.
$$
Then the extended Fourier-Laplace transformation~\eqref{flt}
yields an isomorphism between the weighted Hardy-Sobolev space
$\mathsf H^\ell_\zeta(K^\varphi_w; X)\bigl (\subset\mathsf
W_\zeta^\ell(e^{-i\psi}\Bbb R+w;X)\bigr)$ and the weighted Hardy
space $\mathscr H^\ell_w(\mathcal K^\varphi_\zeta; X)\bigl(\subset
\mathscr W_w^\ell(e^{i\psi}\Bbb R+\zeta;X)\bigr )$. This enables
us to adapt methods of the theory  of ordinary differential
equations with
 unbounded operator coefficients  developed in the scale of
Sobolev spaces $\mathsf W_\zeta^\ell(\Bbb R;X)$ (see
e.g.~\cite{R17}) to the case of the Hardy-Sobolev spaces $\mathsf
H^\ell_\zeta(K^\varphi_w; X)$.

The main purpose of this paper is to study the weighted  spaces
$\mathscr H^\ell_w(\mathcal K^\varphi_\zeta; X)$ and $\mathsf
H^\ell_\zeta(K^\varphi_w; X)$. Nevertheless, 
we demonstrate how to treat the complex scaling of differential
equations with unbounded operator coefficients in terms of the
weighted Hardy-Sobolev spaces. We also give  examples of
applications to boundary value problems in cylindrical domains
 and in domains with quasicylindrical end. 
 Before presenting our results on the complex scaling, we give some
preliminaries.

Let $X_j$ denote a Hilbert space with the norm $\|\cdot\|_j$. We
introduce a set $ \{X_j\}_{j=0}^m$ of Hilbert spaces such that
$\|u\|_{j}\leq \|u\|_{j+1}$ and
   $X_{j+1}$ is dense in $X_{j}$ for all $j=0,\dots,m-1$. Let
   $\{A_j\in\mathscr B(X_j,X_0)\}_{j=0}^m$ be a set of operators, where
$\mathscr B(X_j,X_0)$ stands for the set of all linear bounded
operators $A_j:X_j\to X_0$. We start with the differential
equation with constant operator coefficients
\begin{equation}\label{equation1}
\sum_{j=0}^m A_{m-j}D_t^ju(t)=F(t), \quad t\in\Bbb R,
\end{equation}
where $D_t=-i\partial_t$.
 Assume that the operator $\mathfrak A(\lambda)=\sum_{j=0}^m A_{m-j}\lambda^j$
 from  $\mathscr B(X_m,X_0)$ is Fredholm
for all $\lambda\in\Bbb C$ and is invertible for at least one
value of $\lambda$.  Under these assumptions the operator
$\mathfrak A(\lambda)$ is invertible for all $\lambda\in\Bbb C$
except for isolated eigenvalues of the operator pencil $\Bbb
C\ni\lambda\mapsto\mathfrak A(\lambda) \in \mathscr B(X_m,X_0)$.
These eigenvalues are of finite algebraic multiplicities and can
accumulate only at infinity. We also assume that  there exists
$R>0$ such that
 for all $f\in X_0$ the estimate
\begin{equation}\label{M cond}
\sum_{j=0}^m |\lambda|^j\|\mathfrak A^{-1}(\lambda)f\|_{m-j}\leq c
\|f\|_0, \quad \lambda\in\Bbb R,\  |\lambda|>R,
\end{equation}is fulfilled.
The assumptions we made are widely met in the theory of
differential equations with operator coefficients, they are
satisfied in many applications to boundary value
  problems for partial differential equations; see e.g.~\cite{R17,ref7} and  references
 therein. The assumption~\eqref{M cond} in particular guaranties the existence of an angle
 $\vartheta\in(0,\pi/2)$ such that
  for any $\zeta\in\Bbb C$ and any $\varphi\in(0,\vartheta)$
 the closed dual cone $\overline{\mathcal K^\varphi_\zeta}$ contains at most finitely
 many eigenvalues of the operator pencil $\mathfrak A$, the estimate~\eqref{M cond}
 remains valid for all $\lambda\in \overline{\mathcal K^\varphi_\zeta}$, $|\lambda|>R$; see~\cite[Proposition~2.2.1]{R17}. One can find $\zeta\in\Bbb C$ and $\varphi\in(0,\vartheta)$
 so that $\overline{\mathcal K^\varphi_\zeta}$
 is  free from the
 spectrum of $\mathfrak A$. In this paper dealing with complex scaling we restrict ourselves by the assumption that $\overline{\mathcal K^\varphi_\zeta}$
 is  free from the
 spectrum of $\mathfrak A$, a more general case  will be considered elsewhere.

 In the remaining part of this subsection we assume that
 $\varphi\in(0,\vartheta)$.

 We introduce two scales of Banach spaces
\begin{equation}\label{space1}
\mathsf D^m_\zeta( \Bbb R)=\bigcap_{j=0}^m \mathsf
W^{m-j}_\zeta(\Bbb R; X_j),\quad \|u; \mathsf D^m_\zeta( \Bbb
R)\|=\sum_{j=0}^m \|u;\mathsf W^{m-j}_\zeta( \Bbb R; X_j)\|;
\end{equation}
\begin{equation}\label{space2}
\mathsf D^m_\zeta( K_0^\varphi)=\bigcap_{j=0}^m \mathsf
H^{m-j}_\zeta(K_0^\varphi; X_j),\ \|u; \mathsf
D^m_\zeta(K_0^\varphi)\|=\sum_{j=0}^m \|u;\mathsf
H^{m-j}_\zeta(K_0^\varphi; X_j)\|.
\end{equation}
We can identify the functions $u\in\mathsf D^m_\zeta(
K_0^\varphi)$ with their non-tangential boundary limits
$u\!\!\upharpoonright_{\Bbb R}\in\mathsf D^m_\zeta( \Bbb R)$. The
space $\mathsf D^m_\zeta( K_0^\varphi)$ viewed as a space of
functions on $\Bbb R$ is dense  in $\mathsf D^m_\zeta( \Bbb R)$,
the space $\mathsf D^m_\zeta( \Bbb R)$ does not depend on
$\Re\zeta$. As is well-known~\cite[Theorem 2.4.1 and Remark
2.4.2]{R17}, the operator $\mathfrak A(D_t):\mathsf D^m_\zeta(
\Bbb R)\to \mathsf W^{0}_\zeta( \Bbb R;X_0)$ of the
equation~\eqref{equation1} yields an isomorphism if and only if
the line $\Bbb R+\zeta$ is free from the eigenvalues of the pencil
$\Bbb C\ni\lambda\mapsto\mathfrak A(\lambda) \in \mathscr
B(X_m,X_0)$; if there is an eigenvalue of the pencil on the line
$\Bbb R+\zeta$ then the range of the operator $\mathfrak
A(D_t):\mathsf D^m_\zeta( \Bbb R)\to \mathsf W^{0}_\zeta( \Bbb
R;X_0)$ is not closed.

 Note that  if the line $\Bbb R+\zeta$ is free from
the eigenvalues of the pencil $\mathfrak A$ then for a
sufficiently small angle $\varphi$ the closed dual cone
$\overline{\mathcal K^\varphi_\zeta}$ is also free from the
eigenvalues of the pencil $\mathfrak A$. Now we are in position to
formulate some results on the complex scaling.

\begin{thm}\label{T1} Assume that the closed dual cone $\overline{\mathcal
K^\varphi_\zeta}$ is free from the eigenvalues of the pencil $\Bbb
C\ni\lambda\mapsto\mathfrak A(\lambda) \in \mathscr B(X_m,X_0)$.
Then the following assertions hold. {\rm (i)} The operator
$\mathfrak A(D_t):\mathsf D^m_\zeta( K_0^\varphi)\to \mathsf
H^{0}_\zeta(K_0^\varphi;X_0)$ of the equation~\eqref{equation1}
yields an isomorphism; here the analytic in $K_0^\varphi$
functions $u\in \mathsf D^m_\zeta( K_0^\varphi)$, $F\in \mathsf
H^{0}_\zeta(K_0^\varphi;X_0)$ are identified with their boundary
limits $u\!\!\upharpoonright_{\Bbb R}\in \mathsf D^m_\zeta(\Bbb
R)$, $F\!\!\upharpoonright_{\Bbb R}\in\mathsf W^0_\zeta(\Bbb
R;X_0)$. {\rm (ii)}  Let  $u\in \mathsf D^m_\zeta( \Bbb R)$ be a
unique solution to the equation~\eqref{equation1} with right hand
side $F\in \mathsf H^{0}_\zeta(K_0^\varphi;X_0)$. Then $u\in
\mathsf D^m_\zeta( K_0^\varphi)$ and the function $\Bbb R\ni
t\mapsto v(t)\equiv u(e^{-i\varphi}t)\in X_m$ is a unique solution
$v\in\mathsf D_{e^{-i\varphi}\zeta}^{m}(\Bbb R)$ to the scaled
equation $\mathfrak A(e^{i\varphi}D_t)v(t)=F(e^{-i\varphi}t)$,
$t\in\Bbb R$.
\end{thm}
As an example of application of Theorem~\ref{T1} we consider the
Dirichlet problem in the cylinder $\mathcal C=\{(y,t):
y\in\Omega,t\in\Bbb R\}$, where  $\Omega$ is a domain in $\Bbb
R^{n-1}$ with compact closure and a smooth boundary. We introduce
the differential operator  $\mathscr L(y,D_y,D_t)=\sum_{j=0}^m
A_{m-j}(y,D_y) D_t^j$, where $A_{m-j}(y,D_y)$ are differential
operators with smooth coefficients. Consider the Dirichlet problem
\begin{equation}\label{example1}
\begin{aligned}
&\mathscr L(y,D_y,D_t)u(y,t)=F(y,t),\quad (y,t)\in\mathcal
C,\\&\partial_\nu^j u(y,t)=0,\quad (y,t)\in \partial \mathcal C,\
j=0,\dots,k-1;
\end{aligned}
\end{equation}
here  $\partial_\nu=\partial/\partial\nu$ and $\nu$ is the outward
normal. Suppose that $2k/m$ is an integer and the order of $A_{j}$
does not exceed $2kj/m$ for $j=0,\dots,m$. To the
problem~\eqref{example1} there corresponds the
 operator pencil
\begin{equation}\label{pencilEx}
\Bbb C\ni\lambda\mapsto \mathfrak A(\lambda)=\mathscr L(y,D_y,
\lambda)\in \mathscr B(X_m,X_0),
\end{equation}
 where $X_0= L_2(\Omega)$, $X_m=\{U\in W^{2k}(\Omega): \partial_\nu^j U=0 \text{ on }\partial \Omega, j=0,\dots,
 k-1\}$;
 here $W^{2k}(\Omega)$ stands for the Sobolev space of functions in
 $\Omega$.  Assume that the operator $\mathscr
L(y,D_y,D_t)$ is $(2k,m)-$ elliptic, then the operator $\mathfrak
A(\lambda)$ is Fredholm for all $\lambda\in\Bbb C$, the
estimate~\eqref{M cond} is fulfilled for a sufficiently large
$R>0$ and for all $f\in X_0$; see~\cite[Section 2.5]{R17}, where
the unique solvability of the problem~\eqref{example1} is studied.
Recall that $(2k,2k)-$elliptic operators are standard elliptic
operators of order $2k$, the class of $(2k,1)-$elliptic operators
includes parabolic operators.

We set $X_j=W^{2kj/m}(\Omega)\cap X_m$.
 The coefficients $A_j(y,D_y)$ satisfy the inclusions  $A_j\in\mathscr
 B(X_j;X_0)$. In this case the space $\mathsf W^0_\zeta (\Bbb R; X_0)$ is
 the weighted space of square summable functions with the norm $\|e_\zeta F; L_2(\mathcal
 C)\|$; here $e_\zeta$ is the same exponential weight function as before.  The Hardy space $\mathsf H^0_\zeta(K^\varphi_0; X_0)$
consists of analytic functions $K^\varphi_0\ni z\mapsto F(z)\equiv
F(\cdot,z)\in L_2(\Omega)$, the elements of $\mathsf
H^0_\zeta(K^\varphi_0; X_0)$ are extended to almost all points of
the boundary $\partial K^\varphi_0$. For any angle
$\psi\in[0,\varphi]$ the set of functions $\{\Bbb R\ni t\mapsto
F(e^{-i\psi}t): F\in \mathsf H^0_\zeta(K^\varphi_0; X_0)\}$ is
dense in the space $\mathsf W^0_{e^{-i\psi}\zeta} (\Bbb R; X_0)$ .
The space $\mathsf D^m_\zeta (\Bbb R)$
 introduced in \eqref{space1} consists of all functions $u$ such
 that $e_\zeta D_y^\alpha D_t^j u \in L_2(\mathcal C)$ for $|\alpha|+2kj/m\leq
 2k$ and $\partial_\nu^j u=0$ on $\partial \mathcal C$ for
 $j=0,\dots,k-1$. The elements of the space $\mathsf D^m_\zeta (K^\varphi_0)$ defined in
 \eqref{space2} are analytic functions $K^\varphi\ni z\mapsto
 u(z)\equiv u(\cdot,z)\in
 X_m$ extended to almost all points of the boundary $\partial
 K^\varphi_0$ by the non-tangential limits. For any
 $\psi\in[0,\varphi]$ the set of functions $\{\Bbb R\ni t\mapsto u(e^{-i\psi}t): u\in \mathsf
 D^m_\zeta(K^\varphi_0)\}$ is dense in $\mathsf D^m_{e^{-i\psi}\zeta}(\Bbb
 R)$.

 Suppose that the right hand side $F$ of the boundary value
problem~\eqref{example1}
  is in the Hardy space $\mathsf H^0_\zeta(K^\varphi_0;X_0)$.
   The Dirichlet problem~\eqref{example1}
   has a unique solution  $u\in\mathsf D^m_\zeta(\Bbb R)$ for any $F\in\mathsf W^0_\zeta(\Bbb R; X_0)$ if and only if
   the line $\Bbb R+\zeta$ is free from the spectrum of the
 operator pencil~\eqref{pencilEx}; recall that the spaces $\mathsf D^m_\zeta(\Bbb R)$
  and $\mathsf W^0_\zeta(\Bbb R; X_0)$ do not depend on $\Re\zeta$. If the line $\Bbb R+\zeta$
   is free from the spectrum of $\mathfrak A$ then the results of Theorem~\ref{T1}
   are valid provided that $\varphi$ is a sufficiently small
    angle. By Theorem~\ref{T1} a
 solution $u\in\mathsf D^m_\zeta(\Bbb R)$ belongs to the space $\mathsf D^m_\zeta
 (K^\varphi_0)\subset \mathsf H^0_\zeta(K^\varphi_0;X_0)$; moreover, the function $v(y,t)=u(y,e^{-i\varphi}t)$ is a
 unique solution $v\in \mathsf D^m_{e^{-i\varphi}\zeta}(\Bbb R)$ to
 the  Dirichlet problem
\begin{equation}\label{example1+}
\begin{aligned}
&\mathscr L(y,D_y,e^{i\varphi}D_t)v(y,t)=F(y,e^{-i\varphi}t),\
(y,t)\in\mathcal C,\\&
\partial_\nu^j v(y,t)=0,\ (y,t)\in  \partial \mathcal C,\
j=0,\dots,k-1.
\end{aligned}
\end{equation}
The problem~\eqref{example1+} is related to the original
problem~\eqref{example1} by the complex scaling with the scaling
coefficient $e^{-i\varphi}$. Similarly we can consider complex
scaling of problems in the cylinder $\mathcal
C=\{(y,t):y\in\Omega,t\in\Bbb R\}$, where $\overline{\Omega}$ is a
smooth compact manifold with or without boundary.



Now we consider the differential equation with variable operator
coefficients
\begin{equation}\label{equation2}
\sum_{j=0}^m \bigl(A_{m-j}+Q_{m-j}(t)\bigr)D_t^ju(t)=F(t), \quad
t\in\Bbb R.
\end{equation}
Here the coefficients  $A_j$ are the same as in the
equation~\eqref{equation1}, the coefficients $Q_j(t)$ are operator
functions $\Bbb R\ni t\mapsto Q_j(t)\in \mathscr B(X_j,X_0)$ that
are dilation analytic in the following sense: \vspace{0.2cm}

\noindent{\it {\rm i.} for a large $T>0$ and some $\alpha>0$ the
coefficients $Q_0,\dots,Q_m$ can be extended to holomorphic
operator functions $$ \{z\in\Bbb C: -\alpha\leq\arg(z-T)\leq
0\}\ni z\mapsto Q_j(z)\in\mathscr B(X_j,X_0);$$

\vspace{0.2cm} \noindent {\rm ii.}  the values $\| Q_j(z);
\mathscr B(X_j,X_0)\|$, $j=0,\dots, m$, uniformly tend to zero as
$z$ goes to infinity inside the sector $ \{z\in\Bbb C:
-\alpha\leq\arg(z-T)\leq 0\}$.} \vspace{0.2cm}

 As before we
denote $\mathfrak A(\lambda)=\sum_{j=0}^m A_{m-j}\lambda^j$. The
next theorem presents the results of Section~\ref{devc} in a
simplified form.
\begin{thm}\label{T2} Assume that the right hand
side of the equation~\eqref{equation2} is a   compactly supported
function $\Bbb R\ni t\mapsto F(t)\in X_0$, and there are no
eigenvalues of the operator pencil $\Bbb C\ni\lambda\mapsto
\mathfrak A(\lambda)\in\mathscr B(X_m,X_0)$ in the closed dual
cone $\overline{\mathcal K_\zeta^\varphi}$, where $\zeta\in\Bbb C$
and  $\varphi$ does not exceed the angle $\alpha$ from the
conditions {\rm i}, {\rm ii}. Let $T$ be a sufficiently large
positive number.
 Then a solution $u\in \mathsf D^m_\zeta(\Bbb R)$ to the
equation~\eqref{equation2} can be extended to a holomorphic
function
$$
\{z\in\Bbb C: -\varphi\leq\arg(z-T)\leq 0\}\ni z
\mapsto u(z)\in X_m
$$
 satisfying the uniform in
$\psi\in[-\varphi,0]$ estimate
$$
\sum_{j=0}^m\int_0^\infty \|\exp(-i\zeta e^{i\psi}t) D_z^j
u(e^{i\psi}t+T)\|_{m-j}^2 dt \leq Const;
$$
here $D_z=-\frac i 2(\partial_{\Re z}-i\partial_{\Im
 z})$ is the complex derivative.
\end{thm}

In order to illustrate Theorem~\ref{T2} we consider an elliptic
boundary value  problem in a domain $\mathcal G\subset \Bbb R^{n}$
with quasicylindrical end. We assume that the boundary
$\partial\mathcal G$ of $\mathcal G$ is smooth, the set $\{x\in
\mathcal G:x_n<1\}$ is bounded, and the set $\{x\in\mathcal G:
x_n>1\}$ coincides with the horn-like quasicylindrical end
$\{x=(x',x_n)\in\Bbb R^{n}: x_n^{-a}x'\in \Omega, x_n>1\}$, where
$a<0$ and $\Omega$ is a domain in $\Bbb R^{n-1}$ with compact
closure and smooth boundary. We introduce the Dirichlet problem
\begin{equation}\label{example2}
\mathscr L(D_x)v(x)=f(x),\ x\in \mathcal G;\quad \partial_\nu^j
v(x)=0,\ x\in\partial \mathcal G, j=0,\dots,k-1,
\end{equation}
for  a $2k$ order elliptic differential operator $\mathscr
L(D_x)=\mathscr L(D_{x'},D_{x_n})$ with constant coefficients. The
quasicylindrical end $\{x\in\mathcal G: x_n>1\}$ can be
transformed into the half-cylinder $\mathcal
C_a=\{(y,t):y\in\Omega, t> (1-a)^{-1}\}$ by the diffeomorphism
$(y,t)=\varkappa(x',x_n)=(x_n^{-a} x',(1-a)^{-1}x_n^{1-a})$. Let
$\chi$ be a smooth cutoff function such that $\chi(x)=0$ for
$x_1<2$ and $\chi(x)=1$ for $x_1>3$.  If $v$ is a solution to the
problem~\eqref{example2} then the function $u(y,t)=(\chi
v)\circ\varkappa^{-1}(y,t)$ extended from $\mathcal C_a$ to the
 remaining part
of the cylinder $\mathcal C=\{(y,t):y\in\Omega,t\in\Bbb R\}$ by
zero satisfies the  Dirichlet problem
\begin{equation}\label{pro2}
\begin{aligned}
(\mathcal L(D_y,D_t)+\mathcal Q(y,t,D_y,D_t))u(y,t)=F(y,t),\
(y,t)\in\mathcal C,\\
\partial_\nu^j u(y,t)=0,\ (y,t)\in\partial \mathcal C,
j=0,\dots,k-1.
\end{aligned}
\end{equation}
 Here $\mathcal L(D_y,D_t)$ denotes the homogeneous part (of order $2k$) of $\mathscr
 L(D_y,D_t)$, the operator $\mathscr L(D_{x'},D_{x_n})$
written in the coordinates $(y,t)\in\mathcal C_a$ is decomposed
into
  $\mathcal L(D_y,D_t)$ and $\mathcal
Q(y,t,D_y,D_t)$; without loss of generality we can assume that
$\mathcal Q(y,t,D_y,D_t)=0$ in $\mathcal C\setminus\mathcal C_a$.
 The subsidiary problem~\eqref{pro2} can
 be represented in the form of the equation~\eqref{equation2}.
 The coefficients $Q_j(y,z,D_y)\in\mathscr B(X_j,X_0)$ of the operator
 $\mathcal Q(y,t,D_y,D_t)=\sum_{j=0}^m Q_{m-j}(y,t,D_y) D^j_t$
  are dilation analytic in the sector $\{z\in\Bbb C:|\arg(z-T)|\leq\alpha\}$,
 where $\alpha>0$ is an angle,
the spaces $X_0,\dots,X_m$ are introduced in exactly the same way
as for the
 problem~\eqref{example1}.
 The operator $\mathfrak A(\lambda)=\mathcal L(D_y,\lambda):X_m\to X_0$
  is Fredholm for all $\lambda\in\Bbb C$ and the estimate~\eqref{M cond}
  is valid for a sufficiently large $R>0$ because $\mathscr L$ is elliptic. We suppose that the right hand side $f$ of the
problem~\eqref{example2} is a smooth compactly supported function.
Then $F$ in~\eqref{pro2} is also smooth and compactly supported,
Theorem~\ref{T2} can be applied. As a result we have: if for a
fixed $\zeta\in\Bbb C$ the line $\Bbb R+\zeta$ is free from the
eigenvalues of the pencil $\Bbb C\ni\lambda\mapsto \mathfrak
A(\lambda)$ and  $v(x)=v(x',x_n)$ is a solution  to the
problem~\eqref{example2} such that
\begin{equation}\label{1es}
\sum_{j=0}^m \int_{(1-a)^{-1}}^{+\infty} \|\exp(-i\zeta t)D^j_t
v\circ\varkappa^{-1}(\cdot,t)\|_{m-j}^2\, dt <\infty
\end{equation}
 then for a sufficiently large number $T>0$
and an angle $\varphi>0$ the solution  $v$ can be continued to a
holomorphic function
\begin{equation}\label{hol}
\{z\in\Bbb C: |\arg(z-T)|\leq \varphi \}\ni z \mapsto
v\circ\varkappa^{-1}(\cdot,z)\in X_m
\end{equation}
satisfying the uniform in $\psi\in[-\varphi,\varphi]$ estimate
\begin{equation}\label{2es}
\sum_{j=0}^m \int_{(1-a)^{-1}}^{+\infty} \|\exp(-i\zeta
e^{i\psi}t)D^j_z
v\circ\varkappa^{-1}(\cdot,e^{i\psi}t+T)\bigr)\|_{m-j}^2 dt\leq C
\end{equation}
with  a constant $C$. The condition~\eqref{1es} provides the
inclusion $u\in\mathsf D^m_\zeta(\Bbb R)$
  for the solution of~\eqref{pro2}. Since the function~\eqref{hol} is holomorphic we can perform the
complex scaling of the problem~\eqref{example2}
 with respect to the coordinate $z$ in the conical neighbourhood $\{z\in\Bbb C: |\arg(z-T)|\leq \varphi\}$ of
infinity. The estimate~\eqref{2es} with the weight number
 $\zeta\in\Bbb C$ controls the behavior of
the scaled solution  at infinity. From properties of the
 holomorphic function~\eqref{hol} satisfying~\eqref{2es} it follows that
for any $\phi>0$ and $j\geq 0$ the value $|\exp(-i\zeta
z)||z|^{1/2+j}\|D_z^j v\circ\varkappa^{-1}(\cdot,z)\|_{m}$
uniformly tends to zero as $z\to \infty$,
$|\arg(z-T)|\leq\varphi-\phi$.


  Statements of elliptic
problems in domains with quasicylindrical ends in scales of
weighed Sobolev spaces and properties of solutions were studied in
many works~\cite{Feigin,BagirovFeigin,MP1,MP2,R17,ref7} and
others. In particular, from the known results it follows that if
the line $\Bbb R+\zeta$ is free from the eigenvalues of the pencil
$\Bbb C\ni\lambda\mapsto \mathfrak A(\lambda)$  and the right hand
side $f\in C_0^\infty(\overline{\mathcal G})$ of the
problem~\eqref{example2} is subjected to a finite number of
orthogonality conditions then the problem has at least one
solution (but not more than a finite number of solutions) $v\in
W^m_{loc}(\mathcal G)$  satisfying the conditions~\eqref{1es};
here $W^m_{loc}(\mathcal G)$ denotes the space of functions that
are locally in the Sobolev space $W^m(\mathcal G)$. The results on
the analytic properties of the solutions to the
problem~\eqref{example2} are new. Let us note that the assumption
$f\in C_0^\infty(\mathcal G)$ on the right hand side of the
problem~\eqref{example2} is made for simplicity only, one can
formulate a weaker assumption in terms of the Hardy-Sobolev spaces
using Theorem~\ref{---} instead of Theorem~\ref{T2}, see also
Remark~\ref{r---}. Moreover, quasicylindrical ends of many other
geometric shapes (e.g. quasicylindrical ends that approach a
cylinder at infinity) can be
 considered by choosing suitable diffeomorphisms that map quasicylinders
 onto the half-cylinder $\{(y,t):y\in\Omega, t>0\}$, c.f. \cite{MP2} or
\cite[Chapter 9]{ref7}. For other examples of applications of
differential equations with operator coefficients to boundary
value problems we refer to the books~\cite{R17,ref7}, some of
these examples can also be
 considered in the scale of  Hardy-Sobolev spaces.



\section{Spaces of analytic vector valued functions}\label{sec1}

This section deals with weighted spaces of vector valued functions
analytic in a dual complex cone or in a half-plane. In
Subsection~\ref{ss2.1} we introduce the weighted Hardy classes
$\mathscr H^\ell_w(\mathcal K^\varphi_\zeta; X)$ and study their
basic properties.  In Subsection~\ref{ss2.2} we deal with some
weighted Hardy classes in a half-plane, in particular we prove
that these classes coincide with classical Hardy classes if the
weight numbers are zeros. We also clarify relations between the
Hardy classes in cones and in half-plains. In
Subsection~\ref{ss2.3} we proceed to study the classes $\mathscr
H^\ell_w(\mathcal K^\varphi_\zeta; X)$; here we rely heavily on
the material of Subsection~\ref{ss2.2}. Finally, in
Subsection~\ref{OP} we consider Fredholm polynomial operator
pencils in spaces of analytic vector valued functions.
\subsection{Weighted Hardy classes in  cones}\label{ss2.1} Let $X$ be a Hilbert
space and
 let
$e^{i\psi}\Bbb R+\zeta=\{\lambda\in\Bbb
C:\lambda=e^{i\psi}\xi+\zeta,\xi\in\Bbb R\}$ be the line in the
complex plane $\Bbb C$,  where $\psi$ is an angle and $\zeta$ is a
fixed complex parameter.
 For $\ell\in\Bbb R$ and $w\in\Bbb C$ we
introduce the weighted $L_2$ space $\mathscr
W^\ell_w(e^{i\psi}\Bbb R+\zeta;X)$ of (classes of)
   functions
$e^{i\psi}\Bbb R+\zeta\ni\lambda\mapsto \mathcal F(\lambda)\in X$
with the finite norm
\begin{equation}\label{NW}
\|\mathcal F;\mathscr W^\ell_w (e^{i\psi}\Bbb
R+\zeta;X)\|^2=\int_{e^{i\psi}\Bbb R+\zeta}
|\exp\{2iw\lambda\}|(1+|\lambda|^2)^\ell\|\mathcal
F(\lambda)\|^2\,|d\lambda|,
\end{equation}
where $\|\cdot\|$ is the Hilbert norm in $X$. The space $\mathscr
W^\ell_w(e^{i\psi}\Bbb R+\zeta;X)$ with the norm \eqref{NW} is a
Hilbert space \cite{R8}. It is clear that the space $\mathscr
W^\ell_w(e^{i\psi}\Bbb R+\eta;X)$ and its norm does not change
while $\eta$ travels along the line $e^{i\psi}\Bbb R+\zeta$. Let
us also note that
\begin{equation}\label{chnorm}
\|\mathcal F;\mathscr W^\ell_v(e^{i\psi}\Bbb R+\zeta;X)\|=e^{\Im\{
\zeta(w-v)\}}\|\mathcal F;\mathscr W^\ell_w(e^{i\psi}\Bbb
R+\zeta;X)\|,\  v\in e^{-i\psi}\Bbb R+w.
\end{equation}
 The embedding $ \mathscr W^\ell_w(e^{i\psi}\Bbb
R+\zeta;X)\subset \mathscr W^{s}_w(e^{i\psi}\Bbb R+\zeta;X)$ is
continuous for $\ell> s$.

By $ \mathcal K^\varphi_\zeta$ we denote the open double-napped
cone
\begin{equation*}
 \mathcal
K^\varphi_\zeta=\{\lambda\in\Bbb C: \lambda=e^{i\psi}\xi+\zeta,
\xi\in\Bbb R\setminus\{0\}, 0<\psi<\varphi\}
\end{equation*}
with the vertex $\zeta$ and the angle $\varphi\in(0,\pi]$.
 If $\varphi=\pi$, then
$\mathcal K^\varphi_\zeta=\Bbb C\setminus(\Bbb R+\zeta)$.
\begin{defn}\label{d1}
We introduce the weighted Hardy class $\mathscr H_w^\ell (\mathcal
K^\varphi_\zeta;X)$ as the set of all analytic functions $\mathcal
K^\varphi_\zeta\ni\lambda\mapsto \mathcal F(\lambda)\in X$
satisfying the uniform with respect to $\psi\in(0,\varphi)$
estimate
\begin{equation}\label{no1}
\|\mathcal F;\mathscr W^\ell_w(e^{i\psi}\Bbb R+\zeta;X)\|\leq
C(\mathcal F).
\end{equation}
\end{defn}


 Before we proceed further,
 we cite the following proposition, that contains some
facts from the theory of Fourier transforms of analytic functions.
\begin{prop}\label{p1}
 {\rm 1.} Let $\Phi$ be an analytic in the strip
 $$
 \Pi=\{s\in\Bbb C: s=r+i\psi, r\in \Bbb R, 0<\psi<\varphi\}
 $$
  function taking the
values in a Hilbert space $X$. We set
$\Phi_\psi(r)\equiv\Phi(r+i\psi)$. Suppose that $\Phi$ satisfies
the estimate
\begin{equation}\label{es}
\|\Phi_\psi; L_2(\Bbb R;X)\|\leq Const, \quad\psi\in(0,\varphi);
\end{equation}
here $L_2(\Bbb R;X)\bigl(\equiv \mathscr W^0_0(\Bbb R;X)\bigr)$
denotes the space of square-summable functions $\Phi_\psi:\Bbb
R\to X$. The following assertions are true.

(i) The function $\Phi$ has  boundary limits
$\Phi_0,\Phi_\varphi\in L_2(\Bbb R;X)$ in the sense that for
almost all $r\in \Bbb R$ we have $\| \Phi(s)-\Phi_0(r) \|\to 0$ as
$s$ tends to $r$ by a non-tangential to $\Bbb R$ path, and $
\|\Phi(s)-\Phi_\varphi(r)\|{\to}  0$ as $s$ tends to $r+i\varphi$
by a non-tangential to $\Bbb R+i\varphi$ path; moreover,
$$
\|\Phi_0-\Phi_\psi; L_2(\Bbb R;X)\|\to 0,\quad \psi\to 0 +;\ \
\|\Phi_\varphi-\Phi_\psi; L_2(\Bbb R;X)\|\to 0,\quad \psi\to
\varphi -.
$$

 (ii) For all $\psi\in [0,\varphi]$ the following estimate holds
 $$
\|\Phi_\psi; L_2(\Bbb R;X)\|\leq \|\Phi_0;
 L_2(\Bbb R;X)\|+\|\Phi_\varphi; L_2(\Bbb R;X)\|.
$$

 (iii) For every compact set $\mathfrak
 K\subset \Pi$ there is an independent of $\Phi$ constant $C(\mathfrak K)$ such that
$$ \|\Phi(s)\|\leq C(\mathfrak K)(\|\Phi_0;
 L_2(\Bbb R;X)\|+\|\Phi_\varphi; L_2(\Bbb R;X)\|), \quad s\in\mathfrak K.
$$

(iv) The value $\|\Phi (s)\|$ uniformly tends to zero as
  $s$ goes to infinity in the strip $\{s\in \Pi: \Im s\in[\phi,\varphi-\phi]\}$,
  $\phi>0$.

(v) The set $\{\Phi\!\!\upharpoonright_{\Bbb R+i\psi}:\Phi\in
H(\Pi;X)\}$ is dense in the space $\mathscr W^0_0(\Bbb R+i\psi;
X)$ for any $\psi\in[0,\varphi]$.

{\rm 2.} Let $\Phi$ be an analytic function $\overline{\Pi}\ni
s\mapsto \Phi(s)\in X$ and $\Phi_\psi(r)\equiv\Phi(r+i\psi)$.
Suppose that $\Phi_0,\Phi_\varphi\in L_2(\Bbb R; X)$, and
$\|\Phi(s)\|\leq const$ for all $s\in \overline{\Pi}$. Then the
estimate \eqref{es} is valid.

{\rm 3.} A function $\Pi\ni s\mapsto \Phi(s)\in X$ is an analytic
function satisfying the estimate \eqref{es} if and only if it can
be represented in the form
\begin{equation}\label{CauchyStrip}
\Phi(s)=\frac 1 {2\pi i} \int_{-\infty}^{+\infty} \frac
{\Psi_0(r)}{r-s} \,dr-\frac 1 {2\pi
i}\int_{-\infty}^{+\infty}\frac{\Psi_\varphi(r)}{r+i\varphi-s}\,dr,\quad
s\in \Pi,
\end{equation}
for some $\Psi_0,\Psi_\varphi\in L_2(\Bbb R;X)$.
 \end{prop}
\begin{pf}  These results are well-known in the case of
analytic functions taking their values in $\Bbb C$; see e.g.
\cite{R1,R2,R3}.  The generalization for the case of functions
with values in a Hilbert space $X$ is straightforward due to the
Fourier transformation $\mathscr F:L_2(\Bbb R;X)\to L_2(\Bbb R;X)$
that implements an isometric isomorphism; see e.g. \cite{R15,
R16}.  The Parseval equality and other facts about the Fourier
transformation of vector valued functions can be found  in
\cite{R4, R15}, see also Section~\ref{FL} of this paper.
\qed\end{pf}
\begin{prop}\label{p2} {\rm 1.} Let $\mathcal F\in\mathscr H^\ell_w(\mathcal
K^\varphi_\zeta;X)$ with some $\ell\in\Bbb R$, $w,\zeta\in \Bbb
C$, and $\varphi\in(0,\pi]$. Then the following assertions are
true.

(i)   $\mathcal F$ has boundary limits $\mathcal F_0\in \mathscr
W^\ell_w(\Bbb R +\zeta;X)$ and $\mathcal F_\varphi\in\mathscr
W^\ell_w(e^{i\varphi}\Bbb R+\zeta;X)$ in the sense that for almost
all points $\mu$ of the boundary $\partial\mathcal
K^\varphi_\zeta=(\Bbb R+\zeta)\cup(e^{i\varphi}\Bbb R+\zeta)$ we
have $\|\mathcal F(\lambda)-\mathcal F_0(\mu)\|\to 0$   as
$\lambda$ non-tangentially tends to $\mu\in\Bbb R+\zeta$, and
$\|\mathcal F(\lambda)-\mathcal F_\varphi(\mu)\|\to 0$   as
$\lambda$ non-tangentially tends to $\mu\in e^{i\varphi}\Bbb
R+\zeta$. Moreover,
$$
\begin{aligned}
\|(e_w\mathcal F)\circ\varkappa_{\psi,\zeta}-(e_w\mathcal
F_0)\circ\varkappa_{0,\zeta};\mathscr W^\ell_0(\Bbb R;X)\|&\to
0,\quad\psi\to 0+,
\\\|(e_w\mathcal
F)\circ\varkappa_{\psi,\zeta}-(e_w\mathcal
F_\varphi)\circ\varkappa_{\varphi,\zeta};\mathscr W^\ell_0(\Bbb
R;X)\|&\to 0,\quad\psi\to \varphi-,
\end{aligned}
$$
where $e_w$ denotes the weight function $\lambda\mapsto\exp(i w
\lambda)$ and $\varkappa_{\psi,\zeta}:\Bbb R\to e^{i\psi}\Bbb
R+\zeta$ is the linear transformation
$\varkappa_{\psi,\zeta}(\xi)=e^{i\psi}\xi+\zeta$.

Here and elsewhere we shall suppose that every element $\mathcal
F\in\mathscr H^\ell_w(\mathcal K^\varphi_\zeta;X)$ is extended
 to the boundary $\partial\mathcal
K^\varphi_\zeta$ by the non-tangential limits. In the case
$\varphi=\pi$ we distinguish the banks $\lim_{\psi\to
0+}(e^{i\psi}\Bbb R+\zeta)$ and $\lim_{\psi\to\pi-}(e^{i\psi}\Bbb
R+\zeta)$  in $\partial\mathcal K^\pi_\zeta$.

(ii) The estimate
\begin{equation}\label{no2}
\|\mathcal F;\mathscr W^\ell_w(e^{i\psi}\Bbb R+\zeta;X)\|\leq
C\bigl( \|\mathcal F;\mathscr W^\ell_w(\Bbb
R+\zeta;X)\|+\|\mathcal F;\mathscr W^\ell_w(e^{i\varphi}\Bbb
R+\zeta;X)\|\bigr)
\end{equation}
is valid, where the constant $C$ is independent of
$\psi\in[0,\varphi]$, $\mathcal F$, and $w\in\Bbb C$.

(iii) For every compact set $\mathfrak K\in \mathcal
K^\varphi_\zeta$ there is an independent of $\mathcal F$ and
$w\in\Bbb C$ constant $c(\mathfrak K)$ such that for all
$\lambda\in\mathfrak K$ we have
\begin{equation}\label{no2`}
|\exp \{iw\lambda\}|\|\mathcal F(\lambda)\|\leq c(\mathfrak
K)(\|\mathcal F;\mathscr W^\ell_w(\Bbb R+\zeta;X)\|+\|\mathcal
F;\mathscr W^\ell_w(e^{i\varphi}\Bbb R+\zeta;X)\|).
\end{equation}

(iv) The value
$$
|\exp\{iw \lambda\}|(
1+|\lambda|)^\ell|\lambda-\zeta|^{1/2}\|\mathcal F(\lambda)\|
$$
uniformly tends to zero as $\lambda$ goes to infinity (or
$\lambda$ goes to $\zeta$) in the cone
$$
\{\lambda \in\Bbb C:\lambda=e^{i\psi}t+\zeta, t\in\Bbb
R,\psi\in[\phi,\varphi-\phi]\}, \phi>0.
$$

(v) For any $\psi\in[0,\varphi]$ the set $ \{\mathcal
F\!\!\upharpoonright_{e^{i\psi}\Bbb R+\zeta}:\mathcal F\in\mathscr
H^\ell_w(\mathcal K^\varphi_\zeta; X)\}$ is dense in the space
$\mathscr W^\ell_w(e^{i\psi}\Bbb R+\zeta;X)$.

{\bf\rm 2.} Let $\zeta\in\Bbb C$, $\varphi\in(0,\pi]$, and let
$\mathcal F$ be an analytic in $\overline{\mathcal
K^\varphi_\zeta}\setminus\{\zeta\}$ function with values in $X$.
Suppose that for some $\ell\in\Bbb R$ and $w\in\Bbb C$  the
estimate
$$
|\exp\{iw
\lambda\}|(1+|\lambda|)^\ell|\lambda-\zeta|^{1/2}\|\mathcal
F(\lambda)\| \leq Const,\quad \lambda\in \overline{\mathcal
K^\varphi_\zeta}\setminus\{\zeta\},
$$
holds and the inclusions $\mathcal F\in \mathscr W^{\ell}_w(\Bbb
R+\zeta;X)$, $\mathcal F\in \mathscr W^{\ell}_w(e^{i\varphi}\Bbb
R+\zeta;X)$ are valid. Then $\mathcal F\in\mathscr
H^\ell_w(\mathcal K^\varphi_\zeta;X)$.

\end{prop}
\begin{pf} Let us define an equivalent norm in the space $\mathscr W^\ell_w(e^{i\psi}\Bbb R+\zeta;X)$ by the equality
\begin{equation}\label{eqn}
\begin{aligned}
|\mspace{-2mu}|\mspace{-2mu}|\mathcal F;\mathscr
W^\ell_w(e^{i\psi}\Bbb R&+\zeta;X)|\mspace{-2mu}|\mspace{-2mu}|^2
=\int\limits_{e^{i\psi}\Bbb
R^-+\zeta}|\exp\{2iw\lambda\}||\lambda-\zeta-i|^{2\ell}\|\mathcal
F(\lambda)\|^2\,|d\lambda|\\&+\int\limits_{e^{i\psi}\Bbb
R^++\zeta}|\exp\{2iw\lambda\}|\lambda-\zeta+i|^{2\ell}\|\mathcal
F(\lambda)\|^2\,|d\lambda|;
\end{aligned}
\end{equation}
here $ e^{i\psi}\Bbb R^\pm +\zeta=\{\lambda\in\Bbb C:\lambda=\xi
e^{i\psi}+\zeta, \xi\gtrless0\} $ and $\varphi\in(0,\pi]$. The
norm \eqref{eqn} is equivalent to the norm \eqref{NW} because of
the inequalities
\begin{equation}
(1+|\zeta-i|^2)^{-1}(1+|\lambda|^2)\leq 3|\lambda-\zeta+ i|^2\leq
6(1+|\zeta-i|^2)(1+|\lambda|^2),  \Im(\lambda-\zeta)\geq
0,\label{triv} \end{equation}
\begin{equation*}
(1+|\zeta+i|^2)^{-1}(1+|\lambda|^2)\leq 3|\lambda-\zeta- i|^2\leq
6(1+|\zeta+i|^2)(1+|\lambda|^2), \Im(\lambda-\zeta)\leq
0.\end{equation*}
 Let $\Pi$ be the same strip as in
Proposition~\ref{p1}. For $s\in \Pi$ we set
\begin{equation}\label{FtoP}
\begin{aligned}
\Psi (s)&=\exp\{iw(\zeta- e^{s})+s/2\}(-i-e^s)^\ell \mathcal
F(\zeta-e^{s}),\\
 \Phi(s)&= \exp\{iw (\zeta+e^{s})+s/2\}(i+e^s)^\ell \mathcal
F(\zeta+e^{s});
\end{aligned}
\end{equation}
here we use  analytic in $\overline{\Pi}$ branches of the
functions $(-i-e^s)^\ell$ and $(i+e^s)^\ell$. The functions $\Phi$
and $\Psi$ are analytic in $\Pi$. For all $\psi\in (0,\varphi)$
we have
\begin{equation}\label{FtP}
\begin{aligned}
\int_{-\infty}^\infty\|\Psi(r+i\psi)\|^2\,d
r&+\int_{-\infty}^\infty\|\Phi(r+i\psi)\|^2\,d r\\&=
|\mspace{-2mu}|\mspace{-2mu}|\mathcal F;\mathscr
W^\ell_w(e^{i\psi}\Bbb
R+\zeta;X)|\mspace{-2mu}|\mspace{-2mu}|^2\leq Const(\mathcal
F).\end{aligned}
\end{equation}
 It remains to apply  the items 1,~2 of Proposition~\ref{p1} to the
functions $\Phi,\Psi$, and reformulate the results in terms of
$\mathcal F$. \qed\end{pf}

From Proposition~\ref{p2} it follows that we can take the right
hand side of \eqref{no2} as the constant in the inequality
\eqref{no1}.
\begin{prop}\label{p3`}
The class $\mathscr H^\ell_w(\mathcal K^\varphi_\zeta;X)$ endowed
with the norm
\begin{equation}\label{5}
\|\mathcal F;\mathscr H^\ell_w(\mathcal
K^\varphi_\zeta;X)\|=\|\mathcal F;\mathscr W^\ell_w(\Bbb
R+\zeta;X)\|+\|\mathcal F;\mathscr W^\ell_w(e^{i\varphi}\Bbb
R+\zeta;X)\|
\end{equation} is a Banach space.
\end{prop}
\begin{pf} With the help of \eqref{FtoP} we can split
the class $\mathscr H^\ell_w(\mathcal K^\varphi_\zeta;X)$ into two
Hardy classes $H(\Pi;X)$ in the strip $\Pi$ such that the equality
\eqref{FtP} holds. Since the Hardy class $H(\Pi;X)$ is complete
(see e.g. \cite{R15,R16}), the class $\mathscr H^\ell_w(\mathcal
K^\varphi_\zeta;X)$ is also complete. \qed\end{pf}

The following proposition contains some elementary properties of
the introduced classes $\mathscr H^\ell_w(\mathcal
K^\varphi_\zeta;X)$, it is given without proof.
\begin{prop}\label{p1.5} Let $w,\zeta\in\Bbb C$, $\ell\in\Bbb
R$, and let $\varphi\in(0,\pi]$.

(i)If $\phi\in(0,\varphi]$ and $s\leq\ell$ then the embedding
$\mathscr H^\ell_w(\mathcal K^\varphi_\zeta;X)\subseteq\mathscr
H^s_w(\mathcal K^\phi_\zeta;X)$ is continuous.

(ii) Let $\mathcal F\in\mathscr H^\ell_w(\mathcal
K^\varphi_\zeta;X)$ and $\mathcal K^{\varphi,\pm}_\zeta
=\{\lambda\in\mathcal
K^\varphi_\zeta:\Im\lambda\gtrless\Im\zeta\}$. We set $\mathcal
F^+ =\mathcal F$ on $\lambda\in\mathcal K^{\varphi,+}_\zeta$ and
$\mathcal F^+=0$ on $\mathcal K^{\varphi,-}_\zeta$; let also
$\mathcal F^-=\mathcal F-\mathcal F^+$. The mappings
\begin{equation*}
\mathscr H^\ell_w(\mathcal K^\varphi_\zeta;X)\ni\mathcal F\mapsto
 \mathcal F^+\in \mathscr H^\ell_u(\mathcal
K^\varphi_\zeta;X),
\end{equation*}
\begin{equation*}
\mathscr H^\ell_w(\mathcal K^\varphi_\zeta;X)\ni\mathcal F\mapsto
 \mathcal F^-\in \mathscr H^\ell_v(\mathcal
K^\varphi_\zeta;X)
\end{equation*}
are continuous if
 $\Im(e^{i\psi}u)\geq \Im(e^{i\psi}w)$ and
 $\Im(e^{i\psi}v)\leq \Im(e^{i\psi}w)$ for all $\psi\in[0,\varphi]$.

(iii) Suppose that an analytic function $\overline{\mathcal
K^\varphi_\zeta}\ni\lambda\mapsto p(\lambda)\in\Bbb C$ satisfies
the estimate $| p(\lambda)|\leq C(1+|\lambda|)^{s}$ for all
$\lambda\in \overline{\mathcal K^\varphi_\zeta}$ and some
$s\in\Bbb R$. Then the norm of the multiplication operator
$$ \mathscr H^\ell_w(\mathcal K^\varphi_\zeta;X)\ni \mathcal
F\mapsto  p\,\mathcal F\in \mathscr H^{\ell-s}_w(\mathcal
K^\varphi_\zeta;X)
$$
 is bounded uniformly in $w\in\Bbb C$.

(iv) The following mapping~\eqref{m} is an isomorphism
\begin{equation}\label{m}
\begin{aligned}
\mathscr H^\ell_w(\mathcal K^\varphi_\zeta;X)\ni\mathcal F\mapsto
&(\cdot-\zeta+i)^{s}\mathcal
F^+(\cdot)\\&+(\cdot-\zeta-i)^s\mathcal F^-(\cdot)\in \mathscr
H^{\ell-s}_w(\mathcal K^\varphi_\zeta;X),\ s\in\Bbb R;
\end{aligned}
\end{equation}
 here $\mathcal F^+$ and $\mathcal F^-$ are the
same as in (ii), we use an analytic in $\overline{\mathcal
K^{\varphi,+}_\zeta}$ branch of the function $
(\cdot-\zeta+i)^{s}$ and an analytic in $\overline{\mathcal
K^{\varphi,-}_\zeta}$ branch of $(\cdot-\zeta-i)^s$. The norm of
the mapping~\eqref{m} and the norm of its inverse are uniformly
bounded in $w\in\Bbb C$.
\end{prop}

\begin{prop}\label{p2.6}
(i) Let $\partial \mathcal K^{\varphi,\pm}_\zeta$ denote the
boundary of $\mathcal K^{\varphi,\pm}_\zeta =\{\lambda\in\mathcal
K^\varphi_\zeta:\Im\lambda\gtrless\Im\zeta\}$. For a function
$\mathcal F\in \mathscr H^\ell_w(\mathcal K^\varphi_\zeta;X)$  the
representations
\begin{equation}\label{cau}
 \mathcal F(\lambda)=\int_{\partial \mathcal K^{\varphi,+}_\zeta}\frac
{e^{iw(\mu-\lambda)}(\mu-\eta)^s \mathcal F(\mu)}{2\pi i
(\lambda-\eta)^s(\mu-\lambda)}\, d\mu,\quad \lambda\in \mathcal
K^{\varphi,+}_\zeta,\ \eta\in {\mathcal K^{\varphi,-}_\zeta},\
s\leq\ell,
\end{equation}
\begin{equation}\label{1cau}
 \mathcal F(\lambda)=\int_{
\partial \mathcal K^{\varphi,-}_\zeta}\frac
{e^{iw(\mu-\lambda)}(\mu-\tau)^s \mathcal F(\mu)}{2\pi i
(\lambda-\tau)^s(\mu-\lambda)}\, d\mu,\quad \lambda\in \mathcal
K^{\varphi,-}_\zeta,\ \tau\in{\mathcal K^{\varphi,+}_\zeta},\
s\leq\ell,
\end{equation}
are valid, where we use an analytic in $\overline{\mathcal
K^{\varphi,+}_\zeta}$ branch of the function $(\cdot-\eta)^s$ and
an analytic in $\overline{\mathcal K^{\varphi,-}_\zeta}$ branch of
$(\cdot-\tau)^s$, the contours of integration are oriented such
that  $\lambda$ lies on the left side while $\mu$ travels along a
contour. The integrals
 are absolutely convergent in the space
$X$. The formulas \eqref{cau}, \eqref{1cau} recover a function
$\mathcal F\in \mathscr H^\ell_w(\mathcal K^\varphi_\zeta;X)$ from
its boundary limits.

 (ii) Let $\mathcal F\in \mathscr H^\ell_w(\mathcal K^\varphi_\zeta;X)$. The value $
|\exp\{iw\lambda\}||\lambda|^{\ell} \|\mathcal F(\lambda)\| $
uniformly tends to zero as $\lambda$ goes to infinity in the set
$\{\lambda\in \mathcal K^\varphi_\zeta:\dist\{\lambda,\partial
\mathcal K^\varphi_\zeta\}\geq \epsilon\}$, where $\epsilon>0$.

\end{prop}

The proof of Proposition~\ref{p2.6} is displayed in Appendix, see
page~\pageref{proofpage}. The following proposition shows that in
order to recover a function $\mathcal F\in \mathscr
H^\ell_w(\mathcal K^\varphi_\zeta;X)$ it is not necessary to know
the boundary limits of $\mathcal F$ on the whole boundary
$\partial \mathcal K^\varphi_\zeta$.

\begin{prop}\label{lost} A function $\mathcal F\in \mathscr H^\ell_w(\mathcal
K^\varphi_\zeta;X)$ can be uniquely recovered from its boundary
limits $\mathcal F\!\!\upharpoonright_{\Bbb R+\zeta}$. Similarly,
in order to recover $\mathcal F$ it suffices to know the
non-tangential boundary limits of $\mathcal F$ on the part
$e^{i\varphi}\Bbb R+\zeta$ of the boundary $\partial \mathcal
K^\varphi_\zeta$ (or, equivalently, on the part $e^{i\varphi}\Bbb
R^++\zeta\cup\Bbb R^-+\zeta\subset\partial \mathcal
K^\varphi_\zeta$, or on the part $e^{i\varphi}\Bbb
R^-+\zeta\cup\Bbb R^++\zeta\subset\partial \mathcal
K^\varphi_\zeta$).
\end{prop}
\begin{pf} Here again we split $\mathcal F\in \mathscr H^\ell_w(\mathcal
K^\varphi_\zeta;X)$ into the functions $\Psi$ and $\Phi$ in the
Hardy class  $H(\Pi;X)$, see \eqref{FtoP} and \eqref{FtP}. The
non-tangential boundary limits of the functions $\Psi$,  $\Phi$,
and $\mathcal F$ are related by the formulas~\eqref{FtoP}. As is
well-known (see e.g.~\cite{R1} for the scalar case, and \cite{R15,
R16} for the case of $X$-valued functions), for all $\psi\in
(0,\varphi)$ and $r\in \Bbb R$ we have
\begin{equation}\label{z2}
\Phi(r+i\psi)=\mathscr F^{-1}_{t\to r} e^{t\psi}\mathscr F_{r\to t
}\Phi_0(r),\ \Phi(r+i\psi)=\mathscr F^{-1}_{t\to r}
e^{t(\psi-\varphi)}\mathscr F_{r\to t }\Phi_\varphi(r);
\end{equation}
 here
$\mathscr F:L_2(\Bbb R;X)\to L_2(\Bbb R;X)$ is the Fourier
transformation and $\mathscr F^{-1}$ is its inverse, as in
Proposition~\ref{p1} the boundary limits of $\Phi$ are denoted by
$\Phi_0$ and $\Phi_\varphi$. Obviously, in the same way $\Psi$ can
be recovered from $\Psi_0$ or $\Psi_\varphi$. Knowing $\Phi$ and
$\Psi$ we get $\mathcal F$ from the formulas~\eqref{FtoP}.
\qed\end{pf}

\subsection{Weighted Hardy classes in a half-plane}\label{ss2.2}
 Let $\Bbb C^+$ and $\Bbb C^-$ be the upper and the lower half-plane,
 $\Bbb C^\pm=\{\lambda\in\Bbb C:\Im\lambda\gtrless0\}$.
 We shall use the notations
 $$e^{i\varphi}\Bbb R^\pm
+\zeta=\{\lambda\in\Bbb C:\lambda=\xi e^{i\varphi}+\zeta,
\pm\xi>0\},$$ $$e^{i\varphi}\Bbb C^\pm+\zeta=\{\lambda\in\Bbb
C:\lambda=e^{i\varphi}\mu+\zeta, \mu\in\Bbb C^\pm\},$$ where
$\zeta\in\Bbb C$ and $\varphi$ is an angle. By $\mathscr
W_w^\ell(e^{i\varphi}\Bbb R^\pm +\zeta;X)$ we denote by  the
Hilbert space of all (classes of) functions $ e^{i\varphi}\Bbb
R^\pm +\zeta\ni\lambda\mapsto\mathcal F(\lambda)\in X$  with the
finite norm
\begin{equation*}
\|\mathcal F;\mathscr W_w^\ell(e^{i\varphi}\Bbb R^\pm
+\zeta;X)\|^2=\int\limits_{e^{i\varphi}\Bbb R^\pm +\zeta}|\exp\{
2i w\lambda\}||(1+|\lambda|^2)^\ell\|\mathcal
F(\lambda)\|^2\,|d\lambda|.
\end{equation*}
\begin{defn}\label{defHalfPlane} Let $\zeta,w\in\Bbb C$
and $\ell\in\Bbb R$. We introduce the Hardy class $\mathscr
H^\ell_w(e^{i\varphi}\Bbb C^\pm+\zeta;X)$  as the set of all
analytic functions $e^{i\varphi}\Bbb C^\pm+\zeta\ni\lambda\mapsto
\mathcal F(\lambda)\in X$ satisfying the uniform in
$\psi\in(\varphi,\varphi+\pi)$ estimate $$\|\mathcal F;\mathscr
W^\ell_w(e^{i\psi}\Bbb R^\pm +\zeta;X)\|\leq C(\mathcal F).$$
\end{defn}
It is clear that $\mathscr H^\ell_w(e^{i\varphi}\Bbb
C^++\zeta;X)\equiv \mathscr H^\ell_w(e^{i(\varphi+\pi)}\Bbb
C^-+\zeta;X)$. From Definition~\ref{defHalfPlane} and
Definition~\ref{d1} it is easily seen that a function $\mathcal F$
is in the Hardy class $\mathscr H^\ell_w(\mathcal K^\pi_\zeta;X)$
if and only if $\mathcal F\in \mathscr H^\ell_w(\Bbb
C^++\zeta;X)\cap \mathscr H^\ell_w(\Bbb C^-+\zeta;X)$. Let us also
note that if $\mathcal F\in \mathscr H^\ell_w(\Bbb C^++\zeta;X)$
then after extension of the function $\mathcal F$ to the
half-plane $\Bbb C^-+\zeta$ by zero we have the inclusion
$\mathcal F\in \mathscr H^\ell_w(\mathcal K^\varphi_\zeta;X)$ for
all $\varphi\in(0,\pi]$; moreover, a function $\mathcal F$ is in
the Hardy class $\mathscr H^\ell_w(\Bbb C^++\zeta;X)$ if and only
if the function $\mathcal F$ extended to $\Bbb C^-+\zeta$ by zero
is in the class $ \mathscr H^\ell_w(\mathcal K^\pi_\zeta;X)$. 
The next proposition contains some direct consequences of
Proposition~\ref{p2} and Proposition~\ref{p2.6}, it is presented
without proof.

\begin{prop}\label{pr} Let $\zeta,w\in\Bbb C$,
$\ell\in\Bbb R$, and let $\varphi$ be an angle. For any function
$\mathcal F\in\mathscr H^\ell_w( e^{i\varphi}\Bbb C^++\zeta;X)$
the following assertions are fulfilled.

(i)  $\mathcal F$ has  boundary limit $\mathcal F_\varphi\in
\mathscr W^\ell_w(e^{i\varphi}\Bbb R+\zeta;X)$ in the sense that
for almost all points $\mu$ of the boundary $e^{i\varphi}\Bbb
R+\zeta$ we have $\|\mathcal F(\lambda)-\mathcal
F_\varphi(\mu)\|\to 0$ as $\lambda$  tends to $\mu$ by a
non-tangential to $e^{i\varphi}\Bbb R+\zeta$ path. Moreover
$$\begin{aligned}
\|(e_w\mathcal F)\circ\varkappa_{\psi,\zeta}-(e_w\mathcal
F_\varphi)\circ\varkappa_{\varphi,\zeta};\mathscr W^\ell_0(\Bbb
R^+;X) \|&\to 0,\ \psi\to \varphi+,
\\\|(e_w\mathcal
F)\circ\varkappa_{\psi,\zeta}-(e_w\mathcal
F_\varphi)\circ\varkappa_{\varphi+\pi,\zeta};\mathscr
W^\ell_0(\Bbb R^-;X)\|&\to 0,\ \psi\to
(\varphi+\pi)-,\end{aligned}
$$
where $e_w$ denotes the weight function $\lambda\mapsto\exp(i w
\lambda)$ and $\varkappa_{\psi,\zeta}:\Bbb R\to e^{i\psi}\Bbb
R+\zeta$ is the linear transformation
$\varkappa_{\psi,\zeta}(\xi)=e^{i\psi}\xi+\zeta$.

 Here
and elsewhere we shall suppose that every element $\mathcal
F\in\mathscr H^\ell_w(e^{i\varphi}\Bbb C^++\zeta;X)$ is extended
 to the boundary
$e^{i\varphi}\Bbb R+\zeta$ by its non-tangential limits.

(ii) The estimate
$$
\|\mathcal F;\mathscr W^\ell_w(e^{i\psi}\Bbb R^++\zeta;X)\|\leq
 C\|\mathcal F; \mathscr W^\ell_w(e^{i\varphi}\Bbb R+\zeta;X)\|
$$
holds, where the constant $C$ is independent of
$\psi\in[\varphi,\varphi+\pi]$, $\mathcal F$, and  $w\in\Bbb C$.

 (iii) The value
$|\exp\{iw\lambda\}||\lambda|^{\ell+1/2}\|\mathcal F(\lambda)\|$
uniformly tends to zero as $\lambda$ goes to infinity in the set
$\{\lambda\in\Bbb
C:\arg(\lambda-\zeta)\in[\varphi+\phi,\varphi+\pi-\phi]\}$,
$\phi>0$.

(iv) The value $ |\exp\{iw\lambda\}||\lambda|^{\ell}\|\mathcal
F(\lambda)\|$ uniformly tends to zero as $\lambda$ goes to
infinity in the half-plane  $e^{i\varphi}\Bbb C^++\eta$, where
$\eta\in e^{i\varphi} \Bbb C^++\zeta$.

(v) The representation
\begin{equation}\label{rec}
\mathcal F(\lambda)=\int\limits_{e^{i\varphi}\Bbb R+\zeta} \frac
{e^{iw(\mu-\lambda)}(\mu-\eta)^s\mathcal F(\mu)}{2\pi i
(\lambda-\eta)^s(\mu-\lambda)}\,d\mu, \ \lambda\in
e^{i\varphi}\Bbb C^++\zeta,\ \eta\in e^{i\varphi}{\Bbb
C^-}+\zeta,\ s\leq\ell,
\end{equation}
holds, where we use an analytic in $e^{i\varphi}\overline{\Bbb
C^+}+\zeta$ branch of the function $(\cdot-\eta)^s$,  the contour
of integration is oriented such that  $\lambda$ lies on the left
side while $\mu$ travels along the contour.
 The representation \eqref{rec} recovers a function $\mathcal F\in \mathscr H^\ell_w(e^{i\varphi}{\Bbb
C^+}+\zeta;X)$
  from its non-tangential boundary limits. The integral is absolutely
convergent in $X$.
\end{prop}
As a consequence of Proposition~\ref{p3`}  we get
\begin{prop} The class $\mathscr H^\ell_w(
e^{i\varphi}{\Bbb C^+}+\zeta;X)$ endowed  with  the norm
\begin{equation}\label{NH}
\|\mathcal F;\mathscr H^\ell_w(e^{i\varphi}{\Bbb
C^+}+\zeta;X)\|=\|\mathcal F; \mathscr W^\ell_w(e^{i\varphi}\Bbb
R+\zeta;X)\|
\end{equation}
is a Banach space.
\end{prop}

The next proposition in particular establishes the equivalence of
Definition~\ref{defHalfPlane} and a universally accepted
definition of the Hardy class in a half-plane; this fact is known
for the Hardy class without weight~\cite{Van Winter I}.

\begin{prop}\label{HH} (i) An analytic function $e^{i\varphi}{\Bbb C^+}+\zeta\ni\lambda\mapsto \mathcal
F(\lambda)\in X$ is in the Hardy class ${\mathscr
H}^\ell_w(e^{i\varphi}{\Bbb C^+}+\zeta;X)$ if and only if it
satisfies the uniform  estimate
\begin{equation}\label{**}
\|\mathcal F;\mathscr W^\ell_w(e^{i\varphi}\Bbb R+\eta;X)\|\leq
C(\mathcal F),\quad\eta \in e^{i\varphi}{\Bbb C^+}+\zeta.
\end{equation}

 (ii) If $\mathcal
F\in\mathscr H^\ell_w(e^{i\varphi}{\Bbb C^+}+\zeta;X)$ then
$$
\|\mathcal F(\cdot+\eta)-\mathcal F(\cdot);\mathscr
W^\ell_w(e^{i\varphi}\Bbb R+\zeta;X)\|\to 0, \quad \eta\to 0,\
\eta\in e^{i\varphi}\Bbb C^+.
$$

 (iii) Let $\mathcal J$ be a function from the space $\mathscr W^\ell_w(e^{i\varphi}\Bbb R+\zeta;X)$. In the right hand side of
the equality~\eqref{rec} we replace $\mathcal F$ by $\mathcal J$.
Then the equality~\eqref{rec} defines a function $\mathcal
F\in{\mathscr H}^k_w(e^{i\varphi}\Bbb C^++\zeta;X)$, where
$k<s+1/2$  if $\ell-s>1/2$, and $k\leq s$ if $\ell-s\leq 1/2$.
(Certainly this does not necessarily mean that the prescribed on
$e^{i\varphi}\Bbb R+\zeta$ function $\mathcal J$ coincides with
the
 boundary limit of the analytic in $e^{i\varphi}\Bbb C^++\zeta$ function $\mathcal
 F$.) The estimate
\begin{equation}\label{est1}
\|\mathcal F;{\mathscr H}^k_w(e^{i\varphi}\Bbb C^++\zeta;X)\|\leq
C \|\mathcal J;\mathscr W^\ell_w(e^{i\varphi}\Bbb R+\zeta;X)\|
\end{equation}
 holds with an independent of   $w\in\Bbb C$ and $\mathcal J\in \mathscr W^\ell_w(e^{i\varphi}\Bbb R+\zeta;X)$ constant
$C$.

\end{prop}

\begin{pf}
 Without loss of generality we can suppose that the
parameters  $\zeta$,$w$,$\ell$, and the angle $\varphi$ are equal
to zero; if it is not the case then we identify $\mathcal
F\in\mathscr H^\ell_w(e^{i\varphi}\Bbb C^++\zeta;X)$ with
 $\mathcal G\in \mathscr H^0_0(\Bbb C^+;X)$ by the rule
 \begin{equation}\label{ident}
 \mathcal G(\lambda)=\exp\{iw e^{i\varphi}\lambda\}
 (\lambda+i)^{\ell}\mathcal F(e^{i\varphi}\lambda+\zeta),\quad\lambda\in\Bbb
 C^+,
 \end{equation}
 where we use an analytic in $\Bbb C^+$ branch of the function
 $(\cdot+i)^\ell$; cf. \eqref{triv}.

  As it was already mentioned, the assertion (i) for the Hardy class $\mathscr H^0_0(\Bbb
  C^+;X)$ is proved in~\cite{Van Winter I}. Since this result is not widely known we
   give an independent proof.

  Due to Proposition~\ref{pr}, (v) we can recover a
function $\mathcal F\in \mathscr H^0_0(\Bbb C^+;X)$ from its
non-tangential boundary limits  by the Cauchy integral taken along
the real axis. As is well-known~\cite{R15,R16}, this immediately
leads to the estimate
\begin{equation}\label{re2}
\|\mathcal F(\cdot+\eta); L_2(\Bbb R;X)\|\leq C(\mathcal F),\
\eta\in\Bbb C^+,
\end{equation}
that is the same as \eqref{**} in the case $w=\zeta=0$, $\ell=0$,
and $\varphi=0$. The necessity of the condition~\eqref{**}  is
established.

 Let $\mathcal F\in\mathscr H^0_0(\Bbb C^+;X)$ and let $\Pi$ denote the strip
$$\Pi=\{s\in\Bbb C: s=r+i\psi, r\in\Bbb R, 0<\psi<\pi\}.$$ We define
the analytic function $\Pi\ni s\mapsto \Phi(s)\in X$ by the
formula $ \Phi(s)=\mathcal F(e^s) e^{s/2}$. It is easy to see that
$$
\int_{-\infty}^{+\infty}\|\Phi(r+i\psi)\|^2\,dr=\|\mathcal
F;\mathscr W_0^0(e^{i\psi}\Bbb R^+;X)\|^2\leq Const, \quad\psi\in
[0,\pi].
$$
From the assertion $3$ of Proposition~\ref{p1} we conclude that
$\mathcal F\in\mathscr H^0_0(\Bbb C^+;X)$ if and only if the
 function $\Phi$ can be represented in the form \eqref{CauchyStrip} for some $\Psi_0,\Psi_\varphi\in L_2(\Bbb R;X)$.
  We express the representation \eqref{CauchyStrip} in terms of $\mathcal F$ and  get the
equality
\begin{equation}\label{Cau}
\sqrt{\lambda}\mathcal F(\lambda)=\frac 1{2\pi
i}\int\limits_{-\infty}^{+\infty} \frac {\mathcal G(\xi)}{(\log
\xi-\log\lambda)\sqrt{\xi}}\,d\xi,\quad\lambda\in\Bbb C^+,
\end{equation}
 with
some $\mathcal G\in  L_2(\Bbb R;X)$; here we use  analytic in
$\overline{\Bbb C^+}\setminus\{0\}$ principal branches of
logarithm and square root. Thus we can say that $\mathcal F\in
\mathscr H^0_0(\Bbb C^+;X)$ if and only if the equality
\eqref{Cau} is valid for some $\mathcal G\in L_2(\Bbb R;X)$. In
order to prove the sufficiency of the condition~\eqref{**}, we
shall show that for any analytic function $\Bbb
C^+\ni\lambda\mapsto \mathcal F(\lambda)\in X$ satisfying the
uniform estimate \eqref{re2}
 there exists $\mathcal G\in
L_2(\Bbb R;X)$ such that the equality \eqref{Cau} is fulfilled.
This part of the proof is arranged as Lemma~\ref{ApLemma2} in
Appendix.

Now we prove the assertion (ii). One of the universally accepted
definitions of Hardy classes in $\Bbb C^+$ reads: the Hardy class
$H(\Bbb C^+;X)$ is the set of all analytic functions  $\Bbb
C^+\ni\lambda\mapsto \mathcal F(\lambda)\in X$ satisfying the
uniform in $\eta$ estimate \eqref{re2}. Therefore the
rule~\eqref{ident} allows us to identify the classes $H(\Bbb
C^+;X)$ and $\mathscr H^\ell_w(e^{i\varphi}\Bbb C^++\zeta;X)$. If
$w=\zeta=0$, $\ell=0$, and $\varphi=0$ then the classes are
coincident and the assertion (ii) is known. For nonzero parameters
it is easily seen from~\eqref{ident}.

Let us prove the assertion (iii). If the parameters $\ell$, $s$,
$w$,  are zeros then $\mathcal F$ is the Cauchy integral of
$\mathcal J$, it is known  that $\mathcal F\in \mathscr
H^0_0(e^{i\varphi}\Bbb C^++\zeta;X)$ and
\begin{equation}\label{90}
\|\mathcal F;\mathscr H^0_0(e^{i\varphi}\Bbb C^++\zeta;X)\|\leq
c\|\mathcal J; \mathscr W^0_0(e^{i\varphi}\Bbb R+\zeta;X)\|;
\end{equation}
see e.g. \cite{R16}. If the parameters are not zeros then
$\exp\{iw \cdot\}(\cdot-\eta)^s \mathcal F(\cdot)$ is the Cauchy
integral of $\exp\{iw \cdot\}(\cdot-\eta)^s \mathcal J(\cdot)\in
\mathscr W^0_0(e^{i\varphi}\Bbb R+\zeta;X)$. We perform the
corresponding changes of the functions in~\eqref{90}, after
obvious estimates we arrive at the inequality
\begin{equation*}
\|\mathcal F;\mathscr H^s_w(e^{i\varphi}\Bbb C^++\zeta;X)\|\leq
C\|\mathcal J; \mathscr W^s_w(e^{i\varphi}\Bbb R+\zeta;X)\|.
\end{equation*}
This immediately leads to the estimate \eqref{est1}, where $k\leq
s\leq \ell$; see Proposition~\ref{p1.5}, (i). At this stage we can
only guarantee the validity of the estimate \eqref{est1} with
$k\leq s$. This completely proves the assertion~(iii) for the case
$\ell-1/2\leq s\leq\ell$. In the case $s<\ell-1/2$ the estimate
\eqref{est1} is to be improved so that $k\in (s, s+1/2)$. We shall
 do this in Section~\ref{sec}, see
Corollary~\ref{credit1} and Remark~\ref{creditproof}. Here we take
this fact for granted. \qed\end{pf}

\begin{cor}\label{emb} Let $\zeta,w\in\Bbb C$, $\ell\in\Bbb R$, and let
$\varphi$ be an angle.

(i)  The space $\mathscr H^\ell_w(e^{i\varphi}\Bbb C^++\eta;X)$
and its norm do not change while $\eta$ travels along the line
$e^{i\varphi}{\Bbb R}+\zeta$.

(ii) For all $\eta\in e^{i\varphi}{\Bbb C^+}+\zeta$  the embedding
$ \mathscr H^\ell_w(e^{i\varphi}\Bbb C^++\zeta;X)\subset \mathscr
H^\ell_w(e^{i\varphi}\Bbb C^++\eta;X)$ and the estimate
\begin{equation*}
\|\mathcal F;\mathscr H^\ell_w(e^{i\varphi}\Bbb C^++\eta;X)\|\leq
C \|\mathcal F; \mathscr H^\ell_w(e^{i\varphi}\Bbb C^++\zeta;X)\|
\end{equation*}
are valid. Here the constant $C$ is independent of $w$, $\eta$,
and $\mathcal F\in\mathscr H^\ell_w(e^{i\varphi}\Bbb
C^++\zeta;X)$.

(iii) Every function $\mathcal F\in\mathscr
W^\ell_w(e^{i\varphi}\Bbb R+\zeta;X)$ can be represented as the
sum $(\mathcal F^++\mathcal
F^-)\!\!\upharpoonright_{e^{i\varphi}\Bbb R+\zeta}$ of the
boundary limits of functions $\mathcal F^\pm\in \mathscr
H^s_w(e^{i\varphi}\Bbb C^\pm+\zeta;X)$, where
\begin{equation}\label{s}
\begin{aligned}
&s<1/2  &\text{ if }\ \  \ell>1/2;\phantom{\ \ \text{ and }\ \
[\ell]\leq\ell\leq[\ell]+1/2}\\
&s\leq[\ell]&\text{ if }\ \  \ell\leq 1/2\ \ \text{ and }\ \
[\ell]\leq\ell\leq[\ell]+1/2;\\
&s<[\ell]+1/2 &\text{ if }\ \  \ell\leq 1/2\ \  \text{ and }\ \
[\ell]+1/2<\ell<[\ell].
\end{aligned}
\end{equation}
In the case $s\geq -1/2$  the representation is unique. If
$s<-1/2$ then the entire  functions $$\Bbb
C\ni\lambda\mapsto\lambda^j e^{-i\lambda w}f \in
 X,\quad j=0,1,\dots,-[s+3/2],$$
 are
 in  the spaces $\mathscr H^s_w(e^{i\varphi}\Bbb
 C^\pm+\zeta;X)$ for any $f\in X$; if  $\mathcal F=(\widetilde{\mathcal F}^++\widetilde{\mathcal
F}^-)\!\!\upharpoonright_{e^{i\varphi}\Bbb R+\zeta}$ with some
functions $\widetilde{\mathcal F}^\pm\in \mathscr
H^s_w(e^{i\varphi}\Bbb C^\pm+\zeta;X)$ then for
 coefficients  $f_j\in X$ we have
 \begin{equation}\label{repr}
 \widetilde{\mathcal
F}^\pm (\lambda)=\mathcal
 F^\pm(\lambda)\pm\sum_{j=0}^{-[s+3/2]} f_j \lambda^j
e^{-i\lambda w},\quad\lambda\in e^{i\varphi}\Bbb C^\pm+\zeta.
\end{equation}
\end{cor}
\begin{pf}
The assertions (i) and (ii) are readily apparent from the
equalities~\eqref{NW}, \eqref{NH} together with
Proposition~\ref{HH}, (i), (ii). Let us demonstrate (iii). Without
loss of generality we can suppose that $\ell<1$. We define
$\mathcal G\in\mathscr W^{\ell-[\ell]}_0(e^{i\varphi}\Bbb
R+\zeta;X)$ by the formula
 $$
 \mathcal
G(\lambda)=\exp(iw\lambda)(\lambda-\eta)^{[\ell]}\mathcal
F(\lambda),\quad\lambda\in e^{i\varphi}\Bbb R+\zeta,\ \eta\notin
e^{i\varphi}\Bbb R+\zeta.
$$
  As is well-known, any function $\mathcal G\in\mathscr
W^0_0(e^{i\varphi}\Bbb R+\zeta;X)$ can be  uniquely represented in
the form
 $\mathcal G=(\mathcal G^++\mathcal G^-)\!\!\upharpoonright_{e^{i\varphi}\Bbb
R+\zeta}$, where $\mathcal G^\pm\in\mathscr H^0_0(e^{i\varphi}\Bbb
C^\pm+\zeta;X)$, the functions $\mathcal G^\pm$ are defined as the
Cauchy integral of $\mathcal G$; see e.g. \cite{R16}. In the case
$\ell-[\ell]>1/2$ the inclusion $\mathcal G\in\mathscr
W^{\ell-[\ell]}_0(e^{i\varphi}\Bbb R+\zeta;X)$ together with
Proposition~\ref{HH}, (iii)  allows us to see that $\mathcal
G^\pm\in\mathscr H^k_0(e^{i\varphi}\Bbb C^\pm+\zeta;X)$ for any
$k<1/2$. It remains to set
$$
 \mathcal F^\pm(\lambda)=\exp(-iw\lambda)(\lambda-\eta)^{-[\ell]}
 {\mathcal G}^\pm(\lambda),\quad\lambda\in e^{i\varphi}\Bbb C^\pm+\zeta.
$$
Then the inclusions $\mathcal F^\pm \in \mathscr
H^{s}_w(e^{i\varphi}\Bbb C^\pm+\zeta;X)$ hold, where $s$ is the
same as in~\eqref{s}.
 It is easy to see the validity of the remark
 concerning the
 uniqueness of the representation $\mathcal F=(\mathcal F^++\mathcal
F^-)\!\!\upharpoonright_{e^{i\varphi}\Bbb R+\zeta}$, we do not
cite the proof here. \qed\end{pf}

For $s\in\Bbb R$ and  $v\in  \Bbb C$ we introduce the operator
$\mathscr P^s_{\eta,v}(\phi,\zeta)$ acting by the rule
\begin{equation}\label{proj11}
\bigl(\mathscr P^s_{\eta,v} \mathcal
G\bigr)(\lambda)=\int\limits_{e^{i\phi}\Bbb R+\zeta}\frac
{e^{iv(\mu-\lambda)}(\mu-\eta)^s\mathcal G(\mu)}{2\pi i
(\lambda-\eta)^s(\mu-\lambda)}\,d\mu, \quad \lambda\in
e^{i\phi}\Bbb C^-+\zeta,\eta\in e^{i\phi}{\Bbb C^+}+\zeta,
\end{equation}
where we use an analytic in $ e^{i\phi}\overline{\Bbb C^-}+\zeta$
branch of the function $(\cdot-\eta)^s$; here and elsewhere we
 shall omit the parameters $\phi$ and
$\zeta$ in the notations of the operators $\mathscr
P^s_{\eta,v}(\phi,\zeta)$ when it can be done without ambiguity.

Let $s\leq \ell$, $k=s$ if $\ell-s\leq 1/2$, and $k\in[s,s+1/2)$
if $\ell-s>1/2$. From the equality~\eqref{chnorm} and
Proposition~\ref{HH}, (iii) we see that the operator
\begin{equation}\label{proj}
\mathscr P^s_{\eta,v}: \mathscr W^\ell_w(e^{i\phi}\Bbb R+\zeta;
X)\to \mathscr H^k_v(e^{i\phi}\Bbb C^-+\zeta; X),\quad v\in
e^{-i\phi} \Bbb R+w,
\end{equation}
is continuous; moreover,  the estimate
\begin{equation}\label{1proj}
\|\mathscr P^s_{\eta,v} \mathcal G; \mathscr H^k_v(e^{i\phi}\Bbb
C^-+\zeta; X)\|\leq  Ce^{\Im\{\zeta (w-v)\}}   \|\mathcal
G;\mathscr W^\ell_w(e^{i\phi}\Bbb R+\zeta; X)\|
\end{equation}
holds with an independent of $w\in\Bbb C$ and $v\in e^{-i\phi}
\Bbb R+w$ constant $C$. Proposition~\ref{pr}, (v) allows us to
identify the elements of the space $\mathscr H^k_v(e^{i\phi}\Bbb
C^-+\zeta; X)$ with their boundary limits in $\mathscr
W^k_w(e^{i\phi}\Bbb R+\zeta; X)$. Then $\mathscr
P^r_{\eta,v}\mathcal F=\mathcal F$ for any $\mathcal F\in \mathscr
H^k_v(e^{i\phi}\Bbb C^-+\zeta; X)$ and $r\leq k$;
cf.~\eqref{proj11},~\eqref{rec}. Therefore~\eqref{proj} is a
projection operator with the property
 \begin{equation}\label{2proj}
 \mathscr
P^r_{\eta,w}\mathscr P^s_{\eta,v}=\mathscr P^s_{\eta,v},\quad
r\leq s,\ v\in e^{-i\phi}\overline{\Bbb R^+}+w.
\end{equation}
 It is easy to see that in the
case of an integer nonpositive $s$ the operator
 \begin{equation}\label{3proj}
 (I-\mathscr P^s_{\eta,v}):\mathscr W^\ell_w(e^{i\phi}\Bbb R+\zeta; X)\to
\mathscr H^k_v(e^{i\phi}\Bbb C^++\zeta; X)
 \end{equation}
is also continuous (see the proof of
 Corollary~\ref{emb}); here $I$ denotes the operator of the embedding $\mathscr W^\ell_w(e^{i\phi}\Bbb R+\zeta;
X)\subseteq\mathscr W^k_v(e^{i\phi}\Bbb R+\zeta; X)$, we assume
that the elements of the spaces $\mathscr H^k_v(e^{i\phi}\Bbb
C^\pm+\zeta; X)$ are identified with their boundary limits in
$\mathscr W^k_v(e^{i\phi}\Bbb R+\zeta; X)$.

In Subsection~\ref{sec} we shall prove that the Fourier-Laplace
transformation yields an isomorphism between $\mathscr
H^k_v(e^{i\phi}\Bbb C^-+\zeta; X)$ and a Sobolev space of
$X$-valued distributions supported on the half-line
$\overline{e^{-i\phi}\Bbb R^++v}$ (see Paley-Wiener
Theorem~\ref{P-W thm}), then we shall use the
operator~\eqref{proj} to define a projection operator that maps
distributions defined on the line $e^{-i\phi}\Bbb R+v$ to
distributions supported on $\overline{e^{-i\phi}\Bbb R^++v}$.

\subsection{Further properties of  Hardy classes in cones}\label{ss2.3}
 In the next proposition we extend the result of Proposition~\ref{HH},
(iii) to the Hardy classes in cones.
\begin{prop}\label{col} Let $\mathcal G$ be a function defined
on the boundary $\partial \mathcal K^\varphi_\zeta=(\Bbb
R+\zeta)\cup(e^{i\varphi}\Bbb R+\zeta)$ of a cone $\mathcal
K^\varphi_\zeta$, and let $\mathcal G\in\mathscr W^\ell_w(\Bbb
R+\zeta; X)$, $\mathcal G\in \mathscr W^\ell_w(e^{i\varphi}\Bbb
R+\zeta; X)$. In the right hand sides of the
equalities~\eqref{cau} and~\eqref{1cau} we replace $\mathcal F$ by
$\mathcal G$. Then the equalities~\eqref{cau} and~\eqref{1cau}
define a function $\mathcal F\in{\mathscr H}^k_w(\mathcal
K^\varphi_\zeta;X)$, where $k<s+1/2$  if $\ell-s>1/2$, and $k\leq
s$ if $\ell-s\leq 1/2$.  (This does not necessarily mean that the
prescribed on the boundary $\partial \mathcal K^\varphi_\zeta$
function $\mathcal G$  coincides with boundary limits of the
analytic in $\mathcal K^\varphi_\zeta$ function $\mathcal F$.) The
estimate
\begin{equation}\label{est}
\|\mathcal F; \mathscr H^k_w(\mathcal K^\varphi_\zeta; X)\|\leq C
\bigl(\|\mathcal G; \mathscr W^\ell_w(\Bbb R+\zeta;X)\|+\|\mathcal
G;\mathscr W^\ell_w(e^{i\varphi}\Bbb R+\zeta;X)\|\bigr)
\end{equation}
holds, where the constant $C$ is independent of $\mathcal G$ and
$w\in\Bbb C$.
 \end{prop}
\begin{pf} We rewrite the representation~\eqref{cau} in the form
\begin{equation}\label{ew}
\mathcal F(\lambda)=\int\limits_{\Bbb R+\zeta}\frac
{e^{iw(\mu-\lambda)}(\mu-\eta)^s\mathcal G_1(\mu)}{2\pi i
(\lambda-\eta)^s(\mu-\lambda)}\,d\mu
+\int\limits_{e^{i\varphi}\Bbb R+\zeta}\frac
{e^{iw(\mu-\lambda)}(\mu-\eta)^s\mathcal G_2(\mu)}{2\pi i
(\lambda-\eta)^s(\mu-\lambda)}\,d\mu,
\end{equation}
where $\mathcal G_1(\mu)=\mathcal G(\mu)$ for $\mu\in\Bbb
R^++\zeta$ and $\mathcal G_1(\mu)=0$ for $\mu\in\Bbb R^-+\zeta$,
$\mathcal G_2(\mu)=\mathcal G(\mu)$ for $\mu\in e^{i\varphi}\Bbb
R^++\zeta$ and $\mathcal G_2(\mu)=0$ for $\mu\in e^{i\varphi}\Bbb
R^-+\zeta$. It is clear that $\mathcal G_1\in\mathscr
W^\ell_w(\Bbb R+\zeta;X)$ and $\mathcal G_2\in\mathscr
W^\ell_w(e^{i\varphi}\Bbb R+\zeta;X)$. By Proposition~\ref{HH},
(iii) the first integral in \eqref{ew} defines a function
$\mathcal F_1\in\mathscr H^k_w(\Bbb C^++\zeta;X)$, the second
integral defines a function $\mathcal F_2\in\mathscr
H^k_w(e^{i\varphi}\Bbb C^-+\zeta;X)$. The sum $\mathcal F=\mathcal
F_1+\mathcal F_2$ is analytic in $\mathcal K^{\varphi,+}_\zeta$.
We have
 the uniform in
$\psi\in(0,\varphi)$ estimates
$$
\begin{aligned}
\|\mathcal F;\mathscr W^k_w(e^{i\psi}\Bbb R^+&+\zeta;X)\|\\&\leq
\|\mathcal F_1;\mathscr W^k_w(e^{i\psi}\Bbb
R^++\zeta;X)\|+\|\mathcal F_2;\mathscr W^k_w(e^{i\psi}\Bbb
R^++\zeta;X)\|
\\ &\leq C_0\bigl(
\|\mathcal F_1;\mathscr H^k_w(\Bbb C^++\zeta;X)\|+\|\mathcal
F_2;\mathscr H^k_w(e^{i\varphi}\Bbb C^-+\zeta;X)\|\bigr)
\\ &\leq
C_1\bigl(\|\mathcal G_1;\mathscr W^\ell_w(\Bbb
R+\zeta;X)\|+\|\mathcal G_2;\mathscr W^\ell_w(e^{i\varphi}\Bbb
R+\zeta;X)\|\bigr)
\\ &\leq
C_1\bigl(\|\mathcal G;\mathscr W^\ell_w(\Bbb
R+\zeta;X)\|+\|\mathcal G;\mathscr W^\ell_w(e^{i\varphi}\Bbb
R+\zeta;X)\|\bigr);
\end{aligned}
$$
 see Proposition~\ref{pr}, (ii), the definition~\eqref{NH} of the
 norm  in $\mathscr H^\ell_w(e^{i\varphi}\Bbb C^++\zeta;X)$, and Proposition~\ref{HH}, (iii).
 We apply a
similar argument to the representation~\eqref{1cau} and conclude
that it defines an analytic in  $\mathcal K^{\varphi,-}_\zeta$
function $\mathcal F$,    the norm  $\|\mathcal F;\mathscr
W^k_w(e^{i\psi}\Bbb R^-+\zeta;X)\|$ is  bounded by the right hand
side of the estimate \eqref{est} uniformly in
$\psi\in(0,\varphi)$. Hence $\mathcal F\in\mathscr H^k_w(\mathcal
K^\varphi_\zeta;X)$ and the estimate~\eqref{est} holds.
 \qed\end{pf}

\begin{prop}\label{pr2} Let $\zeta,w\in\Bbb C$, $\ell\in\Bbb R$, and
let $\varphi\in(0,\pi]$. If $\mathcal F\in\mathscr
H^\ell_w(\mathcal K^\varphi_\zeta;X)$ then
\begin{equation}\label{esti1}
\|\mathcal F; \mathscr W^\ell_w(e^{i\psi}\Bbb R^-+\mu;X)\|\leq C
\|\mathcal F; \mathscr H^\ell_w(\mathcal
K^\varphi_\zeta;X)\|,\quad \mu\in\overline{\mathcal
K^{\varphi,-}_\zeta}, \psi\in[0,\varphi],
\end{equation}
\begin{equation}\label{esti2}
 \|\mathcal F; \mathscr W^\ell_w(e^{i\psi}\Bbb R^++\mu;X)\|\leq
C \|\mathcal F; \mathscr H^\ell_w(\mathcal
K^\varphi_\zeta;X)\|,\quad \mu\in\overline{\mathcal
K^{\varphi,+}_\zeta}, \psi\in[0,\varphi],
\end{equation}
where $\mathcal K^{\varphi,\pm}_\zeta =\{\lambda\in\mathcal
K^\varphi_\zeta:\Im\lambda\gtrless\Im\zeta\}$, the constant $C$ is
independent of $w$, $\eta$,  and $\mathcal F$.
\end{prop}
\begin{pf} Let us prove the estimate~\eqref{esti2}. The
estimate~\eqref{esti1} can be proved in a similar way. We use the
representation~\eqref{cau} with $s=\ell$ and the argument from the
proof of Proposition~\ref{col}, where we take the boundary limits
of $\mathcal F\in\mathscr H^\ell_w(\mathcal K^\varphi_\zeta;X)$ as
$\mathcal G$. This provides us with the representation $\mathcal
F=\mathcal F_1+\mathcal F_2$ in $\overline{\mathcal
K^{\varphi,+}_\zeta}$, where the functions $\mathcal
F_1\in\mathscr H^\ell_w(\Bbb C^++\zeta;X)$ and $\mathcal
F_2\in\mathscr H^\ell_w(e^{i\varphi}\Bbb C^-+\zeta;X)$ satisfy the
estimate
\begin{equation}\label{esti3}
\|\mathcal F_1;\mathscr H^\ell_w(\Bbb C^++\zeta;X)\|+\|\mathcal
F_2;\mathscr H^\ell_w(e^{i\varphi}\Bbb C^-+\zeta;X)\|\leq
C\|\mathcal F;\mathscr H^\ell_w(\mathcal K^\varphi_\zeta;X)\|;
\end{equation}
here $C$ is independent of $\mathcal F$ and $w$. Let $\mu\in
\overline{\mathcal K^{\varphi,+}_\zeta}$. For all
$\psi\in[0,\varphi]$ we have
$$
\|\mathcal F_1; \mathscr W^\ell_w(e^{i\psi}\Bbb R^++\mu; X)\|\leq
C_1\|\mathcal F_1;\mathscr H^\ell_w(\Bbb C^++\mu;X)\|\leq
C_2\|\mathcal F_1;\mathscr H^\ell_w(\Bbb C^++\zeta;X)\|,
$$
where $C_1$ and $C_2$ are independent of $\mu$, $\psi$, and
$\mathcal F_1$; see Proposition~\ref{pr}, (ii) and
Corollary~\ref{emb}, (i, ii). Similar estimates are valid for
$\mathcal F_2$. This together with~\eqref{esti3} finishes the
proof. \qed\end{pf}
\begin{cor}  Let $\zeta,w\in\Bbb C$, $\ell\in\Bbb R$, and
let $\varphi\in(0,\pi]$. The value
$$
\begin{aligned}
 \sup_{\substack{
\psi\in(0,\varphi)\\\eta\in{\mathcal
K^{\varphi,-}_\zeta}}}\|\mathcal F; \mathscr
W^\ell_w(e^{i\psi}\Bbb
R^-+\eta;X)\|+\sup_{\substack{\psi\in(0,\varphi)\\\eta\in{\mathcal
K^{\varphi,+}_\zeta}}} \|\mathcal F; \mathscr
W^\ell_w(e^{i\psi}\Bbb R^++\eta;X)\|
\end{aligned}
$$
can be taken as  an equivalent norm
$|\mspace{-2mu}|\mspace{-2mu}|\mathcal F;\mathscr
H^\ell_w(\mathcal K^\varphi_\zeta;X)
|\mspace{-2mu}|\mspace{-2mu}|$ in the Banach space $\mathscr
H^\ell_w(\mathcal K^\varphi_\zeta;X)$. Moreover, the inequalities
$$
|\mspace{-2mu}|\mspace{-2mu}|\mathcal F;\mathscr H^\ell_w(\mathcal
K^\varphi_\zeta;X) |\mspace{-2mu}|\mspace{-2mu}|\leq c\|\mathcal
F;\mathscr H^\ell_w(\mathcal K^\varphi_\zeta;X) \|\leq 2c
|\mspace{-2mu}|\mspace{-2mu}|\mathcal F;\mathscr H^\ell_w(\mathcal
K^\varphi_\zeta;X) |\mspace{-2mu}|\mspace{-2mu}|
$$
hold with an independent of $w$ constant  $c$.
\end{cor}
\begin{pf} The inequality $|\mspace{-2mu}|\mspace{-2mu}|\mathcal F;\mathscr H^\ell_w(\mathcal
K^\varphi_\zeta;X) |\mspace{-2mu}|\mspace{-2mu}|\leq c\|\mathcal
F;\mathscr H^\ell_w(\mathcal K^\varphi_\zeta;X) \|$  readily
apparent from Proposition~\ref{pr2}. It is easy to see that the
norm $\|\mathcal F;\mathscr H^\ell_w(\mathcal K^\varphi_\zeta;X)
\|$ does not exceed  the value
$|\mspace{-2mu}|\mspace{-2mu}|\mathcal F;\mathscr
H^\ell_w(\mathcal K^\varphi_\zeta;X)
|\mspace{-2mu}|\mspace{-2mu}|$. Indeed, the value $\|\mathcal F;
\mathscr W^\ell_w(e^{i\psi}\Bbb R^-+\eta;X)\|$ comes arbitrarily
close to $\|\mathcal F; \mathscr W^\ell_w(e^{i\varphi}\Bbb
R^-+\zeta;X)\|$ (or to $\|\mathcal F; \mathscr W^\ell_w(\Bbb
R^-+\zeta;X)\|$) as $\eta$ tends  to $\zeta$ in ${\mathcal
K^{\varphi,-}_\zeta}$ and $\psi\to \varphi-$ (or $\psi\to 0+$).
Similarly, the value $\|\mathcal F; \mathscr
W^\ell_w(e^{i\psi}\Bbb R^++\eta;X)\|$ comes arbitrarily close to
$\|\mathcal F; \mathscr W^\ell_w(e^{i\varphi}\Bbb R^++\zeta;X)\|$
(or to $\|\mathcal F; \mathscr W^\ell_w(\Bbb R^++\zeta;X)\|$) as
$\eta$ tends  to $\zeta$ in ${\mathcal K^{\varphi,+}_\zeta}$ and
$\psi\to \varphi-$ (or $\psi\to 0+$); cf. Proposition~\ref{p2},
1.(i). Due to the definition~\eqref{5} of the norm in $\mathscr
H^\ell_w(\mathcal K^\varphi_\zeta;X)$ this completes the proof.
\qed\end{pf}

\begin{cor}\label{Scol}  Let $\zeta,w\in\Bbb C$, $\ell\in\Bbb R$, and
 $\varphi\in(0,\pi]$. Assume that $\mathcal F\in\mathscr
H^\ell_w(\mathcal K^\varphi_\zeta;X)$.

(i) Assume in addition that  $\mathcal F$ is analytic in a
neighbourhood $\mathscr O_\zeta$ of the vertex $\zeta$ of the cone
$\mathcal K^\varphi_\zeta$. Then for any cone $\mathcal
K^\phi_\eta$, $\overline{\mathcal K^\phi_\eta}\subset \mathscr
O_\zeta\cup{\mathcal K^\varphi_\zeta}$, we have $\mathcal
F\in\mathscr H^\ell_w(\mathcal K^\phi_\eta;X)$.

(ii) Let $\mathcal K^\phi_\eta$ be a cone such that $\mathcal
K^{\phi,+}_\eta\subseteq\mathcal K^{\varphi,+}_\zeta$ (or
$\mathcal K^{\phi,-}_\eta\subseteq\mathcal K^{\varphi,-}_\zeta$).
We set $\mathcal G(\lambda)=\mathcal F(\lambda)$ for $\lambda\in
\mathcal K^{\phi,+}_\eta$ and $\mathcal G(\lambda)=0$ for
$\lambda\in\mathcal K^{\phi,-}_\eta$ (in the case $\mathcal
K^{\phi,-}_\eta\subseteq\mathcal K^{\varphi,-}_\zeta$ we set
$\mathcal G(\lambda)=\mathcal F(\lambda)$ for $\lambda\in \mathcal
K^{\phi,-}_\eta$ and $\mathcal G(\lambda)=0$ for
$\lambda\in\mathcal K^{\phi,+}_\eta$). Then $\mathcal G\in\mathscr
H^\ell_w(\mathcal K^\phi_\eta;X)$ and the inequality
\begin{equation}\label{esti4}
\|\mathcal G;\mathscr H^\ell_w(\mathcal K^\phi_\eta;X)\|\leq
C\|\mathcal F;\mathscr H^\ell_w(\mathcal K^\varphi_\zeta;X)\|
\end{equation}
holds, where the constant $C$ is independent of $w$, $\mathcal F$,
and $\mathcal K^\phi_\eta$.
\end{cor}

\begin{pf} (i)  Let $\mathscr U$ be an open neighbourhood  of the set $\overline{\mathcal
K^\phi_\eta}\setminus\mathcal K^\varphi_\zeta$ and let
$\overline{\mathscr U}\subset\mathscr O_\zeta$. We denote
 the characteristic function of $\mathscr U$ by $\chi$ (i.e.
 $\chi(\lambda)=1$ if $\lambda\in\mathscr U$ and $\chi(\lambda)=0$
 if $\lambda\notin\mathscr U$). Since $\mathcal F$ is analytic in
 $\mathscr O_\zeta$ the value
 $
 \|\chi \mathcal F;\mathscr W^\ell_w(e^{i\psi}\Bbb R+\eta;X)\|$ is bounded uniformly in $\psi\in [0,\phi]$.
The estimate
$$
 \|(1-\chi) \mathcal F;\mathscr W^\ell_w(e^{i\psi}\Bbb R+\eta;X)\|\leq
 C\|\mathcal F;\mathscr H^\ell_w(\mathcal K^\varphi_\zeta;X)\|, \psi\in [0,\phi],
$$
follows from Proposition~\ref{pr2}. Thus $\mathcal F\in \mathscr
H^\ell_w(\mathcal K^\phi_\eta;X)$ by Definition~\ref{d1}.

(ii) It is clear that $\mathcal G$ is analytic in $\mathcal
K^\phi_\eta$. The estimate \eqref{esti4} follows from
Proposition~\ref{pr2}. Hence $\mathcal G\in\mathscr
H^\ell_w(\mathcal K^\phi_\eta;X)$. \qed\end{pf}

\begin{prop}\label{restriction} Let $\zeta,w\in\Bbb C$, $s\in\Bbb R$, and
 $0<\varphi<\pi$. Assume that  $v\in e^{i\phi}\Bbb R+w$,
$\phi\in[0,\varphi]$, and  $\mathcal G=\mathscr P_{\eta,v}^s
(\mathcal F\!\!\upharpoonright_{e^{i\phi}\Bbb R+\zeta})$, where
$\mathscr P^s_{\eta,v}$ is the projection
operator~\eqref{proj11},~\eqref{proj}, and $\mathcal F$ is a
function in the space $\mathscr H^\ell_w(\mathcal
K^\varphi_\zeta;X)$, $\ell\geq s$. We also suppose that the
parameter $\eta$ in the definition~\eqref{proj11} of the operator
$\mathscr P^s_{\eta,v}$ does not belong to the union
$(e^{i\phi}\overline{\Bbb C^-}+\zeta)\cup \overline{\mathcal
K^\varphi_\zeta}$. Then the function $\mathcal G\in \mathscr
H^k_v(e^{i\phi}\Bbb C^-+\zeta; X)$ can be analytically extended to
a function $\mathcal G\in\mathscr H^k_v(\mathcal
K^\varphi_\zeta;X)$; here $k=s$ if $\ell-s\leq 1/2$, and
$k\in[s,s+1/2)$ if $\ell-s>1/2$. The estimate
\begin{equation}\label{eqd}
\|\mathcal G;\mathscr H^k_v(\mathcal K^\varphi_\zeta;X)\|\leq
Ce^{\Im\{\zeta (w-v)\}}\|\mathcal F;\mathscr H^\ell_w(\mathcal
K^\varphi_\zeta;X)\|
\end{equation}
holds, where the constant $C$ is independent of $w$, $\phi$, and
$v$.
\end{prop}
\begin{pf} We introduce the function
 $\widetilde {\mathcal
F}(\lambda)=(\lambda-\eta)^s\mathcal F(\lambda)$, $\lambda\in
\mathcal K^\varphi_\zeta$, where we use an analytic in
$(e^{i\phi}\overline{\Bbb C^-}+\zeta)\cup \overline{\mathcal
K^\varphi_\zeta}$ branch of the function $(\cdot-\eta)^s$. The
estimate
\begin{equation}\label{0001}
\|\widetilde {\mathcal F};\mathscr H^{\ell-s}_w(\mathcal
K^\varphi_\zeta;X)\|\leq C \|\mathcal F;\mathscr
H^{\ell}_w(\mathcal K^\varphi_\zeta;X)\|
\end{equation}
 holds. By Corollary~\ref{emb},(iii)  the function
$\widetilde {\mathcal F}\!\!\upharpoonright_{e^{i\phi}\Bbb
R+\zeta}\in \mathscr W^{\ell-s}_v(e^{i\phi}\Bbb R+\zeta;X)$  is
uniquely representable as the sum $({\mathcal F}^+ +{\mathcal
F}^-)\!\!\upharpoonright_{e^{i\phi}\Bbb R+\zeta}$ of the boundary
limits of some functions $\mathcal F^\pm\in\mathscr
H^{k-s}_v(e^{i\phi}\Bbb C^\pm+\zeta;X)$.  The functions $\mathcal
F^\pm$ have analytic continuations to the cone $\mathcal
K^\varphi_\zeta$ because $\widetilde{\mathcal F}=({\mathcal F}^+
+{\mathcal F}^-)$ on the line ${e^{i\phi}\Bbb R+\zeta}$,
$\widetilde{\mathcal F}$ is analytic in $\mathcal
K^\varphi_\zeta$, and $\mathcal F^\pm$ are analytic in
$e^{i\phi}\Bbb C^\pm+\zeta$. It is clear that $\mathcal
F^-=\mathscr P_{\eta,v}^0(\widetilde {\mathcal
F}\!\!\upharpoonright_{e^{i\phi}\Bbb R+\zeta})$ and $\mathcal
G(\lambda)=(\lambda-\eta)^{-s}\mathcal F^-(\lambda)$ for all
$\lambda\in e^{i\phi}\Bbb C^-+\zeta$. Thus $\mathcal G$ has an
analytic continuation to the cone $\mathcal K^\varphi_\zeta$. We
intend to demonstrate the inclusion $\mathcal F^-\in\mathscr
H^{k-s}_v(\mathcal K^\varphi_\zeta;X)$ and the estimate
\begin{equation}\label{0004}
\|\mathcal F^-;\mathscr H^{k-s}_v(\mathcal
K^\varphi_\zeta;X)\|\leq Ce^{\Im\{\zeta (w-v)\}} \|\mathcal
F;\mathscr H^\ell_w(\mathcal K^\varphi_\zeta;X)\|.
\end{equation}

Taking into account the inequality~\eqref{1proj} together with the
estimate \eqref{no2} and the norm~\eqref{5} in $\mathscr
H^\ell_w(\mathcal K^\varphi_\zeta;X)$ we get
\begin{equation}\label{000}
\begin{aligned}
\|\mathcal F^- ; \mathscr W^{k-s}_v(e^{i\psi}\Bbb R^++\zeta;
X)\|&\leq
Ce^{\Im\{\zeta (w-v)\}}   \|\widetilde{\mathcal F};\mathscr H^{\ell-s}_w(\mathcal K^\varphi_\zeta;X)\|,\ \psi\in[0,\phi],\\
\|\mathcal F^- ; \mathscr W^{k-s}_v(e^{i\psi}\Bbb R^-+\zeta;
X)\|&\leq C e^{\Im\{\zeta (w-v)\}}   \|\widetilde{\mathcal
F};\mathscr H^{\ell-s}_w(\mathcal K^\varphi_\zeta;X)\|,\
\psi\in[\phi,\varphi],
\end{aligned}
\end{equation}
where the constant $C$ is independent of $\psi$, $w$, and $v\in
e^{i\phi}\Bbb R+w$. Then  we have
\begin{equation}\label{0003}
\begin{aligned}
\|\widetilde{\mathcal F};\mathscr W^{k-s}_v(e^{i\psi}\Bbb
R^-+\zeta;X)\|\leq e^{\Im\{\zeta (w-v)\}}\|\widetilde{\mathcal
F};\mathscr W^{k-s}_w(e^{i\psi}\Bbb R^-+\zeta;X)&\|\\ \leq
e^{\Im\{\zeta (w-v)\}}\|\widetilde{\mathcal F};\mathscr
H^{\ell-s}_w(\mathcal K^\varphi_\zeta;X)\|, \ \psi\in[&0,\phi];
\\
\|\widetilde{\mathcal F};\mathscr W^{k-s}_v(e^{i\psi}\Bbb
R^++\zeta;X)\|\leq e^{\Im\{\zeta (w-v)\}}\|\widetilde{\mathcal
F};\mathscr W^{k-s}_w(e^{i\psi}\Bbb R^++\zeta;X)&\|\\ \leq
e^{\Im\{\zeta (w-v)\}}\|\widetilde{\mathcal F};\mathscr
H^{\ell-s}_w(\mathcal K^\varphi_\zeta;X)\|,\
\psi\in[&\phi,\varphi].
\end{aligned}
\end{equation}
 By analogy with \eqref{000} we derive the estimates
\begin{equation}\label{0002}
\begin{aligned}
\|\mathcal F^+ ; \mathscr W^{k-s}_v(e^{i\psi}\Bbb R^-+\zeta;
X)\|\leq
Ce^{\Im\{\zeta (w-v)\}}   \|\widetilde{\mathcal F};\mathscr H^{\ell-s}_w(\mathcal K^\varphi_\zeta;X)\|,\quad\psi\in[0,\phi],\\
\|\mathcal F^+ ; \mathscr W^{k-s}_v(e^{i\psi}\Bbb R^++\zeta;
X)\|\leq Ce^{\Im\{\zeta (w-v)\}}   \|\widetilde{\mathcal
F};\mathscr H^{\ell-s}_w(\mathcal
K^\varphi_\zeta;X)\|,\quad\psi\in[\phi,\varphi].
\end{aligned}
\end{equation}
The equality $\mathcal F^-(\lambda)=\widetilde{\mathcal
F}(\lambda)-\mathcal F^+(\lambda)$,  $\lambda\in\mathcal
K^\varphi_\zeta$, and  the estimates \eqref{000}--\eqref{0002}
establish the uniform in $\psi\in[0,\varphi]$ estimate
$$
\|\mathcal F^- ; \mathscr W^{k-s}_v(e^{i\psi}\Bbb R+\zeta;
X)\|\leq Ce^{\Im\{\zeta (w-v)\}}   \|\widetilde{\mathcal
F};\mathscr H^{\ell-s}_w(\mathcal K^\varphi_\zeta;X)\|.
$$
This together with~\eqref{0001} leads to~\eqref{0004}. Since
$\mathcal G(\lambda)=(\lambda-\eta)^{-s}\mathcal F^-(\lambda)$ we
conclude that
$$
\|\mathcal G;\mathscr H^k_v(\mathcal K^\varphi_\zeta;X)\|\leq C
\|\mathcal F^- ;\mathscr H^{k-s}_v(\mathcal K^\varphi_\zeta;X)\|.
$$
Then the estimate~\eqref{eqd} follows from~\eqref{0004}.
\qed\end{pf}

By analogy with $\mathscr H^\ell_w(\mathcal K^{\varphi}_\zeta;X)$
we introduce the  Hardy class $\mathscr H^\ell_w(\mathcal
K^{-\varphi}_\zeta;X)$
 in the cone
$ \mathcal K^{-\varphi}_\zeta=\{\lambda\in\Bbb
C:\lambda=e^{-i\varphi}\mu+\zeta, \mu\in\mathcal K^\varphi_0\}$;
here $\varphi\in(0,\pi]$. The class $\mathscr H^\ell_w(\mathcal
K^{-\varphi}_\zeta;X)$ consists of all analytic functions
$\mathcal K^{-\varphi}_\zeta\ni\lambda\mapsto \mathcal
F(\lambda)\in X$ satisfying the uniform in $\psi\in (-\varphi,0)$
estimate~\eqref{no1}. Let $\varrho_{\varphi,\zeta}:\mathcal
K^\varphi_\zeta\to\mathcal K^{-\varphi}_\zeta$ denote the linear
transformation
$\varrho_{\varphi,\zeta}(\lambda)=e^{-i\varphi}(\lambda-\zeta)+\zeta$.
It is easily seen that we can identify the classes $\mathscr
H^\ell_w(\mathcal K^{-\varphi}_\zeta;X)$ and $\mathscr
H^\ell_v(\mathcal K^{\varphi}_\zeta;X)$  by the rule
 $\mathcal G=\mathcal
F\circ \varrho_{\varphi,\zeta}$, where $\mathcal F\in \mathscr
H^\ell_w(\mathcal K^{-\varphi}_\zeta;X)$, $\mathcal G\in\mathscr
H^\ell_v(\mathcal K^{\varphi}_\zeta;X)$, and $v=e^{-i\varphi}w$.
Thus far, after obvious changes, everything said about the classes
$\mathscr H^\ell_w(\mathcal K^{\varphi}_\zeta;X)$ is applicable to
$\mathscr H^\ell_w(\mathcal K^{-\varphi}_\zeta;X)$. In particular,
a function $\mathcal F\in \mathscr H^\ell_w(\mathcal
K^{-\varphi}_\zeta;X)$ has boundary limits $\mathcal
F_0\in\mathscr W^\ell_w(\Bbb R+\zeta;X)$ and $\mathcal
F_{-\varphi}\in\mathscr W^\ell_w(e^{-i\varphi}\Bbb R+\zeta;X)$.

\begin{prop}\label{H^s} Let $\zeta,w\in\Bbb C$, $\ell\in\Bbb R$, and
 let $0<\varphi<\pi/2$. Then every function $\mathcal F\in\mathscr
W^\ell_w(\Bbb R+\zeta;X)$ can be represented as $\mathcal
F^+_0+\mathcal F^-_0$, where $\mathcal F^+_0$ and $\mathcal F^-_0$
are the boundary limits of some functions $\mathcal F^+\in\mathscr
H^\ell_w(\mathcal K^\varphi_\zeta;X)$ and $\mathcal F^-\in\mathscr
H^\ell_w(\mathcal K^{-\varphi}_\zeta;X)$.
\end{prop}
\begin{pf} Let
 $\eta$ be in the set $\Bbb C\setminus(\overline{\mathcal K^\varphi_\zeta}
 \cup\overline{ \mathcal K^{-\varphi}_\zeta})$, this set is not empty because $\varphi<\pi/2$.
  Consider the
  function ${\mathcal G}(\lambda)=
  (\lambda-\eta)^\ell\mathcal F(\lambda)$, where
$\lambda\in \Bbb R+\zeta$, we  use an analytic in
$\overline{\mathcal K^\varphi_\zeta}
 \cup\overline{ \mathcal K^{-\varphi}_\zeta}$ branch of the
 function $(\cdot-\eta)^\ell$. Due to Corollary~\ref{emb}, (iii) the function  $\mathcal G\in\mathscr W^0_w(\Bbb
 R+\zeta;X)$ can be represented as the sum $(\mathcal G^+ +\mathcal G^-)\!\!\upharpoonright_{\Bbb
 R+\zeta}$ of boundary limits of some functions $\mathcal G^\pm\in\mathscr H^0_w(\Bbb
 C^\pm+\zeta;X)$. Therefore we
 can  define the functions $\mathcal F^+\in\mathscr
H^\ell_w(\mathcal K^\varphi_\zeta;X)$ and $\mathcal F^-\in\mathscr
H^\ell_w(\mathcal K^{-\varphi}_\zeta;X)$ by the formulas
$$
\mathcal F^+(\lambda)=\left\{%
\begin{array}{ll}
    (\lambda-\eta)^{-\ell}\mathcal G^+(\lambda), & \lambda\in\mathcal K^{\varphi}_\zeta \cap \Bbb C^++\zeta; \\
    (\lambda-\eta)^{-\ell}\mathcal G^-(\lambda), & \lambda\in\mathcal K^{\varphi}_\zeta \cap \Bbb C^-+\zeta;\\
\end{array}%
\right.
$$
$$
\mathcal F^-(\lambda)=\left\{%
\begin{array}{ll}
    (\lambda-\eta)^{-\ell}\mathcal G^+(\lambda), & \lambda\in\mathcal K^{-\varphi}_\zeta \cap \Bbb C^++\zeta; \\
    (\lambda-\eta)^{-\ell}\mathcal G^-(\lambda), & \lambda\in\mathcal K^{-\varphi}_\zeta \cap \Bbb C^-+\zeta.\\
\end{array}%
\right.
$$
\qed\end{pf}

\subsection{Fredholm polynomial operator pencils in spaces of analytic functions}\label{OP}

Let $X_j$ denote a Hilbert space with the norm $\|\cdot\|_j$. We
introduce a set $ \{X_j\}_{j=0}^m$ of Hilbert spaces such that
$\|u\|_{j}\leq \|u\|_{j+1}$ and
   $X_{j+1}$ is dense in $X_{j}$ for all $j=0,\dots,m-1$. Let
   $\{A_j\in\mathscr B(X_j,X_0)\}_{j=0}^m$ be a set of operators, where
$\mathscr B(X_j,X_0)$ stands for the set of all linear bounded
operators $A_j:X_j\to X_0$. We consider the  operator pencil
\begin{equation}\label{pencil}
\Bbb C\ni\lambda\mapsto \mathfrak A(\lambda)=\sum_{j=0}^m
A_j\lambda^{m-j}\in \mathscr B (X_m, X_0).
\end{equation}
 We assume that the operator $\mathfrak A(\lambda)$ is Fredholm
for all $\lambda\in\Bbb C$ and is invertible for at least one
value of $\lambda$.  Under these assumptions the operator
$\mathfrak A(\lambda)$ is invertible for all $\lambda\in\Bbb C$
except for isolated eigenvalues of the pencil~\eqref{pencil}.
These eigenvalues are of finite algebraic multiplicities and can
accumulate only at infinity. We also assume that the
pencil~\eqref{pencil} satisfies the {\bf condition: }{\it\  there
exist $\vartheta\in(0,\pi/2)$ and $R>0$ such that
 for all $f\in X_0$ the following estimate is fulfilled
\begin{equation}\label{cndtn}
\sum_{j=0}^m |\lambda|^j\|\mathfrak A^{-1}(\lambda)f\|_{m-j}\leq c
\|f\|_0, \quad \lambda\in\overline{\mathcal K^\vartheta_0}\cup
\overline{\mathcal K^{-\vartheta}_0}, |\lambda|>R.
\end{equation}}
The condition~\eqref{cndtn} follows from the assumption~\eqref{M
cond}, see~\cite[Proposition~2.2.1]{R17}. As it was already
mentioned in the introductory part, the assumption~\eqref{M cond}
 is widely met in the theory of operator pencils
and
 is satisfied in many applications of the theory of boundary value
  problems for partial differential equations; see e.g.~\cite{R12},~\cite{R17,ref7} and  references
 therein. It guaranties that for any $\zeta\in\Bbb C$ and any $\varphi\in(0,\vartheta)$
 the closed cones $\overline{\mathcal K^\varphi_\zeta}$ and $\overline{\mathcal K^{-\varphi}_\zeta}$ contain at most finitely
 many eigenvalues of the pencil $\mathfrak A$, one can find $\zeta$ and $\varphi$ so that $\overline{\mathcal K^\varphi_\zeta}$
 and/or $\overline{\mathcal K^{-\varphi}_\zeta}$ are free from the
 spectrum of $\mathfrak A$.

For $\ell\in \Bbb R$, $\varphi\in(0,\pi]$, and $\zeta,w\in\Bbb C$
we introduce the Banach space $\mathfrak D^\ell_w(\mathcal
 K^\varphi_\zeta)$ of analytic functions $\mathcal
 K^\varphi_\zeta\ni\lambda\mapsto u(\lambda)\in X_m$ and the norm in $\mathfrak
D^\ell_w(\mathcal
 K^\varphi_\zeta)$ by the equalities
\begin{equation}\label{sp}
\mathfrak D^\ell_w(\mathcal
 K^\varphi_\zeta)=\bigcap_{j=0}^m \mathscr H^{\ell-j}_w(\mathcal
 K^\varphi_\zeta;X_j); \ \
 \|u;\mathfrak D^\ell_w(\mathcal
 K^\varphi_\zeta) \|=\sum_{j=0}^m\|u;\mathscr H^{\ell-j}_w(\mathcal K^\varphi_\zeta;X_j)\|.
\end{equation}
It is convenient to introduce also the Banach space
\begin{equation}\label{sp'}
\begin{aligned}
\mathfrak D^\ell_w(e^{i\psi}\Bbb R+\zeta)&=\bigcap_{j=0}^m
\mathscr
W^{\ell-j}_w(e^{i\psi}\Bbb R+\zeta;X_j); \\
 \|u;\mathfrak D^\ell_w(e^{i\psi}\Bbb R+\zeta) \|&=\sum_{j=0}^m\|u;\mathscr W^{\ell-j}_w(e^{i\psi}\Bbb R+\zeta;X_j)\|.
\end{aligned}
\end{equation}

\begin{prop}\label{sp prop} Let $\zeta,w\in\Bbb C$, $\ell\in\Bbb R$,
 and $\varphi\in(0,\pi]$. The following assertions are valid.

 (i) Every function $u\in \mathfrak D^\ell_w(\mathcal
 K^\varphi_\zeta)$ has boundary limits
 $u_0\in \mathfrak D^\ell_w(\Bbb R+\zeta)$ and $u_\varphi\in\mathfrak D^\ell_w(e^{i\varphi}\Bbb
 R+\zeta)$
in the sense that for almost all points $\mu$ of the boundary
$\partial\mathcal K^\varphi_\zeta=(\Bbb
R+\zeta)\cup(e^{i\varphi}\Bbb
 R+\zeta)$ we have $\|u(\lambda)-u_0(\mu)\|_m\to 0$ as $\lambda$ non-tangentially tends to
 $\mu\in \Bbb R+\zeta$, and $\|u(\lambda)-u_\varphi(\mu)\|_m\to 0$
as $\lambda$ non-tangentially tends to
 $\mu\in e^{i\varphi}\Bbb R+\zeta$. Moreover,
 \begin{equation}\label{mrvr}
\begin{aligned}
\|(e_w u)\circ\varkappa_{\psi,\zeta}-(e_w
u_0)\circ\varkappa_{0,\zeta};\mathfrak D^\ell_0(\Bbb R)\|&\to
0,\quad\psi\to 0+,
\\
\|(e_w
u)\circ\varkappa_{\psi,\zeta}-(e_wu_\varphi)\circ\varkappa_{\varphi,\zeta};\mathfrak
D^\ell_0(\Bbb R)\|&\to 0,\quad\psi\to \varphi-;
\end{aligned}
\end{equation}
recall that $e_w:\lambda\mapsto\exp(i w \lambda)$ and
$\varkappa_{\psi,\zeta}:\Bbb R\to e^{i\psi}\Bbb R+\zeta$ is the
linear transformation
$\varkappa_{\psi,\zeta}(\xi)=e^{i\psi}\xi+\zeta$. From now on we
suppose that every element $u\in \mathfrak D^\ell_w(\mathcal
 K^\varphi_\zeta)$ is extended
 to the boundary $\partial\mathcal
K^\varphi_\zeta$ by its non-tangential limits.

  (ii) For all $\psi\in[0,\varphi]$ and
$u\in\mathfrak D^\ell_w(\mathcal
 K^\varphi_\zeta)$ the estimate
\begin{equation}\label{ub}
\|u;\mathfrak D^\ell_w(e^{i\psi}\Bbb R+\zeta)\|\leq C\|u;\mathfrak
D^\ell_w(\mathcal
 K^\varphi_\zeta) \|
\end{equation}
holds, where the constant $C$ is independent of $u$, $\psi$, and
$w$.
\end{prop}
\begin{pf}Due to the embedding $\mathfrak D^\ell_w(\mathcal
 K^\varphi_\zeta)\subset \mathscr H^{\ell-m}_w(\mathcal K^\varphi_\zeta;
 X_m)$ and  Proposition~\ref{p2}.1.(i), the function $u\in \mathfrak D^\ell_w(\mathcal
 K^\varphi_\zeta)$ has the boundary limit $u_0$ in the sense that
  for almost all points $\mu\in \Bbb
 R+\zeta$
   we have $\|u(\lambda)-u_0(\mu)\|_m\to
 0$ as $\lambda$ non-tangentially tends to $\mu$.  Under the assumptions made on the spaces
 $X_j$, the convergence of $u(\lambda)$ to $u_0(\mu)$ in  $X_m$
 gives
 $u(\lambda)\to u_0(\mu)$ in  $X_j$ for any $j=0,\dots,m$.
  This together with Proposition~~\ref{p2}.1.(i) and the embedding
$\mathfrak D^\ell_w(\mathcal
 K^\varphi_\zeta)\subset \mathscr H^{\ell-j}_w(\mathcal K^\varphi_\zeta;
 X_j)$ provides us with the inclusion $u_0\in \mathsf W^{\ell-j}_\zeta (\Bbb
 R+w;X_j)$, $j=0,\dots,m$. Hence the boundary limit $u_0$  satisfies the
 inclusion~$u_0\in \mathfrak D^\ell_w(\Bbb R+\zeta)$. A similar argument demonstrates the
inclusion  $u_\varphi\in\mathfrak D^\ell_w(e^{i\varphi}\Bbb
 R+\zeta)$.
 Now the assertion~(ii) is readily
 apparent from Proposition~\ref{p2}.1.(ii) and the
 definition~\eqref{sp} of the norm in $\mathfrak D^\ell_w(\mathcal
 K^\varphi_\zeta)$.
 The first relation in~\eqref{mrvr} is valid since the left hand side of the
 estimate~\eqref{ub} is bounded uniformly in $\psi\in[0,\varphi]$,  and
   $\|u\circ\varkappa_{\psi,\zeta}(\xi)-u_0\circ\varkappa_{0,\zeta}(\xi)\|_j\to
 0$  as $\psi\to 0+$ for $j=0,\dots,m$ and almost all $\xi\in\Bbb
 R$. In the same way one can see the second relation
 in~\eqref{mrvr}.
\qed\end{pf}

\begin{thm}\label{1} Let $\zeta,w\in\Bbb C$ and $\ell\in\Bbb R$. The following assertions are valid.

(i) Let  $\varphi\in(0,\pi]$. The operator pencil~\eqref{pencil}
implements the continuous mapping
\begin{equation}\label{map}
\mathfrak D^\ell_w(\mathcal
 K^\varphi_\zeta)\ni u\mapsto \mathfrak A u=\mathcal F\in\mathscr H^{\ell-m}_w(\mathcal
 K^\varphi_\zeta;X_0);
\end{equation}
here and elsewhere $\mathfrak A u$ stands for the function
$\overline{\mathcal K^\varphi_\zeta}\ni\lambda\mapsto \mathfrak
A(\lambda) u(\lambda)\in
 X_0$.
The estimate
\begin{equation}\label{ese}
\|\mathfrak A u;  \mathscr H^{\ell-m}_w(\mathcal
K^\varphi_\zeta;X_0) \|\leq C\|u;  \mathfrak D^\ell_w(\mathcal
 K^\varphi_\zeta)  \|
\end{equation}
holds with an independent of $w\in\Bbb C$ and $u\in \mathfrak
D^\ell_w(\mathcal
 K^\varphi_\zeta)$ constant $C$.

(ii) Suppose that the pencil $\mathfrak A$ satisfies the
condition~\eqref{cndtn} for some $\vartheta\in(0,\pi/2)$ and
$R>0$. Let the closed cone $\overline{\mathcal K^\varphi_\zeta}$,
$\varphi<\vartheta$, be free from the eigenvalues of the pencil
$\mathfrak A$. Then the mapping~\eqref{map} is an isomorphism, the
estimate
\begin{equation}\label{lp}
\|\mathfrak A^{-1}\mathcal F; \mathfrak D^\ell_w(\mathcal
 K^\varphi_\zeta)  \|\leq C\| \mathcal F;  \mathscr H^{\ell-m}_w(\mathcal K^\varphi_\zeta;X_0) \|
\end{equation}
holds. Here the constants $C$ is independent of $w\in\Bbb C$ and
$\mathcal F\in\mathscr H^{\ell-m}_w(\mathcal
 K^\varphi_\zeta;X_0)$.

(iii) Suppose that the pencil $\mathfrak A$ satisfies the
condition~\eqref{cndtn} for some $\vartheta\in(0,\pi/2)$ and
$R>0$. Let the line $e^{i\psi}\Bbb R+\zeta $, $|\psi|<\vartheta$,
be free from the eigenvalues of the pencil $\mathfrak A$. Then the
 mapping
$$
\mathfrak D^\ell_w(e^{i\psi}\Bbb R+\zeta)\ni u\mapsto \mathfrak A
u=\mathcal F\in\mathscr W^{\ell-m}_w(e^{i\psi}\Bbb R+\zeta;X_0)
$$
 is an isomorphism, and the estimates
\begin{equation*}\label{lp'}\begin{aligned}
\|\mathfrak A^{-1}\mathcal F; \mathfrak D^\ell_w(e^{i\psi}\Bbb
R+\zeta)  \|\leq c_1\| \mathcal F;  \mathscr
W^{\ell-m}_w(e^{i\psi}\Bbb R+\zeta;X_0) \|
\\\leq c_2 \|\mathfrak
A^{-1}\mathcal F; \mathfrak D^\ell_w(e^{i\psi}\Bbb R+\zeta)  \|
\end{aligned}
\end{equation*}
are valid with some independent of $w\in\Bbb C$ and $\mathcal F\in
\mathscr W^{\ell-m}_w(e^{i\psi}\Bbb R+\zeta;X_0)$ constants $c_1$,
$c_2$.
\end{thm}
\begin{pf}
 (i) Recall that the set $\{A_j\in\mathscr B(X_j,X_0)\}_{j=0}^m$
consists of
 bounded, independent of $w\in\Bbb C$ operators. For the
 pencil~\eqref{pencil} we have
\begin{equation}\label{in1}
\|\mathfrak A(\lambda)u(\lambda)\|^2_0\leq\Bigl(\sum_{j=0}^m
|\lambda|^{m-j}\|A_j u(\lambda)\|_0\Bigr)^2\leq C \sum_{j=0}^m(1+
|\lambda|^2)^{(m-j)}\|u(\lambda)\|^2_j,
\end{equation}
 where $C$ is independent of $w\in\Bbb C$ and $u\in \mathfrak
D^\ell_w(\mathcal
 K^\varphi_\zeta) $. Multiplying the inequalities~\eqref{in1} by the factor
 $|\exp\{2iw\lambda\}|(1+|\lambda|^2)^{\ell-m}$ and integrating
 the result with respect to $\lambda\in e^{i\psi}\Bbb R+\zeta$, we
 arrive at the estimate
 \begin{equation*}
\|\mathfrak A u;  \mathscr W^{\ell-m}_w(e^{i\psi}\Bbb R+\zeta;X_0)
\|^2\leq C\sum_{j=0}^m \|u;  \mathscr W^{\ell-j}_w(e^{i\psi}\Bbb
R+\zeta;X_j)  \|^2,\quad \psi\in[0,\varphi].
\end{equation*}
This together with~\eqref{ub} gives the estimate
 \begin{equation}\label{*1}
\|\mathfrak A u;  \mathscr W^{\ell-m}_w(e^{i\psi}\Bbb R+\zeta;X_0)
\|\leq C \|u;  \mathfrak D^\ell_w(\mathcal
 K^\varphi_\zeta)  \|
\end{equation}
with an independent of $w\in\Bbb C$, $\psi\in[0,\varphi]$, and
$u\in \mathfrak D^\ell_w(\mathcal
 K^\varphi_\zeta)$ constant $C$. Therefore the analytic function
 $\mathcal K^\varphi_\zeta\ni\lambda\mapsto \mathfrak A(\lambda) u(\lambda)\in
 X_0$ satisfies the uniform in $\psi\in[0,\varphi]$
 estimate~\eqref{*1}. Thus $\mathfrak Au\in \mathscr H^{\ell-m}_w(\mathcal
 K^\varphi_\zeta;X_0)$ and the estimate~\eqref{ese} is valid.

(ii) The cone $\overline{\mathcal K^\varphi_\zeta}$ is free from
the eigenvalues of the pencil $\mathfrak A$,  consequently the
operator function $\overline{\mathcal
K^\varphi_\zeta}\ni\lambda\mapsto \mathfrak
A^{-1}(\lambda)\in\mathscr B(X_0,X_m)$ is holomorphic. For all
$f\in X_0$ the estimate
$$
\|\mathfrak A^{-1}(\lambda)f\|_m\leq C(R)\| f\|_0,\quad \lambda\in
\overline{\mathcal K^\varphi_\zeta}, |\lambda|<R+1,
$$
is valid; here $R$ is a sufficiently large positive number. Hence,
for all $f\in X_0$  we have
$$
\sum_{j=0}^m|\lambda|^j\|\mathfrak A^{-1}(\lambda)f\|_{m-j}\leq
c(R)\| f\|_0,\quad \lambda\in \overline{\mathcal K^\varphi_\zeta},
|\lambda|<R+1.
$$
This together with the estimate~\eqref{cndtn} gives
\begin{equation}\label{b}
\sum_{j=0}^m|\lambda|^j\|\mathfrak A^{-1}(\lambda)\mathcal
F(\lambda)\|_{m-j}\leq c(R)\| \mathcal F(\lambda)\|_0,\quad
\lambda\in \overline{\mathcal K^\varphi_\zeta},
\end{equation}
 where $\mathcal F\in
\mathscr H^{\ell-m}_w (\mathcal K^\varphi_\zeta;X_0)$. We can
rewrite the inequality~\eqref{b} in the form
\begin{equation}\label{b1}
\sum_{j=0}^m(1+|\lambda|^2)^j\|\mathfrak A^{-1}(\lambda)\mathcal
F(\lambda)\|^2_{m-j}\leq C\| \mathcal F(\lambda)\|^2_0,\quad
\lambda\in \overline{\mathcal K^\varphi_\zeta};
\end{equation}
here the constant $C$ is independent of $w\in\Bbb C$ and $\mathcal
F\in \mathscr H^{\ell-m}_w (\mathcal K^\varphi_\zeta;X_0)$. We
multiply the inequality~\eqref{b1} by the factor
 $|\exp\{2iw\lambda\}|(1+|\lambda|^2)^{\ell-m}$ and integrate
 the result with respect to $\lambda\in e^{i\psi}\Bbb R+\zeta$. We
get
 \begin{equation}\label{z1}
\sum_{j=0}^m \|\mathfrak A^{-1}\mathcal F;  \mathscr
W^{\ell-j}_w(e^{i\psi}\Bbb R+\zeta;X_j)  \|^2\leq C\|\mathcal F;
\mathscr W^{\ell-m}_w(e^{i\psi}\Bbb R+\zeta;X_0) \|^2,\quad
\psi\in[0,\varphi].
\end{equation}
Clearly this leads to the estimate~\eqref{lp}. The assertion (ii)
is proved.

The proof of the assertion (iii) is similar. \qed\end{pf}

\section{Spaces of vector valued distributions}\label{sec2}

In this section we study weighted Hardy-Sobolev spaces $\mathsf
H^\ell_\zeta(K^\varphi_w; X)$ of vector valued distributions. The
space $\mathsf H^\ell_\zeta(K^\varphi_w; X)$ consists of
Fourier-Laplace transforms of all functions from the Hardy space
$\mathscr H^\ell_w(\mathcal K^\varphi_\zeta;X)$. In the case
$\ell\geq 0$  the elements of the space $\mathsf
H^\ell_w(K^\varphi_\zeta; X)$ are $X$-valued functions analytic in
the cone $K^\varphi_\zeta$. In Subsection~\ref{FL} we introduce
the Fourier-Laplace transformation acting in weighted spaces of
vector valued distributions. Then in Subsection~\ref{ss3.2} we
introduce and study the spaces $\mathsf H^\ell_\zeta(K^\varphi_w;
X)$. Subsection~\ref{sec} assembles different aspects related to
the Sobolev, Hardy, and Hardy-Sobolev spaces. In particular we
prove a variant of the Paley-Wiener theorem (Theorem~\ref{P-W
thm}) and a variant of the Paley-Wiener-Schwartz theorem
(Proposition~\ref{PWS}).

\subsection{Fourier-Laplace transformation}\label{FL}
Let $C_0^\infty(\Bbb R; X)$ denote the space of smooth compactly
supported functions with values in $X$. For functions from
$C^\infty_0(\Bbb R;X)$ we can define the Fourier transform
$G=\mathscr F_{\xi\to t}\mathcal G$ and the inverse Fourier
transform $\mathcal G=\mathscr F^{-1}_{t\to\xi}G$ by the formulas
\begin{equation}\label{ft}
G(t)=\frac 1{\sqrt{2\pi}}\int_{-\infty}^{+\infty} e^{i\xi
t}\mathcal G(\xi)\,d\xi;\quad \mathcal G(\xi)=\frac 1
{\sqrt{2\pi}}\int_{-\infty}^{+\infty}e^{-i\xi t} G(t)\,dt.
\end{equation}
The integrals in \eqref{ft} are absolutely convergent in $X$, the
Parseval equality
\begin{equation}\label{PEQ}
\|G; L_2(\Bbb R;X)\|=\|\mathcal G; L_2(\Bbb R;X)\|
\end{equation}
holds. As is well-known \cite{R4,R15}, the Fourier transformation
\eqref{ft} can be continuously extended  to an isometric
isomorphism $\mathscr F_{\xi\to t}:L_2(\Bbb R;X)\to L_2(\Bbb
R;X)$.
\begin{lem}\label{lno1}
Let $\zeta,w\in \Bbb C$ and let $\psi$ be an angle.
 For $\mathcal F\in C_0^\infty(e^{i\psi}\Bbb
R+\zeta;X)$ we define the Fourier-Laplace transformation  $\Bbb
T^\psi_{\zeta,w}:\mathcal F\mapsto  F $ by the equality
\begin{equation}\label{noTi} F(z)=\frac {e^{i\zeta w}} {\sqrt{2\pi}}
\int_{e^{i\psi}\Bbb R+\zeta}e^{i z\lambda}\mathcal F(\lambda) \,
d\lambda,\quad z\in e^{-i\psi}\Bbb R+w,
\end{equation}
where the integration runs from $e^{i\psi}(-\infty)+\zeta$ to
$e^{i\psi}(+\infty)+\zeta$. Then the inverse transformation $(\Bbb
T^\psi_{\zeta,w})^{-1}$
 is given by the
equality
\begin{equation}\label{noT}
 \mathcal F(\lambda)= \frac {e^{-i\zeta w}}
{\sqrt{2\pi}}\int_{e^{-i\psi}\Bbb R+w}e^{-i\lambda z} F(z) \,
dz,\quad \lambda\in e^{i\psi}\Bbb R+\zeta,
\end{equation}
where the integration runs from $e^{-i\psi}(-\infty)+\zeta$ to
$e^{-i\psi}(+\infty)+\zeta$. The Parseval equality
 \begin{equation}\label{PE}\|\mathcal F;\mathscr W^0_w(e^{i\psi}\Bbb
R+\zeta;X)\|=\|F;\mathscr W^0_{-\zeta}(e^{-i\psi}\Bbb R+w;X)\|
\end{equation}
holds.  The transformation $\Bbb T^\psi_{\zeta,w}$ can be
continuously extended  to an isometric isomorphism
\begin{equation}\label{noTsp}\Bbb T^\psi_{\zeta,w}:
\mathscr W^0_w(e^{i\psi}\Bbb R+\zeta;X)\to\mathscr
W^0_{-\zeta}(e^{-i\psi}\Bbb R+w;X).
\end{equation}
(The multipliers $e^{i\zeta w}$ and $e^{-i\zeta w}$ in the
formulas~\eqref{noTi} and~\eqref{noT} are added so that
 the transformation~\eqref{noTsp} is an
isometry.)
\end{lem}
\begin{pf} One can get the formulas \eqref{noTi},
\eqref{noT}, and \eqref{PE} from the formulas for the Fourier
transformation \eqref{ft} and \eqref{PEQ}. Indeed, it suffices to
 set $$\mathcal G(\xi)=\exp\{iw
(e^{i\psi}\xi+\zeta)\}\mathcal F(e^{i\psi}\xi+\zeta),\quad
\xi\in\Bbb R,$$
$$F (e^{-i\psi}t+w)=\exp\{i\psi+i\zeta (e^{-i\psi} t+w)\}G(t),\quad t\in\Bbb R,$$
and $\lambda=e^{i\psi}\xi+\zeta$, $z=e^{-i\psi}t+w$. \qed\end{pf}

  Let $\mathscr S(\Bbb R)$
be the Schwartz space of rapidly decreasing functions on $\Bbb R$.
By $\mathscr S'(\Bbb R; X)$ we denote the space of tempered
distributions with values in $X$, it consists of all linear
continuous mappings $\mathscr S\to X$. For all $\mathcal G\in
\mathscr S'(\Bbb R;X)$ we define the generalized Fourier
transformation $G=\mathcal F_{\xi \to t} \mathcal G$  by the
equality $ \mathcal G(v)=G (\hat v)$, $\forall v\in \mathscr
S(\Bbb R)$, where $\hat v(t)=(2\pi)^{-1/2}\int_{\Bbb R} e^{i\xi t}
v(\xi)\,d\xi$. As is known \cite{R4,R10}, this
 transformation yields an isomorphism of $\mathscr S'(\Bbb
R; X)$ onto itself. Let $D_t=-i\partial_t$. For all $\mathcal G\in
\mathscr S'(\Bbb R; X)$ we have
\begin{equation}\label{nooo4}
\mathscr F_{\xi\to t} \xi^j\mathcal G=D^j_t \mathscr F_{\xi\to
t}\mathcal G;\quad D^j_t \mathscr F_{\xi\to t}\mathcal G
(v)=\mathscr F_{\xi\to t}\mathcal G(D^j_t v), \quad\forall
v\in\mathscr S(\Bbb R).
\end{equation}
 Every element $\mathcal G\in \mathscr
S'(\Bbb R;X)$ can be represented (in non-unique fashion) in the
form $\sum_{m+n\leq\ell}D_\xi^m(1+|\xi|)^n \mathcal G_{mn}$, where
$\ell$ is finite and $\mathcal G_{mn}\in L_2(\Bbb R;X)$; i.e.
$$
\mathcal G(v)=\sum_{m+n\leq\ell}\int_{\Bbb R}\mathcal G_{mn}(\xi)
(1+|\xi|)^n\overline{ D_\xi^m v (\xi)}\, d\xi \quad\forall
v\in\mathscr S(\Bbb R),
$$
the integral is absolutely convergent in $X$. We have the
embedding $L_2(\Bbb R;X)\subset \mathscr S'(\Bbb R;X)$,
$\ell\in\Bbb R$, which means that to any $\mathcal F\in L_2(\Bbb
R;X)$ there corresponds a disribution $\mathcal F\in \mathscr
S'(\Bbb R;X)$ given by the formula $ \mathcal F(v)=\int\mathcal
F(\xi)\overline{v(\xi)}\,d\xi$, $v\in\mathscr S(\Bbb R)$.

By $\mathscr D'(e^{-i\psi}\Bbb R+w;X)$ we denote the space of all
linear continuous mappings $C_0^\infty(e^{-i\psi}\Bbb R+w)\to X$;
here $w\in\Bbb C$ and $\psi$ is an angle. Let $\tau_{\psi,w}:\Bbb
R\to e^{-i\psi}\Bbb R+w$ be the
linear transformation $\tau_{\psi,w}(t)=e^{-i\psi}t+w$. 
We say that
 $\mathcal F\in \mathscr D'(e^{-i\psi}\Bbb R+w; X)$ if and
only if $\mathcal F\circ\tau_{\psi,w}\in \mathscr D'(\Bbb R;X)$,
where $\mathcal F\circ\tau_{\psi,w}(v)=\mathcal F
(v\circ\tau^{-1}_{\psi,w})$, $v\in C_0^\infty(\Bbb R)$. In the
same manner we introduce the space $\mathscr S'(e^{-i\psi}\Bbb
R+w;X)$  of tempered distributions on $e^{-i\psi}\Bbb R+w$. The
embedding $\mathscr S'(e^{-i\psi}\Bbb R+w;X)\subset\mathscr
D'(e^{-i\psi}\Bbb R+w;X)$ holds.

 We define the operator  $
 D_{(\psi)}: \mathscr D'(e^{-i\psi}\Bbb
R+w;X)\to \mathscr D'(e^{-i\psi}\Bbb R+w;X) $ of differentiation
along the line $e^{-i\psi}\Bbb R+w$ by the equality
\begin{equation}\label{D}
 D_{(\psi)} F=e^{i\psi}\bigl( D_t
(F\circ\tau_{\psi,w})\bigr)\circ\tau^{-1}_{\psi,w},\quad F\in
\mathscr D'(e^{-i\psi}\Bbb R+w;X).
\end{equation}
 Let $\mathscr O_x$ be a neighbourhood in $\Bbb C$ of a point
$x\in e^{-i\psi}\Bbb R+w$ and let $\mathscr O_x\ni z\mapsto
G(z)\in X$ be an analytic function. It can be easily seen that in
the intersection $\mathscr O_x\cap(e^{-i\psi}\Bbb
 R+w)$ we have
 $(D_z G)\!\!\upharpoonright_{e^{-i\psi}\Bbb
 R+w}=D_{(\psi)} (G\!\!\upharpoonright_{e^{-i\psi}\Bbb
 R+w}) $, where $D_z=-\frac i 2(\partial_{\Re z}-i\partial_{\Im
 z})$ is the complex derivative. Alternatively, we can define at first the operator $D_{(\psi)}$ in $C_0^\infty (e^{-i\psi}\Bbb
 R+w)$ by the equality $D_{(\psi)} v=(D_z \tilde v)\!\!\upharpoonright_{e^{i\psi}\Bbb
 R+w}$,
  where $\tilde v$ is an almost analytic
extension of $v$ ($\tilde v$ is a smooth extension of $v$ to a
neighbourhood of the line $e^{-i\psi}\Bbb R+w$ such that
$D_{\overline z}\tilde v$ vanishes to infinite order on
$e^{-i\psi}\Bbb R+w$, see e.g. \cite{R11}). Then we extend
$D_{(\psi)}$  to all  $F\in\mathscr D'(e^{-i\psi}\Bbb R+w;X)$.

\begin{lem}\label{lno2} Let $\zeta,w\in \Bbb C$ and let $\psi$ be an angle. We introduce the sets
\begin{equation}\label{nooo1}
\mathscr S'_w(e^{i\psi}\Bbb R+\zeta;X)=\{\mathcal F\in\mathscr
D'(e^{i\psi}\Bbb R+\zeta;X): e_w\mathcal F\in \mathscr
S'(e^{i\psi}\Bbb R+\zeta;X)\},
\end{equation}
\begin{equation}\label{nooo2}
\mathscr S'_\zeta(e^{-i\psi}\Bbb R+w;X)=\{ F\in\mathscr
D'(e^{-i\psi}\Bbb R+w;X): e_\zeta F\in \mathscr S'(e^{-i\psi}\Bbb
R+w;X)\},
\end{equation}
where $e_w$ and $e_\zeta$ denote the exponential weight functions
$$
e_w:\lambda\mapsto\exp(iw\lambda),\ \lambda\in e^{i\psi}\Bbb
R+\zeta;\quad e_\zeta:z\mapsto\exp(-i\zeta z),\ z\in
e^{-i\psi}\Bbb R+w.
$$
The  transformation \eqref{noTsp} can be continuously
 extended to all $\mathcal F\in \mathscr
S'_w(e^{i\psi}\Bbb R+\zeta;X)$ by the formula
\begin{equation}\label{no01}
\Bbb T_{\zeta,w}^\psi\mathcal F=e^{i\psi}e_{-\zeta}(\mathscr
F_{\xi\to t}[(e_w\mathcal F)\circ\varkappa_{\psi,\zeta}])\circ
\tau^{-1}_{\psi,w},
\end{equation}
where $\varkappa_{\psi,\zeta}:\Bbb R\to e^{i\psi}\Bbb R+\zeta$ and
$\tau^{-1}_{\psi,w}:e^{-i\psi}\Bbb R+w\to \Bbb R$ are the linear
transformations $\varkappa_{\psi,\zeta}(\xi)=e^{i\psi}\xi+\zeta$,
 $\tau^{-1}_{\psi,w}(z)=e^{i\psi}(z-w)$. The Fourier-Laplace transformation \eqref{no01} implements an
isomorphism
\begin{equation}\label{nooo3}
\Bbb T_{\zeta,w}^\psi:\mathscr S_w'(e^{i\psi}\Bbb
R+\zeta;X)\to\mathscr S_{\zeta}'(e^{-i\psi}\Bbb R+w;X).
\end{equation}
We have the following differentiation rule
\begin{equation}\label{rule} D_{(\psi)}^j \Bbb T_{\zeta,w}^\psi\mathcal
F=\Bbb T_{\zeta,w}^\psi \lambda^j \mathcal F, \quad \mathcal F\in
\mathscr S_w'(e^{i\psi}\Bbb R+\zeta;X);
\end{equation}
here the differential operator $D_{(\psi)}$ is the same as in
\eqref{D}.
\end{lem}

\begin{pf}
 The elements $\mathcal F\in\mathscr S'_w(e^{i\psi}\Bbb
R+\zeta;X)$ can be naturally identified with the elements
$\mathcal G\in\mathscr S'(\Bbb R;X)$ by the rule $\mathcal
G=(e_w\mathcal F)\circ\varkappa_{\psi,\zeta}$. We also identify
the elements $F\in\mathscr S'_\zeta(e^{-i\psi}\Bbb R+w;X)$ and
$G\in\mathscr S'(\Bbb R;X)$ by the rule $F=e^{i\psi}e_{-\zeta}
G\circ \tau_{\psi,w}^{-1}$. We set $G=\mathscr F_{\xi\to
t}\mathcal G$ and arrive at \eqref{no01}. The transformation
\eqref{nooo3} yields an isomorphism  because $\mathscr F_{\xi\to
t}:\mathscr S'(\Bbb R;X)\to \mathscr S'(\Bbb R;X)$ is an
isomorphism. Since the Fourier transformation in $\mathscr S'(\Bbb
R;X)$  is the continuous extension of the integral transformation
\eqref{ft} and $\mathscr W^\ell_w(e^{i\psi}\Bbb
R+\zeta;X)\subset\mathscr S'(e^{i\psi}\Bbb R+\zeta;X)$ we conclude
that~\eqref{no01} is the continuous extension of the
Fourier-Laplace transformation
 \eqref{noTsp}; cf. the proof of Lemma~\ref{lno1}. The
differentiation rule \eqref{rule} follows from \eqref{nooo4},
\eqref{D}. \qed\end{pf}

\subsection{Weighted Hardy-Sobolev spaces in cones}\label{ss3.2}
Let $\zeta,w\in\Bbb C$, $\ell\in\Bbb R$, and let $\psi$ be an
angle. We introduce the space $\mathsf
W^\ell_{\zeta}(e^{-i\psi}\Bbb R+w;X)$ of all distributions $F\in
\mathscr S'_{\zeta}(e^{-i\psi}\Bbb R+w;X)$ representable in the
form $F=\Bbb T^\psi_{\zeta,w}\mathcal F$ with some $ \mathcal
F\in\mathscr W^\ell_w(e^{i\psi}\Bbb R+\zeta;X)$; here $\mathbb
T^\psi_{\zeta,w}$ is given in \eqref{no01}, and $\mathscr
W^\ell_w(e^{i\psi}\Bbb R+\zeta;X)\subset \mathscr
S'_w(e^{i\psi}\Bbb R+\zeta;X)$, see Lemma~\ref{lno2}. Let the
space $\mathsf W^\ell_{\zeta}(e^{-i\psi}\Bbb R+w;X)$ be endowed
with the norm
\begin{equation}\label{nooo}
\|F;\mathsf W^\ell_{\zeta}(e^{-i\psi}\Bbb R+w;X)\|=\| (\mathbb
T^\psi_{\zeta,w})^{-1}F; \mathscr W^\ell_w(e^{i\psi}\Bbb
R+\zeta;X) \|.
\end{equation}
If $\ell=0$ then the equality \eqref{nooo} coincides with the
Parseval equality~\eqref{PE}. Using the formula \eqref{no01} and
the norm \eqref{NW} of the space $\mathscr W^\ell_w(e^{i\psi}\Bbb
R+\zeta;X)$, one can see that the norm \eqref{nooo} is equivalent
to the norm
\begin{equation}\label{NOS}
|\mspace{-2mu}|\mspace{-2mu}|F;\mathsf
W^\ell_{\zeta}(e^{-i\psi}\Bbb
R+w;X)|\mspace{-2mu}|\mspace{-2mu}|^2=\int_{\Bbb R}
(1+|\xi|^2)^{\ell}\|\mathscr F^{-1}_{t\to\xi}\{(e_\zeta
F)\circ\tau_{\psi,w}\}\|^2\,d\xi;
\end{equation}
moreover,  the constants $c_1$ and $c_2$ in the inequalities
\begin{equation}\label{uneq}
\begin{aligned}
\|F;\mathsf W^\ell_{\zeta}(e^{-i\psi}\Bbb R+w;X)\|\leq c_1
|\mspace{-2mu}|\mspace{-2mu}|F;\mathsf
W^\ell_{\zeta}(e^{-i\psi}\Bbb
R+w;X)|\mspace{-2mu}|\mspace{-2mu}|\\\leq c_2\|F;\mathsf
W^\ell_{\zeta}(e^{-i\psi}\Bbb R+w;X)\|
\end{aligned}
\end{equation} are
independent of $\psi$ and $w\in\Bbb C$.
 Therefore $\mathsf W^\ell_{\zeta}(e^{-i\psi}\Bbb R+w;X)$ is the
weighted Sobolev space of distributions with values in $X$, see
\cite{R4}. Note that in \eqref{NOS} the exponential weight
$$
e_\zeta\circ\tau_{\psi,w}(t)=\exp\{-i\zeta (e^{-i\psi}t+w)\},
\quad t\in\Bbb R,
$$
depends on the weight number $\zeta$ and the angle $\psi$.
 For all $v\in e^{-i\psi}\Bbb R+w$ we have
$$
\|F;\mathsf W^\ell_{\zeta}(e^{-i\psi}\Bbb R+w;X)\|=\|F;\mathsf
W^\ell_{\zeta}(e^{-i\psi}\Bbb R+v;X)\|.
$$
The space $\mathsf W^\ell_{\zeta}(e^{-i\psi}\Bbb R+w;X)$  does not
change while $\zeta$ travels along the line $e^{i\psi}\Bbb
R+\eta$, $\eta\in\Bbb C$, the norm changes for the equivalent one
$$
\|F;\mathsf W^\ell_{\zeta}(e^{-i\psi}\Bbb
R+w;X)\|=e^{\Im\{(\zeta-\eta)w\}}\|F;\mathsf
W^\ell_{\eta}(e^{-i\psi}\Bbb R+w;X)\|,\quad \zeta\in e^{i\psi}\Bbb
R+\eta.
$$
If $\ell\geq 0$ then
  the elements of the
 space $\mathsf W^\ell_{\zeta}(e^{-i\psi}\Bbb R+w;X)$ are (classes of) functions.
 In this case we can define an equivalent norm in $\mathsf W^\ell_{\zeta}(e^{-i\psi}\Bbb
 R+w;X)$ by the equality
\begin{equation}\label{snorm}
\begin{aligned}
&|\mspace{-2mu}|\mspace{-2mu}|F;\mathsf
W^\ell_{\zeta}(e^{-i\psi}\Bbb
R+w;X)|\mspace{-2mu}|\mspace{-2mu}|^2=\int\limits_{e^{-i\psi}\Bbb
R+w}\sum_{j\leq [\ell]}\|e^{-i\zeta
z}D^j_{(\psi)}F(z)\|^2\,|dz| \\
&+ \int\limits_{e^{-i\psi}\Bbb R+w}\int\limits_{e^{-i\psi}\Bbb
R+w}|z-u|^{-2(\ell-[\ell])-1} \|e^{-i\zeta
z}D^{[\ell]}_{(\psi)}F(z)-e^{-i\zeta u}D^{[\ell]}_{(\psi)}F(u)\|^2
\,|dz|\,|du|,
\end{aligned}
\end{equation}
 where $[\ell]$ is the integer part of $\ell$; see
 e.g.~\cite{R13,R14}. If $\ell$ is a nonnegative
  integer then from the rule \eqref{rule} and the
   Parseval equality \eqref{PE} it is easily seen that the
 norm \eqref{nooo} is equivalent to the norm
\begin{equation}\label{nw}
|\mspace{-2mu}|\mspace{-2mu}|F;\mathsf
W^\ell_{\zeta}(e^{-i\psi}\Bbb
R+w;X)|\mspace{-2mu}|\mspace{-2mu}|^2=\int\limits_{e^{-i\psi}\Bbb
R+w}\sum_{j\leq\ell}\| e^{-i\zeta z} D_{(\psi)}^jF(z)\|^2\,|d z|.
\end{equation}
The differentiation in~\eqref{snorm},~\eqref{nw} is understood in
the sense of distributions. For the equivalent
norms~\eqref{snorm},~\eqref{nw} the inequalities \eqref{uneq}
remain valid with some constants $c_1$ and $c_2$ independent of
$w$, $\psi$ (and $\mathcal F$).
 We have
the following corollary of Lemma~\ref{lno2}.
\begin{cor}\label{noCol} The transformation \eqref{no01} realizes
an isometric
 isomorphism
\begin{equation}\label{Tww}
\mathbb T^\psi_{\zeta,w}:  \mathscr W^\ell_w(e^{i\psi}\Bbb
R+\zeta;X)\to \mathsf W^\ell_{\zeta}(e^{-i\psi}\Bbb R+w;X).
\end{equation}
\end{cor}

Let $K^\varphi_w$ denote the open  double-napped cone
$$
K^\varphi_w=\{z\in\Bbb C:z=e^{-i\psi}t+w, t\in \Bbb
R\setminus\{0\}, 0<\psi<\varphi\}
$$
with the vertex $w$ and the angle $\varphi\in(0,\pi]$.
\begin{defn}\label{defH}
Let $w,\zeta\in \Bbb C$, $\ell\geq 0$, and $\varphi\in(0,\pi]$. We
introduce the Hardy-Sobolev class $\mathsf
H^\ell_\zeta(K^\varphi_w;X)$ of order $\ell$ as the set of all
analytic functions $K^\varphi_w\ni z\mapsto F(z)\in X$ satisfying
the uniform in $\psi\in(0,\varphi)$ estimate
\begin{equation}
\|F;\mathsf W^\ell_{\zeta}(e^{-i\psi}\Bbb R+w;X)\|\leq C(F).
\end{equation}
\end{defn}
\begin{thm} \label{p3} Let $\ell\geq 0$, $\varphi\in(0,\pi]$, and $w,\zeta\in\Bbb C$. The following
assertions are fulfilled.

(i) Every function $F\in \mathsf H^\ell_\zeta(K^\varphi_w;X)$ has
 boundary limits $F_0\in \mathsf W^\ell_\zeta(\Bbb
R+w;X)$ and $F_\varphi\in \mathsf W^\ell_\zeta(e^{-i\varphi}\Bbb
R+w;X)$ in the sense that for almost all points $u$ of the
boundary $\partial K^\varphi_w$ we have $\|F(z)-F_0(u)\|\to 0$ as
$z$ non-tangentially tends to $u\in\Bbb R+w$, and
$\|F(z)-F(u)\|\to 0$ as $z$ non-tangentially tends to $u\in
e^{-i\varphi}\Bbb R+w$. Moreover,
\begin{equation}\label{fas}
\begin{aligned}
&\|(e_\zeta F)\circ\tau_{\psi,w}-(e_\zeta
F_0)\circ\tau_{0,w};\mathsf W_0^\ell(\Bbb R;X)\|\to 0,\quad\psi\to
0+,
\\
&\|(e_\zeta F)\circ\tau_{\psi,w}-(e_\zeta
F_\varphi)\circ\tau_{\varphi,w};\mathsf W_0^\ell(\Bbb R;X)\|\to
0,\quad\psi\to \varphi-;
\end{aligned}
\end{equation}
 here the linear transformation $\tau_{\varphi,w}$ and the
exponential weight function $e_\zeta$ are the same as in
Lemma~\ref{lno2}. From now on we suppose that the elements $F\in
\mathsf H^\ell_\zeta(K^\varphi_w;X)$ are extended to the boundary
$\partial K^\varphi_w$ by non-tangential limits. In the case
$\varphi=\pi$ we distinguish the banks $\lim_{\psi\to
0+}(e^{-i\psi}\Bbb R+w)$ and $\lim_{\psi\to\pi-}(e^{-i\psi}\Bbb
R+w)$  in $\partial K^\pi_w$.

(ii)  For all $F\in\mathsf H^\ell_\zeta(K^\varphi_w;X)$ and
$\psi\in[0,\varphi]$ the estimate
\begin{equation}\label{estim}
\|F;\mathsf W^\ell_{\zeta}(e^{-i\psi}\Bbb R+w;X)\|\leq
C\bigl(\|F;\mathsf W^\ell_{\zeta}(\Bbb R+w;X)\|+\|F;\mathsf
W^\ell_{\zeta}(e^{-i\varphi}\Bbb R+w;X)\|\bigr)
\end{equation}
holds, where the constant $C$ is independent of $F$, $\psi$, and
$w$.

(iii) The Hardy-Sobolev class $\mathsf
H^\ell_\zeta(K^\varphi_w;X)$ endowed with the norm
\begin{equation}\label{NormH}
\|F;\mathsf H^\ell_\zeta(K^\varphi_w;X)\|=\|F;\mathsf
W^\ell_{\zeta}(\Bbb R+w;X)\|+\|F;\mathsf
W^\ell_{\zeta}(e^{-i\varphi}\Bbb R+w;X)\|
\end{equation}
is a Banach space.

 (iv) Let us identify $F\in\mathsf H^\ell_\zeta(K^\varphi_w;X)$ with the set
$\{F_\psi:\psi\in[0,\varphi]\}$ of functions $F_\psi\in \mathsf
W^\ell_\zeta(e^{-i\psi}\Bbb R+w;X)$, where $F_\psi$ is the
restriction  of $F$ to the line $e^{-i\psi}\Bbb R+w$. In the same
way we identify $\mathcal F\in\mathscr H^\ell_w(\mathcal
K^\varphi_\zeta;X)$ with the set $\{\mathcal
F_\psi:\psi\in[0,\varphi]\}$, where $\mathcal F_\psi=\mathcal
F\!\!\upharpoonright_{e^{i\psi}\Bbb R+\zeta}\in \mathscr
W^\ell_{w}(e^{i\psi}\Bbb R+w;X)$. The Fourier-Laplace
transformation
\begin{equation}\label{transf}
\begin{aligned}
\mathscr H^\ell_w(\mathcal K^\varphi_\zeta;X)\ni\mathcal
F\equiv&\{\mathcal F_\psi:\psi\in[0,\varphi]\}
\\
&\mapsto \{\Bbb T^\psi_{\zeta,w}\mathcal
F_\psi:\psi\in[0,\varphi]\}\equiv F\in \mathsf
H^\ell_{\zeta}(K^\varphi_w;X)
\end{aligned}
\end{equation}
yields an isometric isomorphism.

(v) For any $\psi\in[0,\varphi]$ the set of functions
$\{F\!\!\upharpoonright_{e^{-i\psi}\Bbb R+w}:F\in \mathsf
H^\ell_\zeta(K^\varphi_w;X)\}$ is dense in the Sobolev space
$\mathsf W^\ell_\zeta(e^{-i\psi}\Bbb R+w;X)$.
\end{thm}

The proof of Theorem~\ref{p3} is preceded by the lemma.
\begin{lem}\label{l4}
 Let  $F\in \mathsf H^0_0(K_0^\varphi;X)$, $\varphi\in (0,\pi]$.
  We denote $K^{\varphi,\pm}_0=\{z\in K^\varphi_0: \Im z\gtrless 0\}$.

(i) Suppose that $F(z)=0$ for  $z\in K^{\varphi,+}_0$. Then there
exists an analytic function
$$
\bigcup_{\psi\in(0,\varphi)} e^{i\psi}\Bbb
C^-\ni\lambda\mapsto\mathcal F(\lambda)\in X
$$ such
that
 for all $\psi\in(0,\varphi)$ the equality $\mathcal F=(\Bbb
T^\psi_{0,0})^{-1}( F\!\!\upharpoonright_{e^{-i\psi}\Bbb R})$
holds almost everywhere on the line $e^{i\psi}\Bbb R$.

(ii) Suppose that $F(z)=0$ for  $z\in K^{\varphi,-}_0$. Then there
exists an analytic function
$$
\bigcup_{\psi\in(0,\varphi)}
e^{i\psi}\Bbb C^+\ni\lambda\mapsto\mathcal F(\lambda)\in X
$$
 such
that
 for all $\psi\in(0,\varphi)$ the equality $\mathcal F=(\Bbb
T^\psi_{0,0})^{-1} (F\!\!\upharpoonright_{e^{-i\psi}\Bbb R})$
holds almost everywhere on the line $e^{i\psi}\Bbb R$.
\end{lem}
\begin{pf} We shall prove the assertion (i). The proof of the assertion (ii) is similar.
Let us  establish the equality
\begin{equation}\label{15}
\int_{e^{-i\psi}\Bbb R^+} e^{-i\lambda z} F (z) \,
dz=\int_{e^{-i\phi}\Bbb R^+} e^{-i\lambda z} F(z) \, dz, \quad
\lambda\in e^{i\psi}\Bbb C^-\cap e^{i\phi}\Bbb C^-,
\end{equation}
where $0<\psi<\phi<\varphi$, the integrations run from $0$ to
$e^{-i\psi}(+\infty)$ and from $0$ to $e^{-i\phi}(+\infty)$. Due
to the inclusion $F\in \mathsf H^0_0(K^\varphi_0;X)$ we have $F\in
\mathsf W^0_0(e^{-i\psi}\Bbb R;X)$ and $F\in\mathsf
W^0_0(e^{-i\phi}\Bbb R;X)$. The integrals in \eqref{15} are
absolutely convergent in $X$. Since $F$ is analytic in
$K^\varphi_0$, we get
\begin{equation}\label{oint}
\oint_{\mathscr C(a,\psi,\phi)} e^{-i\lambda z} F(z)\,dz=0,
\end{equation}
where the integration runs along the closed contour
\begin{equation*}
\begin{aligned}
\mathscr C(a,\psi,\phi)=\{z:z=a
e^{-i\vartheta},\vartheta\in[\psi,\phi]\}\cup\{z:z=te^{-i\phi},1/a\leq
t\leq a\}\\\cup\{z:z=
e^{-i\vartheta}/a,\vartheta\in[\psi,\phi]\}\cup\{z:z=t e^{-i\psi},
1/a\leq t\leq a\},\ a>0.
\end{aligned}
\end{equation*}
Let us note that the value $a^{-1/2}\|F(e^{-i\vartheta}a)\|$ tends
to zero  uniformly in $\vartheta\in [\psi,\phi]$ as $a\to 0$ or
$a\to+\infty$. Indeed, this follows from the assertion 1.(iv) of
Proposition~\ref{p2} since the function $\mathcal
K^\varphi_0\ni\lambda\mapsto F(e^{-i\varphi}\lambda)\in X$ is in
the class $\mathscr H^0_0(\mathcal K^\varphi_0;X)$. Let us also
note that for any $\lambda\in e^{i\psi}\Bbb C^-\cap e^{i\phi}\Bbb
C^-$ the exponent $\exp(-i\lambda a e^{-i\vartheta})$ tends to
zero uniformly in $\vartheta\in[\psi,\varphi]$ as $a\to+\infty$,
and the value $|\exp(-i\lambda  e^{-i\vartheta}/a)|$ remains
uniformly bounded as $a\to+\infty$, $\vartheta\in[\psi,\varphi]$.
Passing in \eqref{oint} to the limit as $a\to+\infty$ we arrive at
\eqref{15}.

It is clear that the integral
\begin{equation}\label{tr}
 \mathcal
F(\lambda)=\frac 1{\sqrt {2\pi}}\int_{e^{-i\psi}\Bbb R^+}
e^{-i\lambda z} F (z) \, dz
\end{equation}
defines an analytic function $e^{i\psi}\Bbb
C^-\ni\lambda\mapsto\mathcal F(\lambda)\in X$. Due to \eqref{15},
where $\psi$ comes arbitrarily close to zero and $\phi$ comes
arbitrarily close to $\varphi$, the function $\mathcal F$ has an
analytic continuation to the set $\cup_{\psi\in(0,\varphi)}
e^{i\psi}\Bbb C^-$. It remains to show that for all
$\psi\in(0,\varphi)$ the equality $\mathcal F=(\Bbb
T^\psi_{0,0})^{-1}( F\!\!\upharpoonright_{e^{-i\psi}\Bbb R})$ is
valid almost everywhere on the line $e^{i\psi}\Bbb R$.

Recall that $F(z)=0$ for $z\in K^{\varphi,+}_0$. Let $\zeta\in
e^{i\psi}\Bbb C^-$. Then we can rewrite the equality \eqref{tr},
where $\lambda\in e^{i\psi}\Bbb R+\zeta$, in the form $\mathcal
F=(\Bbb T^\psi_{\zeta,0})^{-1}
(F\!\!\upharpoonright_{e^{-i\psi}\Bbb R})$ or, equivalently, in
the form $\mathcal F(\cdot+\zeta)=(\Bbb T^\psi_{0,0})^{-1}
(e_\zeta F\!\!\upharpoonright_{e^{-i\psi}\Bbb R})$; cf.
Lemma~\ref{lno1} and Lemma~\ref{lno2}. From Corollary~\ref{noCol}
we have
\begin{equation}\label{n}
\|\mathcal F(\cdot+\zeta)-(\Bbb T^\psi_{0,0})^{-1}
(F\!\!\upharpoonright_{e^{-i\psi}\Bbb R});\mathscr
W^0_0(e^{i\psi}\Bbb R;X)\|=\|e_\zeta F-F;\mathsf
W_0^0(e^{-i\psi}\Bbb R;X)\|.
\end{equation}
The right hand side of the equality \eqref{n} tends to zero as
$\zeta\to 0$, $\zeta\in e^{i\psi}\Bbb C^-$. This is
 because $F(z)=0$
for $z\in e^{-i\psi}\Bbb R^-$, the exponential function $e_\zeta$
is bounded on $e^{-i\psi}\Bbb R^+$ for $\zeta\in e^{i\psi}\Bbb
C^-$, and $e_\zeta\to1$ as $\zeta\to 0$. Thus $\mathcal F(\cdot
+\zeta)$ tends to $(\Bbb T^\psi_{0,0})^{-1}
(F\!\!\upharpoonright_{e^{-i\psi}\Bbb R})$ in $\mathscr
W^0_0(e^{i\psi}\Bbb R;X)$ as $\zeta\to 0$, $\zeta\in e^{i\psi}\Bbb
C^-$, and $\mathcal F=(\Bbb T^\psi_{0,0})^{-1}
(F\!\!\upharpoonright_{e^{-i\psi}\Bbb R})$
 almost everywhere on $e^{i\psi}\Bbb R$.
 \qed\end{pf}

\begin{pf*}{\bf PROOF of Theorem~\ref{p3}.} We first prove the
proposition for the case $\zeta=w=0$.  Let us prove that for every
$F\in\mathsf H^\ell_0( K^\varphi_0;X)$, $\ell\geq 0$, there exists
 $\mathcal F\in \mathscr H^\ell_0(\mathcal K^\varphi_0;X)$ such
that $\mathcal F\!\!\upharpoonright_{e^{i\psi}\Bbb R}=(\Bbb
T^\psi_{0,0})^{-1}(F\!\!\upharpoonright_{e^{-i\psi}\Bbb R})$ for
all $\psi\in(0,\varphi)$. Since $\ell\geq  0$, the inclusion
$\mathscr H^\ell_0(\mathcal K^\varphi_0;X)\subseteqq \mathscr
H^0_0(\mathcal K^\varphi_0;X)$ holds and Lemma~\ref{l4} can be
applied. We set $F^\pm(z)=F(z)$ for $z\in K^{\varphi,\pm}_0$ and
$F^\pm(z)=0$ for $z\in K^{\varphi,\mp}_0$.  By Lemma~\ref{l4}
there exists an analytic function $\mathcal
K^\varphi_0\ni\lambda\mapsto \mathcal F(\lambda)\in X$ such that
for all $\psi\in(0,\varphi)$  the equality
$$
\mathcal F=(\Bbb
T_{0,0}^\psi)^{-1}(F^+\!\!\upharpoonright_{e^{-i\psi}\Bbb
R}+F^-\!\!\upharpoonright_{e^{-i\psi}\Bbb R})=(\Bbb
T_{0,0}^\psi)^{-1}(F\!\!\upharpoonright_{e^{-i\psi}\Bbb R})
$$
 holds almost everywhere on the line ${e^{i\psi}\Bbb R}$.   By Definition~\ref{defH} and
Corollary~\ref{noCol} we have
$$
\|\mathcal F;\mathscr W^\ell_0(e^{i\psi}\Bbb R;X)\|=\|F;\mathsf
W^\ell_0(e^{-i\psi}\Bbb R;X)\|\leq Const, \quad
\psi\in(0,\varphi).
$$
This proves the inclusion $\mathcal F\in\mathscr H^\ell_0(\mathcal
K^\varphi_0;X)$, see Definition~\ref{d1}. In a similar way one can
see that for every $\mathcal F\in\mathscr H^\ell_0(\mathcal
K^\varphi_0;X)$ there exists $F\in \mathsf
H^\ell_0(K^\varphi_0;X)$ such that
$F\!\!\upharpoonright_{e^{-i\psi}\Bbb R}=\Bbb
T^\psi_{0,0}(\mathcal F\!\!\upharpoonright_{e^{i\psi}\Bbb R})$, we
do not cite the proof of this implication.

We set $F_0=\Bbb T_{0,0}^0(\mathcal F\!\!\upharpoonright_{\Bbb
R})$ and $F_\varphi=\Bbb T^\varphi_{0,0}(\mathcal
F\!\!\upharpoonright_{e^{i\varphi}\Bbb R})$. The assertion 1. (i)
of Proposition~\ref{p2} together with Corollary~\ref{noCol} leads
to the relations~\eqref{fas}. It is clear that $F_0$ and
$F_\varphi$ coincide almost everywhere with non-tangential
boundary limits of $F\in\mathsf H^\ell_0(K^\varphi_0;X)$. Indeed,
the embedding $\mathsf W^\ell_0(\Bbb R;X)\subseteq \mathsf
W^0_0(\Bbb R;X)$ is continuous, and for $\ell=0$ the assertion (i)
can be established  by the same argument as we used in the proof
Proposition~\ref{p2}.  The assertion (i) is proved.

The inequality \eqref{no2} provides the estimate \eqref{estim} for
the assertion (ii). If we define the norm in $\mathsf
H^\ell_0(K^\varphi_0;X)$ by the equality \eqref{NormH} then we get
$$
\|F;\mathsf H^\ell_0(K^\varphi_0;X)\|=\|\mathcal F;\mathscr
W^\ell_0(\Bbb R;X)\|+\|\mathcal F;\mathscr
W^\ell_0(e^{i\varphi}\Bbb  R;X)\|=\|\mathcal F;\mathscr
H^\ell_0(\mathcal K^\varphi_0;X)\|.
$$
This equality finishes the proof of the assertion (iv).  The space
$\mathsf H^\ell_0(K^\varphi_0;X)$ is complete because it is
isomorphic to the complete space $\mathscr H^\ell_0(\mathcal
K^\varphi_0;X)$. The assertion (iii)  follows from (iv) and
Proposition~\ref{p3`}. Finally, the assertion (v) is a consequence
of (iv) together with Corollary~\ref{noCol} and
Proposition~\ref{p2}.1,(v). The proposition is proved for the case
$\zeta=w=0$.

Let us consider the general case. We identify the classes $\mathsf
H^\ell_\zeta (K_w^\varphi;X)$ and $ \mathsf H^\ell_0(
K^\varphi_0;X)$ by setting $F(z-w)=exp\{-i\zeta z\}G(z)$, $z\in
K_w^\varphi$, where $G\in \mathsf H^\ell_\zeta(K_w^\varphi;X)$ and
$F\in \mathsf H^\ell_{0} (K^\varphi_0;X)$.  Let us also identify
the classes $\mathscr H^\ell_0(\mathcal K_0^\varphi;X)$ and
$\mathscr H^\ell_w(\mathcal K_\zeta^\varphi;X)$ according to  the
rule $ \mathcal G(\lambda-\zeta)= e^{-i w\lambda}\mathcal
F(\lambda)$, $\lambda\in \mathcal K^\varphi_\zeta $. This allows
us to reformulate the results proved for $F\in \mathsf H^\ell_0(
K_0^\varphi; X)$ and $\mathcal F\in \mathscr H^\ell_0(\mathcal
K_0^\varphi;X)$ in terms of $G\in \mathsf H^\ell_\zeta
(K_w^\varphi;X)$ and $\mathcal G\in \mathscr H^\ell_w(\mathcal
K_\zeta^\varphi;X)$. \qed\end{pf*}

 If $\ell\in
\Bbb N$ then the space $\mathsf H^\ell_\zeta(K^\varphi_w;X)$
consists of all elements $F\in \mathsf H^0_\zeta(K^\varphi_w;X)$
such that the complex derivatives $K^\varphi_w\ni z\mapsto
D_z^jF(z)\in X$, $j=1,\dots,\ell$, are also in the class $\mathsf
H^0_\zeta(K^\varphi_w;X)$,  the norm \eqref{NormH} in $\mathsf
H^\ell_\zeta(K_w^\varphi;X)$ is equivalent to the norm
$$
|\mkern -2mu|\mkern -2mu|F;\mathsf
H^\ell_\zeta(K_w^\varphi;X)|\mkern -2mu|\mkern -2mu|=\| F; \mathsf
H^0_\zeta ( K_w^\varphi;X)\|+\|D_z^\ell F;\mathsf H^0_\zeta (
K_w^\varphi;X)\|.
$$

 The next theorem presents an elementary embedding result.

\begin{thm}\label{THMemb}If  $m\in\Bbb N$ and $\ell>m-1/2$ then
 $ \mathsf H^\ell_\zeta(K^\varphi_w;X)\subset
 C^{m-1}\bigl(\overline{K^{\varphi}_w};X\bigr)$.
\end{thm}

\begin{pf}
Let $G\in \mathsf H^\ell_\zeta (K^\varphi_w;X)$. We  set
$F(z)=e^{-i\zeta z} G(z+w)$, $z\in K^\varphi_0$. Then $F\in
\mathsf H^\ell_0(K^\varphi_0;X)$. It suffices to show that
$F(z)\in C^{m-1}(\overline{K^\varphi_0};X)$. By the Sobolev
embedding theorem we immediately get $F\in\mathsf
W^\ell_0(e^{-i\psi}\Bbb R;X)\subset C^{m-1}(e^{-i\psi}\Bbb R;X)$,
$\psi\in[0,\varphi]$.
 Let $\mathcal F\in
\mathscr H^\ell_0(\mathcal K^\varphi_0;X)$ be the transform
\eqref{transf} of $F$. Since $\ell>m-1/2>0$, the analytic
functions $\mathcal K^\varphi_0\ni\lambda\mapsto \lambda^j\mathcal
F(\lambda)\in X$, $j=0,1,\dots,m-1$, are in the class $\mathscr
H^0_0(\mathcal K_0^\varphi;X)$. For $z\in e^{-i\psi}\Bbb R$,
$\psi\in[0,\varphi]$, we have
\begin{equation}\label{e1}
\begin{aligned}
D_{(\psi)}^j F(z)=\frac 1 {\sqrt{2\pi}}\int_{e^{i\psi}\Bbb R}
e^{i\lambda z} \lambda^j\mathcal F(\lambda)\, d\lambda,\quad
j=0,1,\dots,m-1,
\end{aligned}
\end{equation}
where the integral is absolutely convergent in $X$. Let us recall
that $D_{(\psi)}^j F(z)=D^j_z F(z)$ for $z\in e^{-i\psi}\Bbb
R\setminus\{0\}$, $\psi\in(0,\varphi)$; here $D_z=-\frac i
2(\partial_{\Re z}-i\partial_{\Im z})$. Using the same arguments
as in the proof of Lemma~\ref{l4} we arrive at the equality
\begin{equation}\label{e2}
\int_{e^{i\psi}\Bbb R^+} e^{i\lambda z} \lambda^j\mathcal
F(\lambda)\,d\lambda=\int_{e^{i\phi}\Bbb R^+} e^{i\lambda z}
\lambda^j \mathcal F(\lambda)\,d\lambda,\ j=0,1,\dots,m-1,
\end{equation}
where $0<\psi<\phi<\varphi$ and $z\in e^{-i\psi}\overline{\Bbb
C^+}\cap e^{-i\phi}\overline{\Bbb C^+}$. The assertion 1.(i) of
Proposition~\ref{p2} and the inequality $\ell>j+1/2$ allow us to
pass in \eqref{e2} to the limits as $\psi\to 0+$ and
$\phi\to\varphi-$. As a result we get
\begin{equation}\label{e5}
\int_{e^{i\psi}\Bbb R^+} e^{i\lambda z} \lambda^j\mathcal
F(\lambda)\,d\lambda=\int_{e^{i\varphi}\Bbb R^+} e^{i\lambda z}
\lambda^j \mathcal F(\lambda)\,d\lambda,\ \psi\in[0,\varphi),\
z\in e^{-i\psi}\overline{\Bbb R^+}.
\end{equation}
In a similar way one can obtain the equality
\begin{equation}\label{e6}
\int_{e^{i\psi}\Bbb R^-} e^{i\lambda z} \lambda^j\mathcal
F(\lambda)\,d\lambda=\int_{\Bbb R^-} e^{i\lambda z} \lambda^j
\mathcal F(\lambda)\,d\lambda,\quad\psi\in(0,\varphi],\ z\in
e^{-i\psi}\overline{\Bbb R^-}.
\end{equation}
the formulas \eqref{e1}, \eqref{e5}, and \eqref{e6} give us the
equality
\begin{equation}\label{e3}
\begin{aligned}
D_{(\psi)}^j F(e^{-i\psi}t) =\frac 1
{\sqrt{2\pi}}&\quad\left(\int_{\Bbb R^-} \exp\{{i\lambda
e^{-i\psi}t}\} \lambda^j \mathcal F(\lambda)\,d\lambda\right.
\\&\left.+\int_{e^{i\varphi}\Bbb R^+} \exp\{i\lambda e^{-i\psi}t\}
\lambda^j \mathcal F(\lambda)\,d\lambda\right)
\end{aligned}
\end{equation}
for $j=0,\dots,m-1$, $t\geq 0$ and $\psi\in[0,\varphi]$. Since
$\mathcal F\in \mathscr W^\ell_0(\Bbb R;X)$, $\mathcal
F\in\mathscr W^\ell_0(e^{i\varphi}\Bbb R; X)$, and $\ell>j+1/2$,
the integrals in \eqref{e3} are absolutely convergent in $X$
uniformly with respect to $t\geq 0$ and $\psi\in[0,\varphi]$. Thus
$F\in C^{m-1}\bigl(\overline{\Bbb C^-}\cap
e^{-i\varphi}\overline{\Bbb C^+};X\bigr)$. Similarly we obtain the
inclusion $F\in C^{m-1}\bigl(\overline{\Bbb C^+}\cap
e^{-i\varphi}\overline{\Bbb C^-};X\bigr)$. Summing up we can say
that $F\in C^{m-1}(\overline{K^{\varphi}_0};X)$. \qed\end{pf}

Let $f\in\mathscr D'(\Bbb R\times (0,\varphi);X)$ and $v\in
C_0^\infty(\Bbb R)$. We define $f_{(v)}\in\mathscr
D'((0,\varphi);X)$ by the formula $f_{(v)}(u)=f(uv)$, where $u\in
C_0^\infty(0,\varphi)$. We say that a distribution $f\in\mathscr
D'(\Bbb R\times (0,\varphi);X)$ is continuous with respect to
$\psi$ and write $f\in C((0,\varphi);X)_\psi$ if for any $v\in
C_0^\infty(\Bbb R)$ there exists a function $\widetilde f_{(v)}\in
C((0,\varphi);X)$ such that $$f_{(v)} (u)=\int_0^\varphi
\widetilde f_{(v)}(\psi)\overline{u(\psi)}\,d\psi,\quad\forall
u\in C_0^\infty(0,\varphi).$$ If $f\in \mathscr D'(\Bbb R\times
(0,\varphi);X)\cap C((0,\varphi);X)_\psi$ then the restrictions
$f_\psi\in \mathscr D'(\Bbb R;X)$ of  $f$ are well-defined for all
$\psi\in(0,\varphi)$, and
$$
\widetilde f_{(v)}(\psi)=f_\psi(v), \quad  v\in C_0^\infty(\Bbb
R),\ \psi\in(0,\varphi);
$$
an additional point to emphasize is that
$$
f(\chi)= \int_0^\varphi
f_\psi\bigl(\chi(\cdot,\psi)\bigr)\,d\psi,\quad\forall\chi\in
C^\infty_0(\Bbb R\times (0,\varphi)),
$$
see e.g. \cite{R18}. For any $G\in \mathsf
H^0_\zeta(K^\varphi_w;X)$ we have
$(D^j_{(\psi)}G)\circ\tau_{\psi,w}\in C((0,\varphi);X)_\psi$
because
$$
\|(e_\zeta D^j_{(\psi)}G)\circ\tau_{\psi,w}-(e_\zeta
D^j_{(\phi)}G)\circ\tau_{\phi,w};\mathsf W^{-j}_0(\Bbb R;X)\|\to
0,\  \psi\to\phi,\ j=0,1,2,\dots,
$$
where $\psi$ and $\phi$ are in the interval $[0,\varphi]$;
cf.~\eqref{fas}. Therefore if $G\in \mathsf
H^0_\zeta(K^\varphi_w;X)$ then to the set of distributions
$$\{F_\psi\in \mathsf W^{-j}_0(e^{-i\psi}\Bbb R+w;X):
F_\psi=D^j_{(\psi)}(G\!\!\upharpoonright_{e^{-i\psi}\Bbb R+w}),
\psi\in(0,\varphi)\}$$ there corresponds a unique distribution
$f\in\mathscr D'(\Bbb R\times(0,\varphi);X)\cap
C((0,\varphi);X)_\psi$ such that $ f(\chi)=\int_0^\varphi
F_\psi(\chi(\tau^{-1}_{\psi,w},\psi)\, d\psi$. We are now in
position to consider Hardy-Sobolev classes of negative orders. We
start with the following
 general definition.

\begin{defn}\label{d2} Let $w,\zeta\in \Bbb C$, $\varphi\in(0,\pi]$ and let $\ell\in\Bbb
R$. We introduce the Hardy-Sobolev class $\mathsf H^\ell_\zeta(K^\varphi_w;X)$ 
as the collection
 of all sets of
distributions $
 F\equiv\{F_\psi\in\mathsf W^\ell_\zeta(e^{-i\psi}\Bbb R+w;X):\psi\in(0,\varphi)\}
$ satisfying the conditions: for $\{F_\psi:\psi\in(0,\varphi)\}$
there exist some constant $C(\mathcal F)$, a finite $m\in\Bbb N$,
and a set of functions $\{G_j\in\mathsf
H^0_\zeta(K^\varphi_w;X)\}_{j=0}^m$ such that the estimate and the
representation
\begin{equation}\label{er}
\|F_\psi;\mathsf W^\ell_{\zeta}(e^{-i\psi}\Bbb R+w;X)\|\leq C(
F),\quad F_\psi=\sum_{j=0}^m
D^j_{(\psi)}(G_j\!\!\upharpoonright_{e^{-i\psi}\Bbb R+w})
\end{equation}
 are valid for all $\psi\in(0,\varphi)$; here $D_{(\psi)}$
is the same as in \eqref{D}.
\end{defn}

Since the functions $G_j\in \mathsf H^0_\zeta(K^\varphi_w;X) $ are
analytic in $K^\varphi_w$  we have
\begin{equation}\label{anal}
\bigl(D^j_{(\psi)}(G_j\!\!\upharpoonright_{e^{-i\psi}\Bbb
R+w})\bigr)(z)=D^j_z G_j(z), \  z\in (e^{-i\psi}\Bbb
R+w)\setminus\{w\}, \psi\in(0,\varphi),\
\end{equation}
where $D_z=-\frac i 2(\partial_{\Re z}-i\partial_{\Im z})$ is the
complex derivative.  From \eqref{er} and \eqref{anal} it is easily
seen that  the distribution $F_\psi$ coincides with the analytic
function $\sum D^j_z G_j$ on the set $(e^{-i\psi}\Bbb
R+w)\setminus\{w\}$. If $\{F_\psi:\psi\in(0,\varphi)\}\in \mathsf
H^\ell_\zeta(K^\varphi_w;X)$ then the singular support of $F_\psi$
is empty or it consists of the only vertex $\{w\}$ of the cone
$K^\varphi_w$.

\begin{lem}\label{lno3}
In the case $\ell\geq 0$ Definition~\ref{d2} and
Definition~\ref{defH} are equivalent. If $F$ is a  function from
the class $\mathsf H^\ell_\zeta(K^\varphi_w;X)$ by
Definition~\ref{defH} then the set
\begin{equation}\label{anf1}
 \{F_\psi=F\!\!\upharpoonright_{e^{-i\psi}\Bbb
R+w}\in\mathsf W^\ell_\zeta(e^{-i\psi}\Bbb
R+w;X):\psi\in(0,\varphi)\}
\end{equation}
is in the class $\mathsf H^\ell_\zeta(K^\varphi_w;X)$ by
Definition~\ref{d2}. If $\{F_\psi:\psi\in(0,\varphi)\}\in\mathsf
H^\ell_\zeta(K^\varphi_w;X)$ by Definition~\ref{d2} then in
accordance with Definition~\ref{defH} the class $\mathsf
H^\ell_\zeta(K^\varphi_w;X)$ contains the analytic in $
K^\varphi_w$ function $F$ given by the equality
\begin{equation}\label{anf}
F(z)=F_\psi(z),\quad z\in K^\varphi_w\cap(e^{-i\psi}\Bbb R+w),\
\psi\in(0,\varphi).
\end{equation}
\end{lem}
\begin{pf} If $F\in\mathsf H^\ell_\zeta(K^\varphi_w;X)$ by
Definition~\ref{defH} then we can identify the analytic in
$K^\varphi_w$ function $F$ with the set \eqref{anf1}.
 Obviously, this set is in the
class $\mathsf H^\ell_\zeta(K^\varphi_w;X)$ according to
Definition~\ref{d2}.

Let $\ell\geq 0$ and $\{F_\psi:\psi\in(0,\varphi)\}\in\mathsf
H^\ell_\zeta(K^\varphi_w;X)$ by Definition~\ref{d2}. Then from the
representation in \eqref{er} it follows that the function $F$
given by \eqref{anf} is analytic in the cone $K^\varphi_w$; see
also \eqref{anal}. Taking into account the estimate from
\eqref{er}, we see that $F\in \mathsf H^\ell_\zeta(K^\varphi_w;X)$
by Definition~\ref{defH}. \qed\end{pf}

 Every set of distributions
$\{F_\psi:\psi\in(0,\varphi)\}\in \mathsf
H^\ell_\zeta(K^\varphi_w;X)$ defines an analytic function
$K^\varphi_w\ni z\mapsto F(z)\in X$
by the rule \eqref{anf}. 
 It is easy to see that  in the general case elements of the space
$\mathsf H^\ell_\zeta(K^\varphi_w;X)$ cannot be reconstructed from
the corresponding analytic functions. For instance,  the set of
distributions
 \begin{equation}\label{expl}
  \{F_\psi=e^{i
j\psi}(D_t^{j-1}\delta)\circ\tau^{-1}_{\psi,w}:\psi\in
(0,\varphi)\}\in\mathsf H^{-j}_\zeta (K_w^\varphi;\Bbb C),\quad
j\in\Bbb N,
\end{equation}
 defines the analytic in
$K^\varphi_w$ function $F\equiv 0$; here $\delta$ denotes the
Dirac delta function on the real axis.

In what follows we remain to use the notations
\begin{equation}\label{notations}(e_\zeta F)\circ\tau_{\psi,w},\ \ F\!\!\upharpoonright_{e^{-i\psi}\Bbb R+w},\
\
 \|F;\mathsf W^\ell_\zeta(e^{-i\psi}\Bbb R+w;X)\|,
 \end{equation}
 which
we employed studying the case $\ell\geq 0$. This can be done
without ambiguity if by the notations \eqref{notations} we shall
mean
 $$(e_\zeta F_\psi)\circ\tau_{\psi,w},\ \ F_\psi,\ \  \|F_\psi;\mathsf W^\ell_\zeta(e^{-i\psi}\Bbb R+w;X)\|.$$ The next
theorem in particular generalizes Theorem~\ref{p3} to the case
$\ell\in\Bbb R$.

\begin{thm}\label{p1.13} Let $\ell\in\Bbb R$, $\varphi\in(0,\pi]$, and $w,\zeta\in \Bbb
C$. The following assertions hold.

(i) Every set of distributions
$F\equiv\{F_\psi:\psi\in(0,\varphi)\}\in \mathsf
H^\ell_\zeta(K^\varphi_w;X)$ has  boundary limits $F_0\in \mathsf
W^\ell_\zeta(\Bbb R+w;X)$ and $F_\varphi\in \mathsf
W^\ell_\zeta(e^{-i\varphi}\Bbb R+w;X)$ such that the relations
\eqref{fas} hold. From now on we suppose that the elements
$F\in\mathsf H^\ell_\zeta(K^\varphi_w;X)$ are extended by
continuity to the boundary $\partial K^\varphi_w$, i.e.
$F\equiv\{F_\psi:\psi\in[0,\varphi]\}$.

 (ii) For all $F\in \mathsf H^\ell_\zeta(K^\varphi_w;X)$ and $\psi\in[0,\varphi]$ the
 estimate \eqref{estim} is valid.

(iii) The class $\mathsf H^\ell_\zeta(K^\varphi_w;X)$ endowed with
the norm \eqref{NormH} is a Banach space.

(iv) The transformation \eqref{transf} implements an isometric
isomorphism.

(v) Let $\{F_\psi:\psi\in [0,\varphi]\}\in \mathsf
H^\ell_{\zeta}(K_w^\varphi)$, $\mathcal F\in \mathscr
H^\ell_w(\mathcal K^\varphi_\zeta)$, and let the equality $$
F\!\!\upharpoonright_{e^{-i\psi}\Bbb R+w}=\Bbb
T^\psi_{\zeta,w}(\mathcal F\!\!\upharpoonright_{e^{i\psi}\Bbb
R+\zeta})$$
 be valid for at least one value of
$\psi\in[0,\varphi]$. Then this equality is valid for all $\psi\in
[0,\varphi]$.

(vi) A set of distributions $\{F_\psi:\psi\in
[0,\varphi]\}\in\mathsf H^\ell_\zeta(K^\varphi_w;X)$ can be
uniquely recovered from any distribution $F_\psi$ of this set.

(vii) For any $\psi\in[0,\varphi]$ the set of distributions
$\{F_\psi: F\in\mathsf H^\ell_\zeta(K^\varphi_w;X)\}$ is dense in
the Sobolev space $\mathsf W^\ell_\zeta(e^{-i\psi}\Bbb R+w;X)$.
\end{thm}
\begin{pf} It suffices to consider the case $\ell<0$,
the case $\ell\geq 0$ is covered by Theorem~\ref{p3} and
Lemma~\ref{lno3}. Let us show that for every $F\in \mathsf
H^\ell_\zeta(K^\varphi_w;X)$
 there exists $\mathcal
F\in\mathscr H^\ell_w(\mathcal K^\varphi_\zeta;X)$
 such that $\mathcal F\!\!\upharpoonright_{e^{i\psi}\Bbb R+\zeta}=(\Bbb
 T^\psi_{\zeta,w})^{-1}F_\psi$ for all $\psi\in(0,\varphi)$; here the operator $(\Bbb
 T^\psi_{\zeta,w})^{-1}$ is well-defined due to
 Corollary~\ref{noCol}. Let $\{G^j\in\Bbb
 H^0_\zeta(K^\varphi_w;X)\}_{j=0}^m$ be the set of functions from
 Definition~\ref{d2}. By Theorem~\ref{p3} there exists a unique set $\{\mathcal
 G^j\in\mathcal
 H^0_w(\mathcal K^\varphi_\zeta;X)\}_{j=0}^m$ such that
$\mathcal G^j\!\!\upharpoonright_{e^{i\psi}\Bbb R+\zeta}=(\Bbb
T^\psi_{\zeta,w})^{-1}G^j\!\!\upharpoonright_{e^{-i\psi}\Bbb R+w}$
for all $\psi\in(0,\varphi)$ and $j=0,\dots,m$. Denote $\mathcal
F(\lambda)=\sum_{j=0}^m \lambda^j \mathcal G^j(\lambda)$,
$\lambda\in\mathcal K^\varphi_\zeta$. It is clear that the
function $\mathcal K^\varphi_\zeta\ni\lambda\mapsto \mathcal
F(\lambda)\in X$ is analytic. The rule \eqref{rule} and the
representation in \eqref{er} for $F_\psi$   give
$$
\mathcal F\!\!\upharpoonright_{e^{i\psi}\Bbb R+\zeta}=(\Bbb
T^\psi_{\zeta,w})^{-1} \sum_{j=0}^m
D^j_{(\psi)}G^j\!\!\upharpoonright_{e^{-i\psi}\Bbb R+w}=(\Bbb
T^\psi_{\zeta,w})^{-1}F_\psi, \quad \psi\in(0,\varphi).
$$
The estimate \eqref{er} together with Corollary~\ref{noCol} leads
to the uniform in $\psi\in(0,\varphi)$ estimate
$$
\|\mathcal F; \mathscr W^\ell_w(e^{i\psi}\Bbb R+\zeta;X)\|\leq C.
$$
Thus $\mathcal F\in \mathscr H^\ell_w(\mathcal
K^\varphi_\zeta;X)$. We set $F_0=\Bbb T^\psi_{\zeta,w}(\mathcal
F\!\!\upharpoonright_{\Bbb R+\zeta})$ and $ F_\varphi=\Bbb
T^\psi_{\zeta,w}(\mathcal F\!\!\upharpoonright_{e^{i\varphi}\Bbb
R+\zeta})$. Taking into account Corollary~\ref{noCol} and the
items 1.(i), 1.(ii) of Proposition~\ref{p2}, we arrive at the
relations \eqref{fas} and the estimate \eqref{estim}. Moreover, if
we define the norm in $\mathsf H^\ell_\zeta(K^\varphi_w;X)$ as in
\eqref{NormH} then we have
$$
\|F;\mathsf H^\ell_\zeta(K^\varphi_w;X)\|=\|\mathcal F;\mathscr
H^\ell_w(\mathcal K^\varphi_\zeta;X)\|;
$$
see \eqref{5}. The assertions (i), (ii) are proved.

 To complete
the proof of the assertion (iv) it remains to show that the set
$\{F_\psi:F_\psi=\Bbb T^\psi_{\zeta, w}(\mathcal
F\!\!\upharpoonright_{e^{i\psi}\Bbb
R+\zeta}),\psi\in(0,\varphi)\}$ is in the class $\mathsf
H^\ell_\zeta(K^\varphi_w;X)$ for every $\mathcal F\in\mathscr
H^\ell_w(\mathcal K^\varphi_\zeta;X)$. Let $\mathcal F\in\mathscr
H^\ell_w(\mathcal K^\varphi_\zeta;X)$ and let $m\in\Bbb N$,
$-m<\ell<0$. For all $\lambda\in\mathcal K^\varphi_\zeta$ we set
$$ \mathcal G^-(\lambda)= (
\lambda-\zeta-i)^{-m}\mathcal F^-(\lambda), \quad \mathcal
G^+(\lambda)= ( \lambda-\zeta+i)^{-m}\mathcal F^+(\lambda);$$ here
the notations are the same as in Proposition~\ref{p1.5}. By this
proposition we get $\mathcal G^\pm\in \mathscr H^0_w(\mathcal
K^\varphi_\zeta;X)$. Due to Theorem~\ref{p3} there exist unique
functions $G^+\in \mathsf H^0_\zeta(K^\varphi_w;X)$ and $G^-\in
\mathsf H^0_\zeta(K^\varphi_w;X)$ such that
$G^\pm\!\!\upharpoonright_{e^{-i\psi}\Bbb R+w}=\Bbb
T^\psi_{\zeta,w}\mathcal G^\pm\!\!\upharpoonright_{e^{i\psi}\Bbb
R+\zeta}$, $\psi\in[0,\varphi]$. Now by the rule \eqref{rule} we
have
\begin{equation}\label{rep3}
F_\psi=(D_{(\psi)}-\zeta+i)^m
(G^+\!\!\upharpoonright_{e^{-i\psi}\Bbb
R+w})+(D_{(\psi)}-\zeta-i)^m
(G^-\!\!\upharpoonright_{e^{-i\psi}\Bbb R+w}).
\end{equation}
Obviously, the representation~\eqref{rep3} can be rewritten in the
same form as in \eqref{er}. The estimate in~\eqref{er} follows
from the item 1.(ii) of Proposition~\ref{p2}, the definition
\eqref{5} of the norm in $\mathscr H^\ell_w(\mathcal
K^\varphi_\zeta;X)$,  and Corollary~\ref{noCol}. The assertion
(iv) is proved. The assertion (iii) is readily apparent from (iv)
and Proposition~\ref{p3`}.

Let the assumptions of the assertion (v) be fulfilled. Then by
(iv) there exists $\mathcal G\in \mathscr H^\ell_w(\mathcal
K^\varphi_\zeta;X)$ such that  $(\Bbb
T^\psi_{\zeta,w})^{-1}F_\psi=\mathcal
G\!\!\upharpoonright_{e^{i\psi}\Bbb R+\zeta}$ for all $\psi\in
[0,\varphi]$. For some $\psi\in[0,\varphi]$ we have $\mathcal
G(\lambda)=\mathcal F(\lambda)$, $\lambda\in{e^{i\psi}\Bbb
R+\zeta}$ . If $\psi\in (0,\varphi)$ then the equality $\mathcal
G(\lambda)=\mathcal F(\lambda)$ can be extended by analyticity to
all $\lambda\in\mathcal K^\varphi_\zeta$; if $\psi=0$ or
$\psi=\varphi$ then the equality $\mathcal G(\lambda)=\mathcal
F(\lambda)$, $\lambda\in\mathcal K^\varphi_\zeta$, is a
consequence of Proposition~\ref{lost}. Hence the elements
$\mathcal G$ and $\mathcal F$ of the space $\mathscr
H^\ell_w(\mathcal K^\varphi_\zeta;X)$ are coincident. The
assertion  (v) is proved. The proof of the assertion (vi) is
similar. Finally, the assertion (vii) is a consequence of (iv)
together with Corollary~\ref{noCol} and
Proposition~\ref{p2}.1,(v). \qed\end{pf}



\begin{thm}  Let $w,\zeta\in\Bbb C$, $\ell\in\Bbb
R$, and $\varphi\in(0,\pi]$. Suppose that  $\ell> r$ and the
boundary limits of a set of distributions $F\in\mathsf
H^r_\zeta(K^\varphi_w;X)$ satisfy the inclusions $F_0\in \mathsf
W^\ell_\zeta(\Bbb R+w;X)$ and $F_\varphi\in \mathsf
W^\ell_\zeta(e^{-i\varphi}\Bbb R+w;X)$. Then $F\in\mathsf
H^k_\zeta(K^\varphi_w;X)$ for any $k<\ell$.
\end{thm}

\begin{pf}
 From Theorem~\ref{p1.13}, (iv) and Corollary~\ref{noCol}
we see that $F$ is the
 Fourier-Laplace transform of
 a function $\mathcal F\in \mathscr
H^r_w(\mathcal K^\varphi_\zeta;X)$ with the  boundary limits
$\mathcal F\!\!\upharpoonright_{\Bbb R+\zeta}\in\mathscr
W^\ell_w(\Bbb R+\zeta;X)$ and $\mathcal
F\!\!\upharpoonright_{e^{i\varphi}\Bbb R+\zeta}\in \mathscr
W^\ell_w(e^{i\varphi}\Bbb R+\zeta;X)$. By Proposition~\ref{p2.6},
(i) the function $\mathcal F$ can be represented in the form
\eqref{cau},~\eqref{1cau}, where we take $s\in (r-1/2,r)$
 such that $\ell-s>1/2$.   The boundary limits
$\mathcal F\!\!\upharpoonright_{\Bbb R+\zeta}$ and $\mathcal
F\!\!\upharpoonright_{e^{i\varphi}\Bbb R+\zeta}$ play the role of
the function $\mathcal G$ on $\partial \mathcal K^\varphi_\zeta$
in Proposition~\ref{col}. This proposition gives $\mathcal F\in
\mathscr H^{k}_w(\mathcal K^\varphi_\zeta;X)$, where $k\in
(r,s+1/2)$. Hence $F\in\mathsf H^k_\zeta(K^\varphi_w;X)$ with the
same $k$. On the next step we can apply the same argument to
$F\in\mathsf H^r_\zeta(K^\varphi_w;X)$, where $r=k$. Step by step
the parameter $k$ comes arbitrarily close to $\ell$. \qed\end{pf}

\subsection{Representation of distributions in terms of Hardy-Sobolev spaces}\label{sec}

 By analogy with $\mathsf
H^\ell_\zeta(K^\varphi_w;X)$ we can
 introduce the weighted Hardy-Sobolev space $\mathsf
H^\ell_\zeta(K^{-\varphi}_w;X)$ in the cone
$K^{-\varphi}_w=\{z\in\Bbb C: z= e^{i\varphi}v+w, v\in
K^\varphi_0\}$. We identify the spaces $\mathsf
H^\ell_\zeta(K^{-\varphi}_w;X)$ and $\mathsf
H^\ell_\eta(K^\varphi_w;X)$ by the rule: a set of distributions
$\{F_\psi\in \mathsf W^\ell_\zeta(e^{i\psi}\Bbb R+w;X):\psi\in
[0,\varphi]\}$ is in the space $\mathsf
H^\ell_\zeta(K^{-\varphi}_w;X)$ if and only if $\{F_\psi\circ
\sigma_{\varphi,w}:\psi\in[0,\varphi]\}\in\mathsf
H^\ell_\eta(K^\varphi_w;X)$, where $\eta=e^{i\varphi}\zeta$ and
$\sigma_{\varphi,w}: K^\varphi_w\to K^{-\varphi}_w$ is the linear
transformation $\sigma_{\varphi,w}(z)=e^{i\varphi}(z-w)+w$. We
equip the space $\mathsf H^\ell_\zeta(K^{-\varphi}_w;X)$ with the
norm
$$
\|F;\mathsf H^\ell_\zeta(K^{-\varphi}_w;X)\|=\|F_0; \mathsf
W^\ell_\zeta(\Bbb R+w;X)\|+\|F_\varphi; \mathsf
W^\ell_\zeta(e^{i\varphi}\Bbb R+w;X)\|.
$$

The following proposition is a  direct consequence of
Proposition~\ref{H^s}.
\begin{prop}
Let $\zeta,w\in\Bbb C$, $\ell\in\Bbb R$, and let
$0<\varphi<\pi/2$. Then every distribution $F\in\mathsf
W^\ell_\zeta(\Bbb R+w;X)$ can be represented as the sum
$F_0^++F_0^-$ of the boundary limits $F^\pm_0\in \mathsf
W^\ell_\zeta(\Bbb R+w;X)$ of some sets of distributions
$F^+\in\mathsf H^\ell_\zeta(K^\varphi_w;X)$ and $F^-\in\mathsf
H^\ell_\zeta(K^{-\varphi}_w;X)$.
\end{prop}

We introduce the space  $\mathsf W^\ell_\zeta(e^{-i\phi}\Bbb
R^\pm+w;X)$ as the set of all distributions  $F\in\mathsf
W^\ell_\zeta(e^{-i\phi}\Bbb R+w;X)$ supported on the set
$(e^{-i\phi}{\Bbb R^\pm}+w)\cup\{w\}$; here $\zeta,w\in\Bbb C$,
$\ell\in\Bbb R$, and  $\phi$ is an angle. We equip the space
$\mathsf W^\ell_\zeta(e^{-i\phi}\Bbb R^\pm+w;X)$ with the norm
\begin{equation}\label{nw+}
\|F;\mathsf W^\ell_\zeta(e^{-i\phi}\Bbb R^\pm+w;X)\|=\|F;\mathsf
W^\ell_\zeta(e^{-i\phi}\Bbb R+w;X)\|;
\end{equation}
the norm in the Sobolev space  $\mathsf
W^\ell_\zeta(e^{-i\phi}\Bbb R+w;X)$ is defined in \eqref{NOS}. It
is clear that $\mathsf W^\ell_\zeta(e^{-i\phi}\Bbb
R^++w;X)=\mathsf W^\ell_\zeta(e^{-i(\phi+\pi)}\Bbb R^-+w;X)$. As a
direct consequence of the Sobolev embedding theorem we have
$D^j_{(\phi)}F(w)=0$, $j=0,\dots,m-1$, for all $F\in \mathsf
W^\ell_\zeta(e^{-i\phi}\Bbb R^\pm+w;X)$, where $\ell>m-1/2$,
$m\in\Bbb N$. For the proof of the next proposition we refer to
\cite[Proposition~7.1]{R12}.

\begin{prop}\label{pr1}
Let  $\ell\geq 0$ and let $\mathcal W^\ell_\zeta(e^{-i\phi}\Bbb
R^++w;X)$ denote the Sobolev space of functions $e^{-i\phi}\Bbb
R^++w\ni z\mapsto f(z)\in X$; we define the norm in $\mathcal
W_\zeta^\ell(e^{-i\phi}\Bbb R^++w;X)$ by setting the value
$\|F;\mathcal W_\zeta^\ell(e^{-i\phi}\Bbb R^++w;X)\|^2$ to be
equal to the right hand side of the equality \eqref{snorm}, where
$\Bbb R$ is replaced by $\Bbb R^+$. Suppose that
$\ell\neq[\ell]+1/2$, where $[\ell]$ is the integer part of
$\ell$. Then a function $F\in \mathcal W_\zeta^\ell(e^{-i\phi}\Bbb
R^++w;X)$ extended  to the half-line $e^{-i\phi}\Bbb R^-+w$ by
zero falls into the space $\mathsf W^\ell_\zeta(e^{-i\phi}\Bbb
R^++w;X)$ if and only if
$$
D^j_{(\phi)}F(w)=0,\quad 0\leq j\leq[\ell]-1,
$$
\begin{equation}\label{cite}
\int_{e^{-i\phi}\Bbb R^++w}|z-w|^{-2(\ell-[\ell])}\|e^{-i\zeta
z}D^{[\ell]}_{(\phi)}F(z)\|^2\,|dz|<+\infty, \quad \ell>[\ell].
\end{equation}
(If $\ell>[\ell]+1/2$ then the condition~\eqref{cite} gives
$D^{[\ell]}_{(\phi)}F(w)=0$.)

 For all $F\in \mathsf W^\ell_\zeta(e^{-i\phi}\Bbb
 R^++w;X)$ the estimates
 $$\begin{aligned}
 \|F; \mathcal W_\zeta^\ell(e^{-i\phi}\Bbb
 R^++w;X)\|\leq c_1\|F;\mathsf W_\zeta^\ell(e^{-i\phi}\Bbb
 R+w;X)\|\\\leq c_2\|F; \mathcal W_\zeta^\ell(e^{-i\phi}\Bbb
 R^++w;X)\|
 \end{aligned}
 $$
 are valid, where the constants $c_1$ and $c_2$ are independent of
 $F$.
\end{prop}

In the sequel we shall need the following variant of the
Paley-Wiener theorem.
\begin{thm}[Paley-Wiener]\label{P-W thm} Let $\zeta,w\in\Bbb C$, $\ell\in\Bbb R$, and let $\phi$ be an angle.
If we identify the functions $\mathcal F$ from the Hardy class
$\mathscr H^\ell_w(e^{i\phi}\Bbb C^++\zeta;X)$ with their
 boundary limits $\mathcal
F\!\!\upharpoonright_{e^{i\phi}\Bbb R+\zeta}$ then the
Fourier-Laplace transformation \eqref{no01} yields an isometric
isomorphism
\begin{equation}\label{PW}
\Bbb T^\phi_{\zeta,w}:\mathscr H^\ell_w(e^{i\phi}\Bbb
C^++\zeta;X)\to \mathsf W^\ell_\zeta(e^{-i\phi}\Bbb R^-+w;X).
\end{equation}
\end{thm}
\begin{pf} An alternative proof of this theorem for the case $\ell\geq 0$ can be found in~\cite{R12}.
 Proposition~\ref{pr},(v)
 allows us to
 identify the functions
$\mathcal F\in\mathscr H^\ell_w(e^{i\phi}\Bbb C^++\zeta;X)$ with
their  boundary limits $\mathcal F_\phi\in \mathscr
W^\ell_w(e^{i\phi}\Bbb R+\zeta;X)$. The mapping \eqref{PW}  is
isometric due to Corollary~\ref{noCol} and the
definitions~\eqref{NH}, \eqref{nw+} of the norms in $\mathscr
H^\ell_w(e^{i\phi}\Bbb C^++\zeta;X)$ and in $\mathsf
W^\ell_\zeta(e^{i\phi}\Bbb R^-+\zeta;X)$. Let us prove
that~\eqref{PW} is an isomorphism. Without loss of generality we
can suppose that $\phi=0$.

 {\it Epimorphism.} Let us show that $\Bbb T^0_{\zeta,
w}\mathcal F\in \mathsf W^\ell_\zeta(\Bbb R^-+w;X)$ if $\mathcal
F\in \mathscr H^\ell_w(\Bbb C^++\zeta;X)$. Due to
Corollary~\ref{noCol} the inclusion $\Bbb T^0_{\zeta,w}\mathcal
F\in \mathsf W^\ell_\zeta(\Bbb R+w;X)$ holds. It remains to show
that the distribution $\Bbb T^0_{\zeta,w}\mathcal F$  is supported
on the set $ (\Bbb R^-+w)\cup\{w\}$. We extend the function
$\mathcal F$ to the half-plane $\Bbb C^-+\zeta$ by zero. Then the
inclusion $\mathcal F\in\mathscr H^\ell_w(\mathcal K^\pi_\zeta;X)$
is valid and $\mathcal F$ has
 boundary limits $\mathcal F_0, \mathcal
F_\pi\in\mathscr W^\ell_w(\Bbb R+\zeta; X)$ on $\partial \mathcal
K^\pi_\zeta$; see Proposition~\ref{p2}. It is clear that $\mathcal
F_0(\lambda)=0$ for $\lambda\in \Bbb R^-+\zeta$ and $\mathcal
F_\pi(\lambda)=0$ for $\lambda\in\Bbb R^++\zeta$. Moreover,
$\mathcal F=\mathcal F_0+\mathcal F_\pi$ on $\Bbb R+\zeta$.
Therefore, $\Bbb T^0_{\zeta,w}=-\Bbb T^\pi_{\zeta,w}$; this is
because the transformation $\Bbb T^\psi_{\zeta,w}$ is an extension
of the transformation \eqref{noTi}, where the integration runs in
opposite directions for $\psi=0$ and $\psi=\pi$, see
Lemmas~\ref{lno1},~\ref{lno2}. We have $\Bbb T^0_{\zeta,w}\mathcal
F=\Bbb T^0_{\zeta,w}\mathcal F_0-\Bbb T^\pi_{\zeta,w}\mathcal
F_\pi$. Consider for example the transformation $\Bbb
T^{\pi/2}_{\zeta,w}\mathcal F\!\!\upharpoonright_{i\Bbb R+\zeta}$.
It defines an analytic in $i\Bbb C^-+w$ function $F$ because
$\mathcal F\!\!\upharpoonright_{i\Bbb R^-+\zeta}=0$. Clearly, the
function $F$ coincides with $\Bbb T^{\psi}_{\zeta,w}\mathcal
F\!\!\upharpoonright_{e^{i\psi}\Bbb R+\zeta}$ on the set $(i\Bbb
C^-+w)\cap (e^{-i\psi}\Bbb R+w)$ for all $\psi\in[0,\pi]$. Thus
$\Bbb T^0_{\zeta,w}\mathcal F_0=\Bbb T^\pi_{\zeta,w}\mathcal
F_\pi=F$ on $\Bbb R^++w$, the distribution $\Bbb
T^0_{\zeta,w}\mathcal F_0-\Bbb T^\pi_{\zeta,w}\mathcal F_\pi$ is
supported on the set $ (\Bbb R^-+w)\cup\{w\}$, and $\Bbb
T^0_{\zeta, w}\mathcal F\in \mathsf W^\ell_\zeta(\Bbb R^-+w;X)$.

{\it Monomorphism.} Here we prove that for any $F\in \mathsf
W^\ell_\zeta(\Bbb R^-+w;X)$
  there exists $\mathcal F\in\mathscr H^\ell_w(\Bbb
 C^++\zeta;X)$ such that $(\Bbb T^0_{\zeta,w})^{-1}F=\mathcal F\!\!\upharpoonright_{\Bbb
 R+\zeta}$.
  We first consider the classical case $\ell=0$. Since the function $F$ is supported on the set $\Bbb R^{-}+w$, the
inverse Fourier-Laplace transform $(\Bbb T^0_{\zeta,w})^{-1}F$
defines an  analytic function $\Bbb C^++\zeta\ni\lambda\mapsto
\mathcal F(\lambda)\in X$ such that
\begin{equation}\label{eta}
\mathcal F\!\!\upharpoonright_{\Bbb R+\zeta+\eta}=(\Bbb
T^0_{\zeta,w})^{-1}e_\eta F,\quad\forall \eta\in \Bbb C^+,
\end{equation}
where $e_\eta:z\mapsto \exp(-i\eta z)$ is the weight function.
(The equality \eqref{eta} can be easily seen from the integral
representation~\eqref{noT} for $(\Bbb T^0_{\zeta,w})^{-1}F$.)
  The formula
\eqref{eta} together with Corollary~\ref{noCol}  leads to the
equalities
\begin{equation}\label{eq4} \|\mathcal
F(\cdot+\eta);\mathscr W^0_w(\Bbb R+\zeta;X)\|=\|e_\eta F;\mathsf
W^0_\zeta(\Bbb R^-+w;X)\|,
\end{equation}
\begin{equation}\label{eq3}
\|\mathcal F(\cdot+\eta)-(\Bbb T^0_{\zeta,w})^{-1}F;\mathscr
W^0_w(\Bbb R+\zeta;X)\|=\|e_\eta F-F;\mathsf W^0_\zeta(\Bbb
R^-+w;X)\|.
\end{equation}
It is easily seen that the uniform in $\eta\in \overline{\Bbb
C^+}$ estimate
\begin{equation}\label{eq1}
\|e_\eta F;\mathsf W^0_\zeta(\Bbb R^-+w;X)\|\leq C \|F;\mathsf
W^0_\zeta(\Bbb R^-+w;X)\|
\end{equation}
is valid, where $C$ is independent of $F\in \mathsf W^0_\zeta(\Bbb
R^-+w;X)$; moreover,
\begin{equation}\label{eq2}
\|e_\eta F-F;\mathsf W^0_\zeta(\Bbb R^-+w;X)\|\to 0,\quad
\eta\to\zeta, \eta\in {\Bbb C^+}.
\end{equation}
The estimate~\eqref{eq1} together with~\eqref{eq4} proves the
inclusion $\mathcal F\in\mathscr H^0_w(\Bbb C^++\zeta;X)$; see
Proposition~\ref{HH}, (i). Due to \eqref{eq3} and \eqref{eq2} the
inverse Fourier-Laplace transform $(\Bbb T^0_{\zeta,w})^{-1}F$ and
the boundary limits $\mathcal F\!\!\upharpoonright_{\Bbb R+\zeta}$
are coincident in $\mathscr W^0_w(\Bbb R+\zeta;X)$.

Let us consider the case of a nonnegative integer $\ell$. In this
case the norm~\eqref{nooo} is equivalent to the norm~\eqref{nw}
and $D^j_{(0;w)}F\in \mathsf W^0_\zeta(\Bbb R^-+w;X)$,
$j\leq\ell$. Thus far everything said in the previous case is
applicable here. Taking into account the rule~\eqref{rule}  we see
that the functions $\lambda\mapsto \lambda^j\mathcal F(\lambda)$,
$j\leq \ell$, are in the space $\mathscr H^0_w(\Bbb C^++\zeta;X)$.
Consequently $\mathcal F\in \mathscr H^\ell_w(\Bbb C^++\zeta;X)$
and $(\Bbb T^0_{\zeta,w})^{-1}F=\mathcal
F\!\!\upharpoonright_{\Bbb
 R+\zeta}$.

Let $\ell$ be a negative integer. If $F\in \mathsf
W^\ell_\zeta(\Bbb R^-+w;X)$ then $F$ is a distribution supported
on $(\Bbb R^-+w)\cup\{w\}$ and representable in the form
$F=\sum_{j\leq|\ell|}D_{(0;w)}^j G_j$, where $G_j\in\mathsf
W^0_\zeta(\Bbb R^-+w;X)$. Then $(\Bbb
T^0_{\zeta,w})^{-1}G_j=\mathcal G_j\!\!\upharpoonright_{\Bbb
R+\zeta}$
 and $\mathcal G_j\in \mathscr H^0_w(\Bbb C^++\zeta;X)$,
 $j\leq|\ell|$. The function $\mathcal F(\lambda)=\sum_{j\leq|\ell|}\lambda^j\mathcal
G_j(\lambda)$ is in the space $\mathscr H^\ell_w(\Bbb
C^++\zeta;X)$.  From the rule~\eqref{rule} it follows that
 $\mathcal F\!\!\upharpoonright_{\Bbb R+\zeta}=(\Bbb
T^0_{\zeta,w})^{-1}F$.

Now we are in position to consider the general case $\ell\in\Bbb
R$. The embedding  $\mathsf W^\ell_\zeta(\Bbb
R^-+w;X)\subseteq\mathsf W^{[\ell]}_\zeta(\Bbb R^-+w;X)$ is
fulfilled. Hence  we have $(\Bbb T^0_{\zeta,w})^{-1}F=\mathcal
F\!\!\upharpoonright_{\Bbb R+\zeta}$, where $\mathcal F\in
\mathscr H^{[\ell]}_w(\Bbb C^++\zeta;X)$ and $\mathcal
F\in\mathscr W^\ell_w(\Bbb R+\zeta;X)$. Define the function $
\mathcal G(\lambda)=(\lambda-\zeta+i)^{\ell-[\ell]-1}\mathcal
F(\lambda)$. This function satisfies the inclusions $\mathcal
G\in\mathscr H^{2[\ell]-\ell}_w(\Bbb C^++\zeta;X)$ and $\mathcal
G\in\mathscr W^{[\ell]+1}_w(\Bbb R+\zeta;X)$; see the
estimates~\eqref{triv}. The last inclusion and
Corollary~\ref{noCol} give $\Bbb T^0_{\zeta,w}(\mathcal
G\!\!\upharpoonright_{\Bbb R+\zeta})\in \mathsf
W^{[\ell]+1}_\zeta(\Bbb R+w;X)$.  As we already know, the
Fourier-Laplace transform $\Bbb T^0_{\zeta,w}(\mathcal
J\!\!\upharpoonright_{\Bbb R+\zeta})$ of any function $\mathcal
J\in\mathscr H^{2[\ell]-\ell}_w(\Bbb C^++\zeta;X)$ is supported on
$(\Bbb R^-+w)\cup\{w\}$. Therefore the transform $G=\Bbb
T^0_{\zeta,w}(\mathcal G\!\!\upharpoonright_{\Bbb R+\zeta})$ is in
the space $\mathsf W^{[\ell]+1}_\zeta(\Bbb R^-+w;X)$. From the
proved cases we conclude that the transformation $(\Bbb
T^0_{\zeta,w})^{-1}G$ defines a function $\widetilde{\mathcal
G}\in\mathscr H^{[\ell]+1}_w(\Bbb C^++\zeta;X)$ such that
$\widetilde{\mathcal G}\!\!\upharpoonright_{\Bbb
R+\zeta}={\mathcal G}\!\!\upharpoonright_{\Bbb R+\zeta}$. By
Proposition~\ref{pr}, (v) we have ${\mathcal
G}\equiv\widetilde{\mathcal G} \in\mathscr H^{[\ell]+1}_w(\Bbb
C^++\zeta;X)$. Then the function $\mathcal F$ is in the space
$\mathscr H^{\ell}_w(\Bbb C^++\zeta;X)$. \qed\end{pf}

\begin{cor}\label{scolor}  Let $\zeta,w\in\Bbb C$, $\ell\in\Bbb R$, and let $\phi$ be an angle.

 (i) The space $\mathsf W^\ell_\zeta(e^{-i\phi}\Bbb R^\pm+w;X)$ with
the norm \eqref{nw+} is a Banach space.

(ii) For all $\eta\in e^{i\phi}{\Bbb C^+}+\zeta$ we have $\mathsf
W^\ell_\zeta(e^{-i\phi}\Bbb R^-+w;X)\subset\mathsf
W^\ell_\eta(e^{-i\phi}\Bbb R^-+w;X)$ and
\begin{equation*}
\|F;\mathsf W^\ell_\eta(e^{-i\phi}\Bbb R^-+w;X)\|\leq C
e^{\Im\{(\eta-\zeta)w\}}\|F;\mathsf W^\ell_\zeta(e^{-i\phi}\Bbb
R^-+w;X)\|,
\end{equation*}
where the constant $C$ is independent of $w$, $\eta$, and $F\in
\mathsf W^\ell_\zeta(e^{-i\phi}\Bbb R^-+w;X)$.

(iii) For all $v\in e^{-i\phi}\Bbb R^-+w$ we have $\mathsf
W^\ell_\zeta(e^{-i\phi}\Bbb R^-+v;X)\subset\mathsf
W^\ell_\zeta(e^{-i\phi}\Bbb R^-+w;X)$ and
\begin{equation*}
\|F;\mathsf W^\ell_\zeta(e^{-i\phi}\Bbb R^-+w;X)\|=\|F;\mathsf
W^\ell_\zeta(e^{-i\phi}\Bbb R^-+v;X)\|,\quad F\in\mathsf
W^\ell_\zeta(e^{-i\phi}\Bbb R^-+v;X).
\end{equation*}
\end{cor}
\begin{pf}
(i) The space $\mathsf W^\ell_\zeta(e^{-i\phi}\Bbb R^-+w;X)$ is
complete because it is isometrically isomorphic to the Banach
space $\mathscr H^\ell_w(e^{i\phi}\Bbb C^++\zeta;X)$. (ii) For
$F\in\mathsf W^\ell_\zeta(e^{-i\phi}\Bbb R^-+w;X)$ we have the
equality
\begin{equation}\label{neta}
 \mathcal F\!\!\upharpoonright_{e^{i\phi}\Bbb
R+\eta}=e^{i(\eta-\zeta)w}(\Bbb T^\phi_{\eta,w})^{-1}F,
\end{equation}
where $\mathcal F=(\Bbb T^\phi_{\zeta,w})^{-1}F$, $\mathcal
F\in\mathscr H^\ell_w(e^{i\phi}\Bbb C^++\zeta)$. For $\ell\geq 0$
the equality~\eqref{neta} can be easily seen from the integral
representation~\eqref{noT}. The rule~\eqref{rule} allows us to
extend~\eqref{neta} to the case $\ell<0$. Corollary~\ref{emb},(ii)
together with Corollary~\ref{noCol} and the equality~\eqref{neta}
establishes the assertion. (iii) The embedding is obvious. The
equality is valid because the norm of the space $\mathsf
W^\ell_\zeta(e^{i\phi}\Bbb R+v;X)$ does not change while $v$
travels along the line  $e^{i\phi}\Bbb R+w$. \qed\end{pf} Let us
note that the assertions (iii), (iv) of Proposition~\ref{pr}
furnish estimates for the transform $\mathcal F=(\Bbb
T^\phi_{\zeta,w})^{-1}F$ of a function $F\in\mathsf
W^\ell_\zeta(e^{i\phi}\Bbb R^-+w;X)$.

\begin{prop}\label{credit} Let $\ell>1/2$ and $F\in \mathsf W^\ell_\zeta(e^{-i\phi}\Bbb
R+w;X)$. We set $G=F$ on $e^{-i\phi}\overline{\Bbb R^{+}}+w$ and
$G=0$ on $e^{-i\phi}\Bbb R^-+w$. Then for any $k<1/2$ the
inclusion  $G\in \mathsf W^k_\zeta(e^{-i\phi}\Bbb R^++w;X)$ and
the estimate
\begin{equation}\label{ne}
\|G;\mathsf W^k_\zeta(e^{-i\phi}\Bbb R^++w;X)\|\leq C \|F; \mathsf
W^\ell_\zeta(e^{-i\phi}\Bbb R+w;X)\|
\end{equation}
are valid. Here the constant $C$ is independent of $F$.
\end{prop}
\begin{pf} For a function $F\in \mathsf W^\ell_\zeta(e^{-i\phi}\Bbb
R+w;X)$, $\ell>1/2$, the Sobolev theorem gives $F\in
C(e^{-i\phi}\Bbb R+w;X)$. Hence the function $e^{-i\phi}\Bbb
R+w\ni z\mapsto \|F\|$ is bounded in a neighbourhood of the point
$w$, the condition~\eqref{cite} is satisfied with $\ell$ replaced
by $k$, $k<1/2$. Then from Proposition~\ref{pr1} it follows  the
inclusion $G\in \mathsf W^k_\zeta(e^{-i\phi}\Bbb R^++w;X)$ and
 the estimate~\eqref{ne}.
\qed\end{pf} Let us recall that the estimate~\eqref{est1} in
Proposition~\ref{HH}, (iii) remains  unproved  in case
$s<\ell-1/2$, $k\in (s,s+1/2)$. The next corollary  presents a
result which enables us to finalize the proof of
Proposition~\ref{HH}, (iii).
\begin{cor}\label{credit1} Let $\ell>1/2$ and $\mathcal J\in \mathscr W^\ell_0(e^{i\phi}\Bbb
R+\zeta;X)$. We set $$\mathcal F(\lambda)=\int_{e^{i\phi}\Bbb
R+\zeta}(\lambda-\mu)^{-1} \mathcal J(\mu)\,d\mu,\quad \lambda\in
e^{i\phi}\Bbb C^++\zeta.$$ Then for any  $k<1/2$ the inclusion
$\mathcal F\in\mathscr H^k_0(e^{i\phi}\Bbb C^++\zeta; X)$ and the
estimate
\begin{equation}\label{debet}
\|\mathcal F;\mathscr H^k_0(e^{i\phi}\Bbb C^++\zeta; X)\|\leq
c\|\mathcal J;\mathscr W^\ell_0(e^{i\phi}\Bbb R+\zeta;X)\|
\end{equation}
are valid. The constant $c$ is independent of $\mathcal J$.
\end{cor}
\begin{pf} Let $J\in \mathsf
W^\ell_\zeta(e^{-i\phi}\Bbb R;X)$. Consider the Fourier-Laplace
transform $J=\Bbb T^{\phi}_{\zeta,0}\mathcal J$.  By
Proposition~\ref{credit} we can represent $J$ as the sum $F+G$,
where $F\in \mathsf W^k_\zeta(e^{-i\phi}\Bbb R^-;X)$ and $G\in
\mathsf W^k_\zeta(e^{-i\phi}\Bbb R^+;X)$, $k<1/2$; moreover, the
estimate
\begin{equation}\label{debet1}
\|F;\mathsf W^k_\zeta(e^{-i\phi}\Bbb R^-;X) \|\leq C\|J; \mathsf
W^\ell_\zeta(e^{-i\phi}\Bbb R;X)\|
\end{equation}
 holds. Then $\mathcal
J=\mathcal F+\mathcal G$, where $\mathcal F=(\Bbb
T^\phi_{\zeta,0})^{-1}F$ and $\mathcal G=(\Bbb
T^\phi_{\zeta,0})^{-1}G$. As a consequence of  Theorem~\ref{P-W
thm} we have  the estimate
$$
\|\mathcal F; \mathscr H^k_0(e^{i\phi}\Bbb C^++\zeta; X)\|\leq
C\|F;\mathsf W^k_\zeta(e^{-i\phi}\Bbb R^-;X)\|.
$$
 This together with~\eqref{debet1} establishes
the estimate~\eqref{debet}. We represented the function $\mathcal
J\in \mathscr W^0_0(e^{i\phi}\Bbb R+\zeta;X)$ as the sum of
boundary limits of functions $\mathcal F\in\mathscr
H^0_0(e^{i\phi}\Bbb C^++\zeta; X)$  and $\mathcal G\in \mathscr
H^0_0(e^{i\phi}\Bbb C^-+\zeta; X)$. Hence $\mathcal F$ is the
Cauchy integral of $\mathcal J$. \qed\end{pf}

\begin{rem}\label{creditproof} To finalize the proof of Proposition~\ref{HH}, (iii)
 it suffices to apply Corollary~\ref{credit} with the
function $\mathcal J$ replaced by  $\exp\{iw
\cdot\}(\cdot-\eta)^s\mathcal J(\cdot)$ and $\mathcal F$ replaced
by $\exp\{iw\cdot\}(\cdot-\eta)^s\mathcal F(\cdot)$; here
$\mathcal F$ and $\mathcal J$ are the same as in
Proposition~\ref{HH}, (iii).
\end{rem}

\begin{prop}\label{pr0}
 Let $\zeta,w\in\Bbb C$,  $\ell\in\Bbb R$, and let $\phi$ be an angle.

(i) Every distribution $F\in\mathsf W^\ell_\zeta(e^{-i\phi}\Bbb
R+w;X)$ can be represented as $F^++F^-$, where $F^\pm\in \mathsf
W^s_\zeta(e^{-i\phi}\Bbb R^\pm+w;X)$ for
\begin{equation}\label{cnd}
\begin{aligned}
&s<1/2  &\text{ if }\ \  \ell>1/2;\phantom{\ \ \text{ and }\ \
[\ell]\leq\ell\leq[\ell]+1/2}\\
&s\leq[\ell]&\text{ if }\ \  \ell\leq 1/2\ \ \text{ and }\ \
[\ell]\leq\ell\leq[\ell]+1/2;\\
&s<[\ell]+1/2 &\text{ if }\ \  \ell\leq 1/2\ \  \text{ and }\ \
[\ell]+1/2<\ell<[\ell].
\end{aligned}
\end{equation}
 If $s<-1/2$ then this representation is
not unique due to the inclusions
$$f\cdot(D^j_t\delta)\circ\tau^{-1}_{\phi,w} \in \mathsf
W^s_\zeta(e^{-i\phi}\Bbb R^++w;X)\cap\mathsf
W^s_\zeta(e^{-i\phi}\Bbb R^-+w;X),\quad 0\leq j\leq-[s+3/2]; $$
here $f$ is a coefficient in $X$, $\delta$ denotes the Dirac delta
function, and $\tau^{-1}_{\phi,w}:e^{-i\phi}\Bbb R+w\to \Bbb R$ is
the linear transformation $\tau^{-1}_{\phi,w}(z)=e^{i\phi}(z-w)$.

(ii) For every $F\in\mathsf W^\ell_\zeta(\Bbb R+w;X)$ there exists
a set of distributions $G\in\mathsf H^s_\zeta(K^\pi_w;X)$ such
that $F=G_0-G_\pi$, where $s$ is the same as in~\eqref{cnd}, $G_0$
and $G_\pi$ are boundary limits of $G$ in the sense of
Theorem~\ref{p1.13}.  In the case $s<-1/2$ the set $G$ is not
unique. If $\widetilde G\in \mathsf H^s_\zeta(K^\pi_w;X)$ and
$F=\widetilde G_0-\widetilde G_\pi$ then the sets of distributions
$\widetilde G\equiv \{\widetilde G_\psi:\psi\in[0,\pi]\}$ and
$G\equiv \{ G_\psi:\psi\in[0,\pi]\}$ satisfy the equality
\begin{equation}\label{street}
(e_\zeta\widetilde G_\psi)\circ \tau_{\psi,w}-(e_\zeta
G_\psi)\circ \tau_{\psi,w}=\sum_{j=0}^{-[s+3/2]}
e^{i(j+1)\psi}f_jD_t^j
 \mathcal P\frac 1 t
\end{equation}
 for all $\psi\in[0,\pi]$ and some coefficients $f_j\in X$; here $e_\zeta:z\mapsto \exp(-i\zeta
 z)$,
the distribution $\mathcal P \frac 1 t (v)$, $v\in \mathscr S(\Bbb
R)$, is defined as the Cauchy principal value of the integral
$\int_\Bbb R t^{-1}{v(t)} \, dt$, and
$\tau_{\psi,w}(t)=e^{-i\psi}t+w$.

(iii)  Let $F\in  \mathsf W^\ell_\zeta(\Bbb R^-+w;X)$. Then there
exists an analytic function
\begin{equation}\label{afun}
 \Bbb C\setminus \overline{\Bbb
R^-+w}\ni z\mapsto G(z)\in X
\end{equation}
such that for any $v\in \Bbb R^++w$  the inclusion
\begin{equation}\label{Sg}
\{G\!\!\upharpoonright_{e^{-i\psi}\Bbb
R+v}:\psi\in(0,\pi)\}\in\mathsf H^\ell_\zeta(K^\pi_v;X)
\end{equation}
 and the
equality $F=G_0-G_\pi$ hold, where the distributions $G_0$ and
$G_\pi$ are  boundary limits of the set \eqref{Sg} in the sense of
Theorem~\ref{p1.13}. Moreover, there exists a set of distributions
$G^-\in\mathsf H^\ell_\zeta(K^\pi_w;X)$  such that
$F=G^-_0-G^-_\pi$  and
$$
G^-_\psi=G\text{ on the set } e^{-i\psi}\Bbb
R+w\setminus\overline{\Bbb R^++w}, \ \psi\in[0,\pi].
$$
A similar assertion is valid for $F\in  \mathsf W^\ell_\zeta(\Bbb
R^++w;X)$.
\end{prop}
\begin{pf}
(i) We represent  $\mathcal F=(\Bbb T^\phi_{\zeta,w})^{-1}F$ in
the form $\mathcal F=(\mathcal F^++\mathcal
F^-)\!\!\upharpoonright_{e^{i\phi}\Bbb R+\zeta}$, where $\mathcal
F^\pm\in\mathscr H^s_w(e^{i\phi}\Bbb C^\pm+\zeta;X)$; see
Corollary~\ref{emb}, (iii). This together with Theorem~\ref{P-W
thm} proves  the assertion.

(ii) Here again we represent $\mathcal F=(\Bbb
T^0_{\zeta,w})^{-1}F$ in the form $\mathcal F=(\mathcal
F^++\mathcal F^-)\!\!\upharpoonright_{\Bbb R+\zeta}$, where
$\mathcal F^\pm\in\mathscr H^s_w(\Bbb C^\pm+\zeta;X)$. Let us
define $\mathcal G\in\mathscr H^s_w(\mathcal K^\pi_\zeta;X)$ by
setting $\mathcal G=\mathcal F^+$ on $\Bbb C^++\zeta$ and
$\mathcal G=\mathcal F^-$ on $\Bbb C^-+\zeta$. It s clear that
$(\mathcal F^++\mathcal F^-)\!\!\upharpoonright_{\Bbb
R+\zeta}=\mathcal G_0+\mathcal G_\pi$, where $\mathcal G_0$ and
$\mathcal G_\pi$ are the
 boundary limits of $\mathcal G$; see
 Propositions~\ref{p2},~\ref{pr}. Let $G\in\mathsf H^s_\zeta(K^\pi_w;X)$ be the set of distributions
  $$G=\{G_\psi=\Bbb T^\psi_{\zeta,w}(\mathcal G\!\!\upharpoonright_{e^{i\psi}\Bbb R+\zeta}):\psi\in[0,\varphi]\}.$$
  It is clear  that
   $F=\Bbb T^0_{\zeta,w}(\mathcal F^++\mathcal F^-)\!\!\upharpoonright_{\Bbb R+\zeta}
   =\Bbb T^0_{\zeta,w} \mathcal G_0-\Bbb T^\pi_{\zeta,w} \mathcal
   G_\pi=G_0-G_\pi$.
    By Corollary~\ref{emb}, (iii) the representation
    $\mathcal F=(\mathcal F^++\mathcal F^-)\!\!\upharpoonright_{\Bbb R+\zeta}$ is not unique
    if
    $\ell<-1/2$. The relations~\eqref{repr} lead to the equlity~\eqref{street}.

(iii) By Theorem~\ref{P-W thm}, for $\mathcal F=(\Bbb
T^\phi_{\zeta,w})^{-1}F$ we have the inclusion $\mathcal
F\in\mathscr H^\ell_w(\Bbb C^++\zeta;X)$. Let $\mathcal G=\mathcal
F$ in $\Bbb C^++\zeta$ and $\mathcal G=0$ in $\Bbb C^-+\zeta$.  It
is clear that $\mathcal G\in \mathscr H^\ell_v(\mathcal
K^\pi_\zeta;X)$ for any $v\in\overline{\Bbb R^+}+w$  and $\mathcal
F\!\!\upharpoonright_{\Bbb R+\zeta}=\mathcal G_0+\mathcal G_\pi$,
where $\mathcal G_0$ and $\mathcal G_\pi$ are boundary limits of
$\mathcal G$; see Proposition~\ref{p2}. We can define
$G^-\in\mathsf H^\ell_\zeta(K^\pi_w;X)$ by setting $G^-_\psi=\Bbb
T^\psi_{\zeta,w}(\mathcal G\!\!\upharpoonright_{e^{i\psi}\Bbb
R+w})$. Therefore $F=\Bbb T^0_{\zeta,w}(\mathcal G_0+\mathcal
G_\pi)=G^-_0-G^-_\pi$. In the same way as in the item {\it
Epimorphism} of the proof of Theorem~\ref{P-W thm} we see that the
distributions $G^-_0$ and $G^-_\pi$ coincide on $\Bbb R^++w$ with
an analytic function. Hence the set $G^-\in\mathsf
H^\ell_\zeta(K^\pi_w;X)$ defines an analytic
function~\eqref{afun}. Let us also consider the set of
distributions
\begin{equation}\label{SB}
\widetilde G\equiv\{e^{i\zeta(w-v)}\Bbb T^\psi_{\zeta,v}(\mathcal
G\!\!\upharpoonright_{e^{i\psi}\Bbb R+\zeta}):\psi\in[0,\pi]\}\in
\mathscr H^\ell_\zeta(K^\pi_v;X),\quad v\in \Bbb R^++w.
\end{equation}
By the same argument as above we conclude that the set~\eqref{SB}
defines an analytic in $\Bbb C\setminus \overline{\Bbb R^-+v}$
function $\widetilde G$. Note that $\Bbb
T^0_{\zeta,w}=e^{i\zeta(w-v)}\Bbb T^0_{\zeta,v}$ and $\Bbb
T^\pi_{\zeta,w}=e^{i\zeta(w-v)}\Bbb T^\pi_{\zeta,v}$. Thus
$\widetilde G_0=G^-_0$ and $\widetilde G_\pi=G^-_\pi$.
Consequently the set of distributions $\widetilde G$ defines the
same analytic function $G$ as the set $G^-$, we have the equality
$G\!\!\upharpoonright_{e^{-i\psi}\Bbb R+v}=\widetilde G_\psi$,
$\psi\in(0,\pi)$. \qed\end{pf}

The next proposition deals with compactly supported distributions,
 the first assertion is an integral form of the
Paley-Wiener-Schwartz theorem.
\begin{prop} \label{PWS} Let $\zeta,w\in\Bbb C$, $\ell\in\Bbb R$, and let $\varphi$ be an
angle. Assume that $v\in e^{i\varphi}\overline{\Bbb R^+}+w$.

(i){\bf (Paley-Wiener-Schwartz)}  Let $F\in\mathscr
S'_\zeta(e^{-i\varphi}\Bbb R+w;X)$. The inclusion
\begin{equation}\label{inc}
F\in \mathsf W^\ell_0(e^{-i\varphi}\Bbb R^++w;X)\cap \mathsf
W^\ell_0(e^{-i\varphi}\Bbb R^-+v;X)
\end{equation}
 holds if and only if the
inverse Fourier-Laplace transform $(\Bbb
T^\varphi_{\zeta,w})^{-1}F$ is an entire function $\Bbb
C\ni\lambda\mapsto \mathcal E(\lambda)\in X$ satisfying the
uniform in $\psi\in (0,\pi)$ estimates
\begin{equation}\label{es+}
\|\mathcal E;\mathscr W^\ell_w(e^{i(\varphi+\psi)}\Bbb
R^-+\zeta;X)\|\leq C;\quad  \|\mathcal E;\mathscr
W^\ell_v(e^{i(\varphi+\psi)}\Bbb R^++\zeta;X)\|\leq C.
\end{equation}

(ii) For any  $F\in  \mathsf W^\ell_0(\Bbb R^++w;X)\cap \mathsf
W^\ell_0(\Bbb R^-+v;X)$ there exists  a unique analytic function
$$
 \Bbb C\setminus\left(\overline{\Bbb R^++w}\cap \overline{\Bbb
R^-+v}\right)\ni z\mapsto G(z)\in X
$$
such that the following inclusion and equalities hold
\begin{equation}\label{Sg1}
\{G\!\!\upharpoonright_{e^{-i\psi}\Bbb
R+u}:\psi\in(0,\pi)\}\in\mathsf H^\ell_0(K^\pi_u;X),\ \  u\in
(\Bbb R^-+w)\cup(\Bbb R^++v);
\end{equation}
\begin{equation}\label{Sg2}
F= G_\pi-G_0\text{\ \  if\ \  } u\in\Bbb R^-+w;\quad
F=G_0-G_\pi\text{\ \ if\ \  } u\in\Bbb R^++v.
\end{equation}
Here the distributions $G_0$ and $G_\pi$ are the boundary limits
of the set \eqref{Sg1} in the sense of Theorem~\ref{p1.13}.

(iii)   Let $F\in  \mathsf W^\ell_0(\Bbb R^++w;X)\cap \mathsf
W^\ell_0(\Bbb R^-+v;X)$ and let $\mathcal E$ be the corresponding
entire function from the  assertion (i). We set $\mathcal
E^\pm=\mathcal E$ on $\Bbb C^\pm$ and $\mathcal E^\pm=0$ on $\Bbb
C^\mp$. Then the sets of distributions
\begin{equation}\label{sets}
\begin{aligned}
G^-\equiv\{G^-_\psi=-\Bbb T^\psi_{0,w}(\mathcal
E^-\!\!\upharpoonright_{e^{i\psi}\Bbb
R}:\psi\in[0,\pi]\}\in\mathsf
H^\ell_0(K^\pi_w;X), \\
G^+\equiv\{G^+_\psi=\Bbb T^\psi_{0,v}(\mathcal
E^+\!\!\upharpoonright_{e^{i\psi}\Bbb
R}:\psi\in[0,\pi]\}\in\mathsf H^\ell_0(K^\pi_v;X)
\end{aligned}
\end{equation}
satisfy the conditions $G_0^--G_\pi^-=-F$, $G_0^+-G_\pi^+=F$,  and
meet the relations
$$G^-_\psi=G\text{
on the set } e^{-i\psi}\Bbb R+w\setminus\left(\overline{\Bbb
R^++w}\cap \overline{\Bbb R^-+v}\right),
\psi\in[0,\pi],$$$$G^+_\psi=G\text{ on the set } e^{-i\psi}\Bbb
R+v\setminus\left(\overline{\Bbb R^++w}\cap \overline{\Bbb
R^-+v}\right), \psi\in[0,\pi],$$ where $G$ denotes the analytic
function from the assertion (ii).
\end{prop}
\begin{pf} (i) {\it Necessity.} If the
inclusion~\eqref{inc} is valid then
\begin{equation}\label{supp}
\supp F\subseteq\overline{\Bbb R^++w}\cap \overline{\Bbb R^-+v}.
\end{equation}
Therefore we have $F\in
\mathsf W^\ell_\zeta(e^{-i\varphi}\Bbb R^++w;X)\cap \mathsf
W^\ell_\zeta(e^{-i\varphi}\Bbb R^-+v;X)$ for any $\zeta\in\Bbb C$.
By
 Theorem~\ref{P-W thm} the transform $(\Bbb T^\varphi_{\zeta,v})^{-1}F$ is the boundary limit of a function
 $\mathcal F^+\in\mathscr H^\ell_v(e^{i\varphi}\Bbb C^++\zeta;X)$, and
 $(\Bbb T^\varphi_{\zeta,w})^{-1}F$  is the boundary limit of $\mathcal
F^-\in\mathscr H^\ell_w(e^{i\varphi}\Bbb C^-+\zeta;X)$. Since the
equality  $(\Bbb T^\varphi_{\zeta,w})^{-1}F=e^{i\zeta(v-w)}(\Bbb
T^\varphi_{\zeta,v})^{-1}F$ holds, and $\mathcal
F^+\!\!\upharpoonright_{\Bbb R+\eta}=e^{i(\eta-\zeta)v}(\Bbb
T^\varphi_{\eta,v})^{-1}F$ for $\eta\in e^{i\varphi}\Bbb
C^++\zeta$ (cf. \eqref{neta}), we conclude that $(\Bbb
T^\varphi_{\zeta,w})^{-1}F$ defines an analytic in $\Bbb C$
function $\mathcal E\in \mathscr H^\ell_w(e^{i\varphi}\Bbb
C^-+\zeta;X)\cap \mathscr H^\ell_v(e^{i\varphi}\Bbb C^++\zeta;X)$
such that $\mathcal E=e^{i\zeta(v-w)}\mathcal F^+$ on
$e^{i\varphi}\Bbb C^++\zeta$ and $\mathcal E=\mathcal F^-$ on
$e^{i\varphi}\Bbb C^-+\zeta$. The estimates \eqref{es+} are
fulfilled due to Definition~\ref{defHalfPlane}.

{\it Sufficiency.} If an analytic in $\Bbb C$ function $\mathcal
E$ satisfies the estimates \eqref{es+} then by
Definition~\ref{defHalfPlane} we have $\mathcal E\in\mathscr
H^\ell_w(e^{i\varphi}\Bbb C^-+\zeta;X)\cap\mathscr
H^\ell_v(e^{i\varphi}\Bbb C^++\zeta;X)$. For $F=\Bbb
T^\varphi_{\zeta,w}\mathcal E$ Theorem~\ref{P-W thm} gives $F\in
\mathsf W^\ell_\zeta(e^{-i\varphi}\Bbb R^++w;X)\cap \mathsf
W^\ell_\zeta(e^{-i\varphi}\Bbb R^-+v;X)$. Hence we
have~\eqref{supp} and the inclusion~\eqref{inc} holds.

 (ii, iii) Let  $\mathcal E\in\mathscr H^\ell_w(\Bbb
C^-;X)\cap \mathscr H^\ell_v(\Bbb C^+;X)$  be the entire function
from the assertion (i). We define the functions $\mathcal
E^+\in\mathscr H^\ell_v(\mathcal K^\pi_0;X)$ and $\mathcal E^-
\in\mathscr H^\ell_w(\mathcal K^\pi_0;X)$ by the equalities
$\mathcal E^\pm=\mathcal E$ on $\Bbb C^\pm$ and $\mathcal E^\pm=0$
on $\Bbb C^\mp$. It is clear that $\mathcal
E\!\!\upharpoonright_{\Bbb R}=\mathcal E^+_0+\mathcal
E^+_\pi=\mathcal E^-_0+\mathcal E^-_\pi$ and $\mathcal
E^\pm_0=\mathcal E_\pi^\mp$; here $\mathcal E^\pm_0$ and $\mathcal
E^\pm_\pi$ are boundary limits of $\mathcal E^\pm$ in the sense of
Proposition~\ref{p2}. Introduce the sets of
distributions~\eqref{sets}. Repeating the argument from the item
{\it Epimorphism} in the proof of Theorem~\ref{P-W thm} we
conclude that  the distributions $G_0^+$ and $G_\pi^+$ are equal
to an analytic function on $\Bbb R^++v$,  the distributions
$G_0^-$ and $G_\pi^-$ are equal to an analytic function on $\Bbb
R^-+w$. Thus the set of distributions $G^-\in\mathsf
H^\ell_0(K^\pi_w;X)$ defines an analytic in $\Bbb
C\setminus\overline{\Bbb R^++w}$ function $G^-$, and
$G^+\in\mathsf H^\ell_0(K^\pi_v;X)$ defines an analytic in $\Bbb
C\setminus\overline{\Bbb R^-+v}$ function $G^+$. Then for any
$u\in\Bbb R+w$ we have $F=\Bbb T^0_{0,u} (\mathcal
E\!\!\upharpoonright_{\Bbb R})=G^+_0-G^+_\pi=G^-_\pi-G^-_0$. It is
clear that $G_0^\pm=G_\pi^\mp$. Therefore the analytic functions
$G^+$ and $G^-$ are coincident on $(\Bbb R^-+w)\cup(\Bbb R^++v)$.
Hence we can define the needed analytic function $G$ as
$G=G^+=G^-$. Indeed, the inclusion $\mathcal E^- \in\mathscr
H^\ell_u(\mathcal K^\pi_0;X)$, where  $u\in \Bbb R^-+w$, and the
inclusion $\mathcal E^+ \in\mathscr H^\ell_u(\mathcal K^\pi_0;X)$,
where $u\in \Bbb R^++v$, allow us to see the relations
$$
\begin{aligned}
\{G\!\!\upharpoonright_{e^{-i\psi}\Bbb R+u}&=-\Bbb
T^\psi_{0,u}(\mathcal E^-\!\!\upharpoonright_{e^{i\psi}\Bbb
R}):\psi\in(0,\pi)\}\in\mathsf H^\ell_0(K^\pi_u;X),\quad u\in
{\Bbb R^-+w},
\\
\{G\!\!\upharpoonright_{e^{-i\psi}\Bbb R+u}&=\Bbb
T^\psi_{0,u}(\mathcal E^+\!\!\upharpoonright_{e^{i\psi}\Bbb
R}):\psi\in(0,\pi)\}\in\mathsf H^\ell_0(K^\pi_u;X),\quad u\in
{\Bbb R^++v}.
\end{aligned}
$$
This proves the inclusions~\eqref{Sg1}. The  equalities $F=-\Bbb
T^0_{0,u}(\mathcal E^-_0+\mathcal E^-_\pi)$ and $F=\Bbb
T^0_{0,u}(\mathcal E^+_0+\mathcal E^+_\pi)$ lead to~\eqref{Sg2}.
The function $G$ is unique due to Proposition~\ref{pr0}, (i).
\qed\end{pf}

\begin{exmp} Consider the Dirac delta function $\delta\in \mathsf W^{\ell}_0 (\Bbb
R;\Bbb C)$, $\ell<-1/2$, as an example of the distribution $F$ in
Proposition~\ref{PWS}, (ii) and (iii). Using the same notations as
in  Proposition~\ref{PWS} we have $\mathcal
E(\lambda)=(2\pi)^{-1/2}$ for all $\lambda\in\Bbb C$. Then the
analytic in $\Bbb C\setminus\{0\}$ function is $G(z)=-i/({2\pi
z})$, the boundary limits $G_0$ and $G_\pi$  of the
set~\eqref{Sg1} are given by the formulas
$$G_0(t)= -\frac 1 2 \delta(t)-\frac i {2\pi}
\mathcal P\frac 1 t\ \text {   and   }\  G_\pi(t)=\frac 1 2
\delta(t)-\frac i {2\pi} \mathcal P\frac 1 t\ \ \text{ if }\ \
u\in\Bbb R^-;$$
$$G_0(t)= \frac 1 2 \delta(t)-\frac i {2\pi}
\mathcal P\frac 1 t\ \text {   and   }\  G_\pi(t)=-\frac 1 2
\delta(t)-\frac i {2\pi} \mathcal P\frac 1 t\ \  \text{ if }\ \
u\in\Bbb R^+.$$ Here the distribution $\mathcal P \frac 1 t (v)$,
$v\in \mathscr S(\Bbb R)$, is defined as the Cauchy principal
value of the integral $\int_\Bbb R t^{-1}{v(t)} \, dt$.
 The distributions $G^\pm_\psi$ from the sets $G^\pm\in\mathsf H^{-1}_0(K^\pi_0;\Bbb
 C)$ are such that $G^+_\psi\circ \tau_{\psi,0} (t)=e^{i\psi}\left(\frac 1 2 \delta(t)-\frac i {2\pi}
\mathcal P\frac 1 t\right )$ and $G^-_\psi\circ
\tau_{\psi,0}(t)=e^{i\psi}\left(-\frac 1 2 \delta(t)-\frac i
{2\pi} \mathcal P\frac 1 t\right)$, where $\tau_{\psi,0}$ is the
linear transformation $\tau_{\psi,0}(t)=e^{-i\psi}t$; cf.
\eqref{no01}.
\end{exmp}

Let $\ell\in\Bbb R$ and  $s\leq\ell$. For $v\in\Bbb C$ and
$\eta\in e^{i\phi}\Bbb C^++\zeta$
 we introduce the operator
\begin{equation}\label{proj1}
\begin{aligned}
 \mathsf P_{\eta,v}^s(\phi,\zeta)=\Bbb T^\phi_{\zeta,v}\mathscr P^s_{\eta,v}(\phi,\zeta)(\Bbb
T^\phi_{\zeta,v})^{-1}:\mathsf W^\ell_\zeta (e^{-i\phi}\Bbb R+v;X)
\to\mathsf W^k_\zeta(e^{-i\phi}\Bbb R^++v;X)
\end{aligned}
\end{equation}
where  $k=s$ if $\ell-s\leq 1/2$ and $k\in[s,s+1/2)$ if
$\ell-s>1/2$. Recall that $\mathscr P^s_{\eta,v}(\phi,\zeta)$
denotes the projection operator~\eqref{proj11},~\eqref{proj}. We
shall omit the parameters $\phi$ and $\zeta$ in the notations of
the operators  $\mathsf P_{\eta,v}^s(\phi,\zeta)$ when it can be
done without ambiguity. From~\eqref{proj11} it is clearly seen
that the operator $\mathsf P_{\eta,v}^0$ does not depend on the
parameter $\eta$. For any $F\in \mathsf W^0_\zeta (e^{-i\phi}\Bbb
R+w;X)$ we have $\mathsf P_{\eta,v}^0 F=F$ on $e^{-i\phi}\Bbb
R^++v$ and $\mathsf P_{\eta,v}^0 F=0$ on $e^{-i\phi}\Bbb R^-+v$.
In the case of an integer nonpositive $s$, any distribution $F\in
\mathsf W^\ell_\zeta (e^{-i\phi}\Bbb R+w;X)$, $\ell\geq s$, can be
uniquely represented in the form $F=(D_{(\phi)}-\eta)^{-s}G$ with
some $G\in \mathsf W^{0}_\zeta (e^{-i\phi}\Bbb R+w;X)$, then
$\mathsf P^s_{\eta,v}F=(D_{(\phi)}-\eta)^{-s}\mathsf
P^0_{\eta,v}G$.  Theorem~\ref{P-W thm} and Corollary~\ref{noCol}
allows us to pass from~\eqref{1proj} to the estimate
\begin{equation}\label{la}
\|\mathsf P_{\eta,v}^s F;\mathsf W^k_\zeta(e^{-i\phi}\Bbb
R^++v;X)\|\leq C\|F;\mathsf W^\ell_\zeta(e^{-i\phi}\Bbb R+v;X)\|,
\end{equation}
where  $C$ is independent of $v\in \Bbb C$ and $F\in\mathsf
W^\ell_\zeta(e^{-i\phi}\Bbb R+v;X)$.  In other words, the norm of
the operator~\eqref{proj1} is bounded uniformly in  $v\in\Bbb C$.
Since the space $\mathsf W^\ell_\zeta(e^{-i\phi}\Bbb R+w;X)$ and
its norm are independent of $w\in e^{-i\phi}\Bbb R+v$ the norm of
the operator
$$
\mathsf P_{\eta,v}^s:\mathsf W^\ell_\zeta (e^{-i\phi}\Bbb R+w;X)
\to\mathsf W^k_\zeta(e^{-i\phi}\Bbb R^++v;X)
$$
is  bounded uniformly in  $v\in\Bbb C$ and $w\in e^{-i\phi}\Bbb
R+v$. As a consequence of the property~\eqref{2proj} we get
\begin{equation}\label{proj3}
\mathsf P_{\eta,w}^r\mathsf P_{\eta,v}^s=\mathsf
P_{\eta,v}^s,\quad r\leq s,\ v\in e^{-i\phi} \overline{\Bbb
R^+}+w.
\end{equation}
 If $s$ is
an integer nonpositive number then also the norm of the operator
 \begin{equation}\label{proj2}
 (\mathsf I-\mathsf P^s_{\eta,v}):\mathsf W^\ell_\zeta(e^{-i\phi}\Bbb R+v; X)\to
\mathsf W^k_\zeta(e^{-i\phi}\Bbb R^-+v; X),\ \ell\geq s,
 \end{equation}
is bounded uniformly in $v$, cf. \eqref{3proj}; here $\mathsf
I:\mathsf W^\ell_\zeta(e^{-i\phi}\Bbb R+v; X)\to\mathsf
W^k_\zeta(e^{-i\phi}\Bbb R+v; X)$ is the continuous embedding
operator.

\begin{prop}\label{effect} Let $\zeta,w\in\Bbb C$, $\ell\in\Bbb R$, and
$0<\varphi<\pi$. Assume that
$F\equiv\{F_\psi:\psi\in[0,\varphi]\}$  is a set of distributions
in $\mathsf H^\ell_\zeta(K^\varphi_w;X)$.  Let $v\in
e^{-i\phi}\overline{\Bbb R^+}+w$ for some $\phi\in[0,\varphi]$,
and let $G_\phi=\mathsf P_{\eta,v}^s F_\phi$, where $s\leq \ell$,
$\eta\notin(e^{i\phi}\overline{\Bbb C^-}+\zeta)\cup
\overline{\mathcal K^\varphi_\zeta}$, and $\mathsf P_{\eta,v}^s$
is the projection operator~\eqref{proj1}. Then the distribution
$G_\phi\in \mathsf W^k_\zeta(e^{-i\phi}\Bbb R+w;X)$ can be
uniquely extended to a set of distributions
$G\equiv\{G_\psi:\psi\in[0,\varphi]\}\in \mathsf
H^k_\zeta(K^\varphi_v;X)$; here  $k=s$ if $\ell-s\leq 1/2$ and
$k\in[s,s+1/2)$ if $\ell-s>1/2$. The set $G$ satisfies the
estimate
\begin{equation}\label{estimate}
\|G;\mathsf H^k_\zeta(K^\varphi_v;X)\|\leq C \|F;\mathsf
H^\ell_\zeta(K^\varphi_w;X)\|,
\end{equation}
where the constant $C$ is  independent of $w$, $\phi$, and $v$.
Moreover, in the case of an integer nonpositive $s$,   the
analytic in $K^\varphi_w$ function $F$ defined by the set
$\{F_\psi:\psi\in[0,\varphi]\}$ and the analytic in $K^\varphi_v$
function $G$ defined by the set $\{G_\psi:\psi\in[0,\varphi]\}$
are coincident on $K^{\varphi,+}_v$.

\end{prop}
\begin{pf} The existence of the set $G\in\mathsf H^k_\zeta(K^\varphi_v;X)$
and the estimate~\eqref{estimate} are consequences of
Proposition~\ref{restriction} and the definition~\eqref{proj1} of
the projection operator $\mathsf P^s_{\eta,v}$. The set $G$ is
unique
 by Theorem~\ref{p1.13}.(vi).

Let $s=0$. Then the analytic in $K^\varphi_w$ function $F$ defined
by the set $\{F_\psi:\psi\in[0,\varphi]\}$ and the analytic in
$K^\varphi_v$ function $G$ defined by the set
$\{G_\psi:\psi\in[0,\varphi]\}$ are coincident on
$K^{\varphi,+}_v$. Indeed, this assertion is trivial if
$\phi\in(0,\varphi)$. For the remaining cases we note that the
assertion is nearly the same as the assertion (ii) of
Corollary~\ref{Scol}, it can be established in a similar way. In
the case of a negative integer $s$ every distribution from the set
$\{F_\psi:\psi\in(0,\varphi)\}\in\mathsf H^s_\zeta(K^\varphi_w;X)$
can be represented in the form
$F_\psi=(D_{(\psi)}-\eta)^{-s}\widetilde F_\psi$, where
$\{\widetilde F_\psi:\psi\in(0,\varphi)\}\in\mathsf
H^0_\zeta(K^\varphi_w;X)$, and the parameter $\eta$ is outside of
the union $(e^{i\phi}\overline{\Bbb C^-}+\zeta)\cup
\overline{\mathcal K^\varphi_\zeta}$; see the proof of
Proposition~\ref{restriction} and the rule~\eqref{rule}. We have
$G_\phi=\mathsf P^s_{\eta,v}
F_\phi=(D_{(\psi)}-\eta)^{-s}\widetilde G_\phi$, where $\widetilde
G_\phi=\mathsf P^0_{\eta,v}\widetilde F_\phi$. There exists a set
$\{\widetilde G_\psi:\psi\in(0,\varphi)\}\in\mathsf
H^0_\zeta(K^\varphi_v;X)$ such that the corresponding analytic
function $\widetilde G$ coincides with the analytic function
$\widetilde F$; here $\widetilde F$ corresponds to  the set
$\{\widetilde F_\psi:\psi\in(0,\varphi)\}$. It remains to note
that $G(z)=(D_z-\eta)^{-s}\widetilde G(z)$ and
$F(z)=(D_z-\eta)^{-s}\widetilde F(z)$ for $z\in K^{\varphi,+}_v$,
where $D_z$ is the complex derivative. \qed\end{pf}

Let $\phi\in[0,\varphi]$ be a fixed angle.
Theorem~\ref{p1.13},(vi) allows us to identify every element
$\{F_\psi:\psi\in[0,\varphi]\}$ of the space $\mathsf
H^\ell_\zeta(K^\varphi_w; X)$ with $F_\phi\in\mathsf
W^\ell_\zeta(e^{-i\phi}\Bbb R+w; X)$. Then  $\mathsf
H^\ell_\zeta(K^\varphi_w; X)$ is dense in $\mathsf
W^\ell_\zeta(e^{-i\phi}\Bbb R+w; X)$ by Theorem~\ref{p1.13}, (v).
All operators defined on $\mathsf W^\ell_\zeta(e^{-i\phi}\Bbb R+w;
X)$ can be restricted to $\mathsf H^\ell_\zeta(K^\varphi_w; X)$.
In particular, the projection operator $\mathsf
P^s_{\eta,v}:\mathsf H^\ell_\zeta(K^\varphi_w; X)\to \mathsf
W^k_\zeta(e^{-i\phi}\Bbb R^++v; X)$ is well-defined for $v\in
e^{-i\phi}\Bbb R+w$; cf.~\eqref{proj1}. If
$\eta\notin(e^{i\phi}\overline{\Bbb C^-}+\zeta)\cup
\overline{\mathcal K^\varphi_\zeta}$ then by
Proposition~\ref{effect} the image of this operator is in the
space $\mathsf H^k_\zeta(K^\varphi_v; X)\subset\mathsf
W^k_\zeta(e^{-i\phi}\Bbb R+v; X)$ . Moreover, by the same
proposition the projection operator
\begin{equation}\label{exP}
\mathsf P^s_{\eta,v}:\mathsf H^\ell_\zeta(K^\varphi_w; X)\to
\mathsf H^k_\zeta(K^\varphi_v; X),
\quad\eta\notin(e^{i\phi}\overline{\Bbb C^-}+\zeta)\cup
\overline{\mathcal K^\varphi_\zeta},
\end{equation}
satisfies the estimate
\begin{equation}\label{exP1}
\|\mathsf P^s_{\eta,v}F;\mathsf H^k_\zeta(K^\varphi_v; X)\|\leq
C\|F;\mathsf H^\ell_\zeta(K^\varphi_w; X)\|,
\end{equation}
where the constant $C$ is independent of $w$, $\phi$, and $v\in
e^{-i\phi}\overline{\Bbb R^+}+w$. In the same manner we can
consider the differential operator
\begin{equation}\label{exD1}
D^j_{(\phi)}:\mathsf W^\ell_\zeta(e^{-i\phi}\Bbb R+w;X)\to\mathsf
W^{\ell-j}_\zeta(e^{-i\phi}\Bbb R+w;X),\quad j\in\Bbb N,
\end{equation}
on the subspace $\mathsf H^\ell_\zeta(K^\varphi_w; X)$ of $\mathsf
W^\ell_\zeta(e^{-i\phi}\Bbb R+w;X)$; here $D^j_{(\phi)}$ is the
same as in~\eqref{D}. Indeed,
 by Proposition~\ref{p1.5}.(iii)  the operator of multiplication
 $$\mathscr H^\ell_w(\mathcal
K^\varphi_\zeta;X)\ni\mathcal F\mapsto (\cdot)^j\mathcal
F(\cdot)\in \mathscr H^{\ell-j}_w(\mathcal
K^\varphi_\zeta;X),\quad j\in\Bbb N$$ is bounded uniformly in
$w\in\Bbb C$. Hence the norm of the mapping
 \begin{equation}\label{3.72}
\begin{aligned}
\mathsf H^\ell_\zeta(K^\varphi_w;X)\ni
\{F_\psi:\psi\in[0,\varphi]\}&\mapsto
\{D^j_{(\psi)}F_\psi:\psi\in[0,\varphi]\}\in\mathsf
H^{\ell-j}_\zeta(K^\varphi_w;X)
\end{aligned}
\end{equation}
is also bounded uniformly in $w\in\Bbb C$;
 see the differentiation
rule~\eqref{rule} and Theorem~\ref{p1.13},(iv). Therefore the
differential operator~\eqref{exD1} maps the subspace $\mathsf
H^\ell_\zeta(K^\varphi_w;X)$ of $\mathsf
W^\ell_\zeta(e^{-i\phi}\Bbb R+w;X)$ to the subspace $\mathsf
H^{\ell-j}_\zeta(K^\varphi_w;X)\subset \mathsf
W^{\ell-j}_\zeta(e^{-i\phi}\Bbb R+w;X)$, the operator
\begin{equation}\label{exD2}
D^j_{(\phi)}:\mathsf H^\ell_\zeta(K^\varphi_w;X)\to\mathsf
H^{\ell-j}_\zeta(K^\varphi_w;X),\quad j\in\Bbb N,\
\phi\in[0,\varphi],
\end{equation}
 satisfies the uniform in $w\in\Bbb C$ and $\phi\in[0,\varphi]$ estimate
\begin{equation}\label{exD3}
\|D^j_{(\phi)}F;\mathsf H^{\ell-j}_\zeta(K^\varphi_w;X)\|\leq
C\|F;\mathsf H^\ell_\zeta(K^\varphi_w;X)\|.
\end{equation}
\section{Complex scaling of differential equations with unbounded operator
coefficients}\label{sec3}

In this section we consider linear ordinary differential equations
in spaces of analytic functions. The primary purpose here is to
motivate the study made in the previous sections and to
demonstrate the main ideas on the treatment of the complex scaling
in terms of the Hardy-Sobolev spaces and the Fredholm polynomial
operator pencils. For this reason we leave aside the important
question of complex scaling in presence of operator pencil
eigenvalues in the dual cone $\mathcal K^\varphi_\zeta$, this
aspect will be detailed elsewhere. In Subsection~\ref{GK} we
introduce the complex scaling method for differential equations
with constant operator coefficients. In Subsection~\ref{devc} we
consider equations with variable coefficients; here, for the sake
of simplicity, we restrict ourselves by the case of spaces of
positive integer orders. Two examples of applications to the
complex scaling of boundary value problems were presented in
Subsection~\ref{ss2}. Examples of applications of differential
equations with operator coefficients to boundary value problems
can be found e.g. in~\cite{R17,ref7}, some of these examples can
also be considered in context of the Hardy-Sobolev spaces.

\subsection{Differential equations with constant coefficients}\label{GK}
We shall use the same notations as in
Section~\ref{OP}. Here again we suppose that the
operator~$\mathfrak A(\lambda)$ is Fredholm for all
$\lambda\in\Bbb C$ and is invertible for at least one value of
$\lambda$. We also assume that the condition~\eqref{cndtn} is
satisfied for some $\vartheta\in(0,\pi/2)$ and $R>0$. By
$\mathfrak A(D_{(\phi)})$ we shall denote the differential
operator $\sum_{j=0}^m A_j D_{(\phi)}^{m-j}$ on the line
$e^{-i\phi}\Bbb R+w$; recall that for all $j=0,\dots,m-1$ we have
$A_j\in\mathscr B(X_j,X_0)$, $\|u\|_{j}\leq \|u\|_{j+1}$, and the
Hilbert space
   $X_{j+1}$ is dense in $X_{j}$.

For $\ell\in\Bbb R$ and $w,\zeta\in\Bbb  C$ we introduce the
Banach space
$$
\mathsf D^\ell_\zeta(e^{-i\phi}\Bbb R+w)=\bigcap_{j=0}^m \mathsf
W^{\ell-j}_\zeta(e^{-i\phi}\Bbb R+w; X_j);
$$
$$
 \|u; \mathsf
D^\ell_\zeta(e^{-i\phi}\Bbb R+w)\|=\sum_{j=0}^m \|u;\mathsf
W^{\ell-j}_\zeta(e^{-i\phi}\Bbb R+w; X_j)\|.
$$
The  operator
\begin{equation}\label{pencilD}
\mathfrak A(D_{(\phi)}): \mathsf D^\ell_\zeta(e^{-i\phi}\Bbb
R+w)\to \mathsf W^{\ell-m}_\zeta(e^{-i\phi}\Bbb R+w; X_0)
\end{equation}
is bounded for any $\zeta,w\in\Bbb C$, $\ell\in\Bbb R$, and
$\phi\in[-\pi,\pi)$. The next theorem is in essence a variant of
Theorem~2.4.1 in~\cite{R17}.

\begin{thm}\label{thm1} Suppose that the operator pencil $\lambda\mapsto \mathfrak A(\lambda)$ satisfies
the condition~\eqref{cndtn} for some $\vartheta\in(0,\pi/2)$ and
$R>0$. Let the line $e^{i\phi}\Bbb R+\zeta$, $|\phi|<\vartheta$,
be free from the spectrum of the operator pencil~$\Bbb
C\ni\lambda\mapsto\mathfrak A(\lambda)$. Then for all $\ell\in\Bbb
R$ and $w\in\Bbb C$ the operator~\eqref{pencilD} yields an
isomorphism and the estimates
$$
\begin{aligned}
\|u; \mathsf D^\ell_\zeta( e^{-i\phi}\Bbb R+w)\|\leq
c_1\|\mathfrak A(D_{(\phi)})u;\mathsf
H^{\ell-m}_\zeta(e^{-i\phi}\Bbb R+w; X_0) \|\\\leq c_2\|u; \mathsf
D^\ell_\zeta(e^{-i\phi}\Bbb R+w)\|
\end{aligned}
$$
hold. The constants $c_1$, $c_2$ does not depend on $w\in\Bbb C$
and $u\in\mathsf D^\ell_\zeta( e^{-i\phi}\Bbb R+w)$.
\end{thm}
\begin{pf} It is easily  seen that the Fourier-Laplace
transformation  implements an isometric isomorphism $\Bbb
T^\phi_{\zeta,w}:\mathfrak D^\ell_w(e^{i\phi}\Bbb
R+\zeta)\to\mathsf D^\ell_\zeta(e^{-i\phi}\Bbb R+w)$; see
Corollary~\ref{noCol} and the definition~\eqref{sp'} of the space
$\mathfrak D^\ell_w(e^{i\phi}\Bbb R+\zeta)$. Now the assertion
follows from Theorem~\ref{1},(iii) and the differentiation
rule~\eqref{rule}. \qed\end{pf} Let us introduce the scale of
Banach spaces
\begin{equation}\label{space} \mathsf D^\ell_\zeta(
K_w^\varphi)=\bigcap_{j=0}^m \mathsf H^{\ell-j}_\zeta(K_w^\varphi;
X_j);\quad \|u; \mathsf D^\ell_\zeta(K_w^\varphi)\|=\sum_{j=0}^m
\|u;\mathsf H^{\ell-j}_\zeta(K_w^\varphi; X_j)\|.
\end{equation}
 From
Theorem~\ref{p1.13} and the definitions~\eqref{sp},~\eqref{space}
of the spaces $\mathfrak D^\ell_w( \mathcal K_\zeta^\varphi)$ and
$\mathsf D^\ell_\zeta( K_w^\varphi)$ we see that the
Fourier-Laplace transformation $\Bbb T^\psi_{\zeta, w}$ yields an
isometric isomorphism between $\mathfrak D^\ell_w( \mathcal
K_\zeta^\varphi)$ and $\mathsf D^\ell_\zeta( K_w^\varphi)$. The
next proposition is a consequence of Proposition~\ref{sp prop}.
\begin{prop}
Let $\zeta,w\in\Bbb C$, $\ell\in\Bbb R$,
 and $\varphi\in(0,\pi]$. The following assertions are valid.

 (i) Every function $u\in \mathsf D^\ell_\zeta(
 K^\varphi_w)$ has boundary limits
 $u_0\in \mathsf D^\ell_\zeta(\Bbb R+w)$ and $u_\varphi\in\mathsf D^\ell_\zeta(e^{-i\varphi}\Bbb
 R+w)$
in the sense that
 \begin{equation*}
\begin{aligned}
\|(e_\zeta u)\circ\tau_{\psi,w}-(e_\zeta
u_0)\circ\tau_{0,w};\mathsf D^\ell_0(\Bbb R)\|&\to 0,\quad\psi\to
0+,
\\
\|(e_\zeta u)\circ\tau_{\psi,w}-(e_\zeta
u_\varphi)\circ\tau_{\psi,w};\mathsf D^\ell_0(\Bbb R)\|&\to
0,\quad\psi\to \varphi-;
\end{aligned}
\end{equation*}
recall that $e_\zeta:z\mapsto\exp(-i\zeta z)$ and
$\tau_{\psi,w}(t)=e^{-i\psi}t+w$.

  (ii) For all $\psi\in[0,\varphi]$ and
$u\in\mathsf D^\ell_\zeta (
 K^\varphi_w)$ the estimate
\begin{equation*}
\|u;\mathsf D^\ell_\zeta(e^{-i\psi}\Bbb R+w)\|\leq C\|u;\mathsf
D^\ell_\zeta(\mathcal
 K^\varphi_w) \|
\end{equation*}
holds, where the constant $C$ is independent of $u$, $\psi$, and
$w$.
\end{prop}

By analogy with the case of the spaces $\mathsf
H^\ell_\zeta(K^\varphi_w;X)$ and $\mathsf
W^\ell_\zeta(e^{-i\phi}\Bbb R+w;X)$ we can identify a set of
distributions $\{u_\psi:\psi\in[0,\varphi]\}\in \mathsf
D^\ell_\zeta(K^\varphi_w)$ with the correspondent element $u_\phi$
of the space $\mathsf D^\ell_\zeta(e^{-i\phi}\Bbb R+w)$,
$\phi\in[0,\varphi]$; the argument is the same, see the
explanation to the formulas~\eqref{exP},~\eqref{exD2}. Then we can
interpret $\mathsf D^\ell_\zeta(K^\varphi_w)$ as a dense subspace
of $\mathsf D^\ell_\zeta(e^{-i\phi}\Bbb R+w)$ and restrict the
operator~\eqref{pencilD} to $\mathsf D^\ell_\zeta(K^\varphi_w)$.
Since the operator $D^j_{(\phi)}$ yields the continuous
mapping~\eqref{exD2}  the operator $\mathfrak A(D_{(\phi)})$
continuously maps $\mathsf D^\ell_\zeta(K^\varphi_w)$ to the space
$ \mathsf H^{\ell-m}_\zeta(K^\varphi_w; X_0)\subset \mathsf
W^{\ell-m}_\zeta(e^{-i\phi}\Bbb R+w; X_0)$.

Theorem~\ref{T1} cited in the introductory part is a consequence
of the following
\begin{thm}\label{t1} Suppose that the operator pencil~$\mathfrak A$
meets the condition~\eqref{cndtn} for some $\vartheta\in(0,\pi/2)$
and $R>0$. Let $\varphi\in(0,\vartheta)$ and $\zeta\in \Bbb C$. If
the closed cone $\overline{\mathcal K^\varphi_\zeta}$ is free from
the spectrum of the operator pencil~$\Bbb
C\ni\lambda\mapsto\mathfrak A(\lambda)$ then the operator
\begin{equation}\label{pencilD'}
\mathfrak A(D_{(\phi)}): \mathsf D^\ell_\zeta(K^\varphi_w)\to
\mathsf H^{\ell-m}_\zeta(K^\varphi_w; X_0),\quad
\phi\in[0,\varphi],
\end{equation}
yields an isomorphism,  and the estimates
$$
\|u; \mathsf D^\ell_\zeta( K_w^\varphi)\|\leq c_1\|\mathfrak
A(D_{(\phi)})u;\mathsf H^{\ell-m}_\zeta(K_w^\varphi; X_0) \|\leq
c_2\|u; \mathsf D^\ell_\zeta( K_w^\varphi)\|
$$
are valid, where the constants $c_1$, $c_2$ are independent of
$w\in\Bbb C$, $\phi\in[0,\varphi]$, and $u\in\mathsf D^\ell_\zeta(
K_w^\varphi)$.
\end{thm}
\begin{pf} The assertion directly follows from Theorem~\ref{1}, (i) and
(ii) because the transformation $\Bbb T^\phi_{\zeta, w}:\mathfrak
D^\ell_w( \mathcal K_\zeta^\varphi)\to\mathsf D^\ell_\zeta(
K_w^\varphi)$ yields an isometric isomorphism. \qed\end{pf}

\begin{cor}\label{coltr1} Let the assumptions of Theorem~\ref{t1} be
fulfilled and let $\phi\in[0,\varphi]$. If $u\in\mathsf
D^\ell_\zeta(e^{-i\phi}\Bbb R+w)$
 satisfies the equation $\mathfrak
A(D_{(\phi)})u=F$, and  the right hand side $F$ is in the subspace
$\mathsf H^{\ell-m}_\zeta(K^\varphi_w;X_0)$ of the space $\mathsf
W^{\ell-m}_\zeta (e^{-i\phi}\Bbb R+w;X_0)$,  then the solution $u$
is in the subspace $\mathsf D^{\ell}_\zeta(K^\varphi_w)$ of the
space $\mathsf D^{\ell}_\zeta(e^{-i\phi}\Bbb R+w)$.
\end{cor}
For the sake of simplicity we restrict
 ourselves in the next theorem to the case of an integer $\ell$, $\ell\geq
 m$. Moreover, we make an additional assumption on the regularity
 of the right hand side $F$ in a neighbourhood of the point $t=0$.
 This allows us ``to localize" the problem
 to the right half-line preserving the analyticity and without recourse to the spaces of negative orders.
\begin{thm}\label{comsc}  Suppose that the operator pencil~$\mathfrak A$ meets
the condition~\eqref{cndtn} for some $\vartheta\in(0,\pi/2)$ and
$R>0$. Let $\ell$, $\ell\geq m$, be an integer number, and let the
parameters $\varphi\in(0,\vartheta)$, $\zeta\in \Bbb C$ be such
that the closed cone $\overline{\mathcal K^\varphi_\zeta}$ is free
from the spectrum of the operator pencil~$\Bbb
C\ni\lambda\mapsto\mathfrak A(\lambda)$. Assume that $F$  is a
function representable in the form $F=J+G\!\!\upharpoonright_{\Bbb
R}$, where $J \in\mathsf W^{\ell-m}_\zeta(\Bbb R^-;X_0)$ and $G\in
\mathsf H^{\ell-m}_\zeta(K^\varphi_0;X_0)$. In addition we assume
that $\chi F\in \mathsf W^\ell_0(\Bbb R;X_0)$, where $\chi\in
C^\infty(\Bbb R)$ is a compactly supported cutoff function,
$\chi=1$ in a neighbourhood of the point $t=0$.
 Then a solution $u\in \mathsf D^{\ell}_\zeta(\Bbb
R)$ of the equation~$\mathfrak A(D_t)u=F$ can be extended to an
analytic function $K^{\varphi,+}_0\ni z\mapsto u(z)\in X_{m}$ such
that: (i)  the extension meets the inclusions
 $u\in C^{\ell-j-1}\bigl(\overline{K^{\varphi,+}_0}\cup\Bbb R; X_j \bigr)$, where
$0\leq j\leq \max\{m,\ell-1\}$; (ii) for  almost all $t>0$
 the value $\|D_z^{\ell-m}u(z)-D_t^{\ell-m}u(t)\|_m$ tends to zero as $z$
goes to $t$
 by a non-tangential to $\Bbb R$ path in $
 K^{\varphi,+}_0$. Moreover, $ \mathfrak
A(D_z)u(z)=G(z)$ for $z\in K^{\varphi,+}_0$,
 and the estimate
\begin{equation}\label{estimate1}
\begin{aligned}
\sum_{j=0}^\ell\int_{e^{-i\psi}\Bbb R^+} \|e^{-i\zeta
z}D^j_zu(z)\|^2_{\ell-j}\,|dz|\leq C \bigl( \|J;\mathsf
W^{\ell-m}_\zeta(\Bbb R^-;X_0)\|^2\\+\|G;\mathsf
H^{\ell-m}_\zeta(K^\varphi_0;X_0)\|^2+\|\chi (J+G);\mathsf
W^{\ell}_0(\Bbb R;X_0)\|^2\bigr)
\end{aligned}
\end{equation}
holds, where   the constant $C$ is independent of
$\psi\in(0,\varphi)$, $J$,  $G$, and $\chi$.

\end{thm}
\begin{pf}
 Let $\rho\in C^\infty
(\Bbb R)$ be a cutoff function such that $\rho=1$ in a
neighbourhood of the point $t=0$, and $\rho\chi=\rho$. Then
$\mathfrak A(D_t)\rho u=\rho\chi F+[\mathfrak A(D_t), \rho]u $,
where $\rho\chi F=\rho F\in\mathsf W^{\ell}_\zeta (\Bbb R; X_0)$
and $\ord [\mathfrak A(D_t),\rho]\leq m-1$. The operator
$\mathfrak A(D_t):\mathsf D^s_\zeta(\Bbb R)\to \mathsf
W^{s-m}_\zeta (\Bbb R; X_0)$ yields an isomorphism for all
$s\in\Bbb R$ (Theorem~\ref{thm1}). Hence we have $\rho u\in\mathsf
D^{\ell+m}_\zeta(\Bbb R)$ and
\begin{equation}\label{n1}
\|\rho u; \mathsf D^{\ell+m}_\zeta(\Bbb R)\|\leq C\bigl( \|\chi
F;\mathsf W^\ell_0(\Bbb R; X_0)\|+\|F;\mathsf
W^{\ell-m}_\zeta(\Bbb R;X_0)\|\bigr).
\end{equation}
 By the Sobolev theorem the function $\rho u\in\mathsf D^{\ell+m}_\zeta(\Bbb R)\bigl(=\bigcap_{j=0}^m \mathsf
W^{\ell+m-j}_\zeta (\Bbb R; X_j)\bigr)$ 
 has traces $\bigl(D_t^j\rho u\bigr)(0)=\bigl
(D^j_t u \bigr)(0)\in X_m$, $j\leq \ell-1$,
 the estimates
\begin{equation}\label{n2}
\|\bigl (D^j_t u \bigr)(0); X_m\|\leq \|\rho u; \mathsf
D^{\ell+m}_\zeta(\Bbb R)\|,\quad j=0,\dots,\ell-1,
\end{equation}
hold. Then there exists an analytic function
$\overline{K^\varphi_0}\ni z\mapsto \Phi(z)\in X_m$ such that
\begin{equation}\label{4.8}
\bigl(D^j_z\Phi\bigr)(0)=\bigl(D_t^j u \bigr)(0)\in X_m,\quad
j=0,\dots,\ell-1,
\end{equation}
\begin{equation}\label{sc}
\begin{aligned}
\sum_{j=0}^{\ell+m} \int_{e^{-i\psi}\Bbb R^+}\|e^{-i\zeta
z}\bigl(D^j_z\Phi \bigr)(z)\|^2_m &\,|dz|\leq C\bigl( \|\chi
F;\mathsf W^\ell_0(\Bbb R; X_0)\|^2\\&+\|F;\mathsf
W^{\ell-m}_\zeta(\Bbb R;X_0)\|^2\bigr), \quad\psi\in[0,\varphi].
\end{aligned}
\end{equation}
 One can find a function satisfying the conditions~\eqref{4.8},~\eqref{sc} in
 the form $e^{\gamma z}\sum_{j=0}^{\ell-1} a_j z^j$, where
 $\gamma\in\Bbb C$ and $a_j\in X_m$; cf.~\eqref{n1} and~\eqref{n2}.
From Proposition~\ref{pr1} together with~\eqref{4.8} we get
 $\theta(u-\Phi)\in \cap_{j=0}^m \mathsf W^{\ell-j}_\zeta(\Bbb R^+;X_j)$
 and, in particular, $\theta(u-\Phi)\in\mathsf D^{\ell}_\zeta(\Bbb R)$;
 here $\theta$ denotes the Heaviside unit step function.
 Due to the relations~\eqref{4.8} the equality $[\mathfrak A(D_t),\theta] (u-\Phi)=0$
holds, where $[\cdot,\cdot]$ stands for the commutator
$[a,b]=ab-ba$. Therefore $\mathfrak A(D_t)\theta (u-\Phi)=\theta
(F-\mathfrak A(D_t)\Phi)=\theta (G-\mathfrak A(D_z)\Phi)$.
Moreover, it is easy to see that $D_t^j(G-\mathfrak
A(D_z)\Phi)(0)=0$ for $j=0,\dots,\ell-m-1$. The embedding result
of Theorem~\ref{THMemb} states that  $G\in
C^{\ell-m-1}(\overline{K^\varphi_0};X_0)$. Thus we have
$D_{(\psi)}^j(G-\mathfrak A(D_z)\Phi)(0)=0$ for
$j=0,\dots,\ell-m-1$ and $\psi\in[0,\varphi]$. This together with
Proposition~\ref{pr1} and the estimate~\eqref{sc} gives
\begin{equation}\label{ema1}
\begin{aligned}
\|\theta (G-\mathfrak A(D_z)\Phi); \mathsf H^{\ell-m}_\zeta
(K_0^\varphi;X_0)\|^2\leq C \bigl( \|J;\mathsf
W^{\ell-m}_\zeta(\Bbb R^-;X_0)\|^2\\+\|G;\mathsf
H^{\ell-m}_\zeta(K^\varphi_0;X_0)\|^2+\|\chi F;\mathsf
W^{\ell}_0(\Bbb R;X_0)\|^2\bigr).
\end{aligned}
\end{equation}
As a consequence of Theorem~\ref{t1} (see also
Corollary~\ref{coltr1}) we get the estimate
\begin{equation}\label{ema2}
\|\theta (u-\Phi); \mathsf D^\ell_\zeta (K^\varphi_0)\|^2\leq c
\|\theta (G-\mathfrak A(D_z)\Phi); \mathsf H^{\ell-m}_\zeta
(K_0^\varphi;X_0)\|^2.
\end{equation}
Recall that $\mathsf D^\ell_\zeta (K^\varphi_\zeta)=\cap_{j=0}^m
\mathsf H^{\ell-j}_\zeta(K^\varphi_0;X_j)$. Hence by
Theorem~\ref{THMemb} we have the inclusions $\theta (u-\Phi)\in
C^{\ell-j-1}(\overline {K^\varphi_0}; X_j)$, where $0\leq j\leq
\max\{m,\ell-1\}$; by the Sobolev embedding theorem with the same
restrictions on $j$ we have $u\in C^{\ell-j-1}(\Bbb R; X_j)$. This
together with the relations~\eqref{4.8} proves the property (i).
We also note that $D^{\ell-m}_t(\theta u-\theta\Phi)$ is in the
subspace $\mathsf H^{0}_\zeta(K^\varphi_0;X_m)$ of the space
$\mathsf W^{0}_\zeta(\Bbb R;X_m)$; see~\eqref{3.72}. The property
(ii) follows from Theorem~\ref{p3}, (i). It is clear that
$\mathfrak A(D_z)u(z)=G(z)$ for all $z\in K^{\varphi,+}_0$. The
estimate~\eqref{estimate1} is readily apparent from the
estimates~\eqref{ema1},~\eqref{ema2} and~\eqref{sc}. \qed\end{pf}

\subsection{Differential equations with variable
coefficients}\label{devc}
 The aim of this section is to give an analog of
Theorem~\ref{comsc} for the case of equations with variable
operator coefficients. We shall consider the equation
\begin{equation}\label{VC}
\mathfrak A(D_t)u(t)-\sum_{j=0}^m Q_j(t)D_t^{m-j} u(t)=F(t),\quad
t\in\Bbb R.
\end{equation}
Here  $\mathfrak A(D_t)$ is the same as before, the coefficients
$Q_0,\dots,Q_m$ are operator functions $\Bbb R\ni t\mapsto
Q_j(t)\in \mathscr B(X_j,X_0)$ satisfying the following
conditions: \vspace{0.2cm}

\noindent{\it {\rm i.} for a large $T>0$ and some $\alpha>0$ the
coefficients $Q_0,\dots,Q_m$ can be extended to holomorphic
operator functions $ \overline{ K^{\alpha,+}_T}\ni z\mapsto
Q_j(z)\in\mathscr B(X_j,X_0)$;

\vspace{0.2cm} \noindent {\rm ii.}  for all $n=0,1,\dots$ the
values $\|D_z^n Q_j(z); \mathscr B(X_j,X_0)\|$, $j=0,\dots, m$,
uniformly tend to zero as $z\to\infty$, $z\in\overline{
K^{\alpha,+}_T}$.}\vspace{0.2cm}

If  the operator functions $\Bbb R\ni t\mapsto D^n_t Q_j(t)\in
\mathscr B(X_j,X_0)$, $j\leq m$, are bounded for all $n=0,1,\dots$
then the operator $ (\mathfrak A(D_t)-\sum_{j=0}^m
Q_j(t)D_t^{m-j}):\mathsf D^\ell_\zeta (\Bbb R)\to \mathsf
W^{\ell-m}_\zeta(\Bbb R; X_0)$ of the equation~\eqref{VC} is
continuous for any $\ell\in\Bbb R$  and $\zeta\in\Bbb C$. For
simplicity we restrict ourselves to the case of an integer
$\ell\geq m$.

We start with the subsidiary equation
\begin{equation}\label{MP}
\mathfrak A(D_{(\phi)})U_\phi(z)- \sum_{j=0}^m Q_j(z)
D^{m-j}_{(\phi)}\mathsf P^{\ell-j}_{\eta, w}U_\phi(z)=G_\phi(z),\
\ \quad z\in e^{-i\phi}\Bbb R+w,
\end{equation}
where  $w\in \overline{ K^{\alpha,+}_T}$ and $\phi\in[0,\alpha]$.
Here $\mathsf P^{\ell-j}_{\eta, w}$ is the projection
operator~\eqref{proj1} with the  parameter $\eta\in
\bigcap_{\phi\in[0,\alpha]}e^{i\phi}\Bbb C^++\zeta$. The
coefficients $Q_0,\dots,Q_m$
 are defined on the set $e^{-i\phi}\Bbb
R+w\cap\overline{ K^{\alpha,+}_T}$ as the analytic extensions of
the correspondent coefficients of the equation~\eqref{VC}. Without
loss of generality we can assume that $Q_j(z)=0$ for $z\notin
\overline{ K^{\alpha,+}_T}$ because $\supp\{ \mathsf P^{-j}_{\eta,
w} U_\phi\}\subset \overline{ K^{\alpha,+}_T}$. Let us note that
the  norms of the operators
$$
D^{m-j}_{(\phi)}\mathsf P^{\ell-j}_{\eta, w}: \mathsf
W_\zeta^{\ell-j}( e^{-i\phi}\Bbb R+w;X_j)\to\mathsf
W_\zeta^{\ell-m}( e^{-i\phi}\Bbb R^++w;X_j), \quad j\leq m,
$$
are bounded uniformly in $w\in\Bbb C$; see~\eqref{la}.
 Since
$Q_0(z),\dots,Q_m(z)$ satisfy the condition ii  the mappings
$$
\mathsf W_\zeta^{\ell-m}( e^{-i\phi}\Bbb R^++w;X_j)\ni F\mapsto
Q_j F\in \mathsf W_\zeta^{\ell-m}( e^{-i\phi}\Bbb R^++w;X_0),
\quad j\leq m,
$$
tend to zero in the operator norms as $w\to\infty$, $w\in
\overline{ K^{\alpha,+}_T}$. Hence the operator
\begin{equation}\label{pertD}
\mathfrak Q^+_w(z, D_{(\phi)})=\sum_{j=0}^mQ_j(z)
D^{m-j}_{(\phi)}\mathsf P^{\ell-j}_{\eta, w},
\end{equation}
\begin{equation}\label{pert}
\mathfrak Q^+_w(z, D_{(\phi)}): \mathsf D^{\ell}_\zeta(
e^{-i\phi}\Bbb R+w)\to\mathsf W^{\ell-m}_\zeta( e^{-i\phi}\Bbb
R^++w;X_0)
\end{equation}
tends to zero in the operator norm as $w\to\infty$, $w\in
\overline{ K^{\alpha,+}_T}$.

\begin{thm}\label{thm2} Let the coefficients $Q_0,\dots,Q_m$ of the equation~\eqref{VC}
meet the conditions i, ii in some cone $\overline{
K^{\alpha,+}_T}$, and  let $w\in \overline{ K^{\alpha,+}_T}$ be a
complex number with  a sufficiently large module. Suppose that
$\zeta\in \Bbb C$, $\phi\in[0,\alpha]$ and $\phi<\vartheta$; here
$\vartheta$ is the angle from the condition~\eqref{cndtn}. If the
line $e^{i\phi}\Bbb R+\zeta$ is free from the spectrum of the
pencil~$\Bbb C\ni\lambda\mapsto \mathfrak A(\lambda)$ then the
operator
\begin{equation}\label{thm2 op}
\mathfrak A(D_{(\phi)})-{\mathfrak Q^+_w}(z, D_{(\phi)}):\mathsf
D^{\ell}_\zeta( e^{-i\phi}\Bbb R+w)\to\mathsf W^{\ell-m}_\zeta(
e^{-i\phi}\Bbb R+w;X_0)
\end{equation}
of the subsidiary  equation~\eqref{MP} yields an isomorphism for
all $\ell\geq m$, $\ell\in\Bbb Z$.
\end{thm}
\begin{pf}  By Theorem~\ref{thm1} the
operator~\eqref{pencilD} with constant coefficients is invertible,
and the norm of the inverse operator is bounded uniformly in
$w\in\Bbb C$. Let us denote the inverse operator by $\mathfrak R$.
Since the norm of the operator~\eqref{pert} tends to zero and the
norm of $\mathfrak R$ remains bounded as $w\to\infty$, $w\in
\overline{ K^{\alpha,+}_T}$ , the norm of the composition
$$
{\mathfrak Q^+_w}(z, D_{(\phi)})\mathfrak R:\mathsf
W^{\ell-m}_\zeta( e^{-i\phi}\Bbb R+w;X_0)\to\mathsf
W^{\ell-m}_\zeta( e^{-i\phi}\Bbb R^++w;X_0)
$$
is strictly less than $1$ as far as  $w\in \overline{
K^{\alpha,+}_T}$ and $|w|$ is a sufficiently large positive
number. Then the operator~\eqref{thm2 op} with variable
coefficients is invertible, its inverse operator is defined as the
Neumann operator series $$ \mathfrak
R\sum_{n=0}^{+\infty}\bigl({\mathfrak Q^+_w}(z,
D_{(\phi)})\mathfrak R\bigr)^n:\mathsf W^{\ell-m}_\zeta(
e^{-i\phi}\Bbb R+w;X_0)\to\mathsf D^{\ell}_\zeta( e^{-i\phi}\Bbb
R+w).$$ \qed\end{pf}

Now we consider the operator~\eqref{thm2 op} on the subspace
$\mathsf D^{\ell}_\zeta( K^\varphi_w)$ of $\mathsf D^{\ell}_\zeta(
e^{-i\phi}\Bbb R+w)$, $\phi\in[0,\varphi]$. Let us verify that the
operator $\mathfrak A(D_{(\phi)})-{\mathfrak Q^+_w}(z,
D_{(\phi)})$ continuously maps $\mathsf D^{\ell}_\zeta(
K^\varphi_w)$ to the space $\mathsf H^{\ell-m}_\zeta( K^\varphi_w;
X_0)\subset \mathsf W^{\ell-m}_\zeta(e^{-i\phi}\Bbb R+w; X_0)$. We
have already studied the operator $\mathfrak A(D_{(\phi)})$ on
$\mathsf D^{\ell}_\zeta( K^\varphi_w)$, therefore we only need to
check the continuity of the operator ${\mathfrak Q^+}(z,
D_{(\phi)}):\mathsf D^{\ell}_\zeta( K^\varphi_w)\to\mathsf
H^{\ell-m}_\zeta( K^\varphi_w; X_0)$.

 First of all we note that the
norms of the operators
\begin{equation}\label{dop}
D^{m-j}_{(\phi)}\mathsf P^{\ell-j}_{\eta, w}: \mathsf
H_\zeta^{\ell-j}( K^\varphi_w;X_j)\to\mathsf H_\zeta^{\ell-m}(
 K^\varphi_w;X_j), \quad j\leq m,
\end{equation}
are bounded uniformly in $w\in\Bbb C$; see~\eqref{exP1}
and~\eqref{exD3}. The operators~\eqref{dop} have their images in
$\mathsf H_\zeta^{\ell-m}(
 K^\varphi_w;X_j)\cap \mathsf W^{\ell-m}_\zeta(e^{-i\phi}\Bbb R^++w; X_j)$.
  Since the coefficients $Q_j(z)$ are subjected
to the condition ii, the operator norms of the mappings
$$
Q_jD^{m-j}_{(\phi)}\mathsf P^{\ell-j}_{\eta, w}: \mathsf
H_\zeta^{\ell-j}( K^\varphi_w;X_j)\to\mathsf H_\zeta^{\ell-m}(
 K^\varphi_w;X_0), \quad j\leq m,
$$
tend to zero  as $w\to\infty$, $w\in \overline{ K^{\alpha,+}_T}$.
Hence the continuous operator
\begin{equation}\label{pert'}
\mathfrak Q^+_w(z, D_{(\phi)}): \mathsf D^{\ell}_\zeta(
K^\varphi_w)\to\mathsf W^{\ell-m}_\zeta( K^\varphi_w;X_0)
\end{equation}
tends to zero in the operator norm as $w\to\infty$, $w\in
\overline{ K^{\alpha,+}_T}$.

\begin{thm} \label{lt} Let the coefficients of the equation~\eqref{VC}
meet the conditions i, ii in some cone $\overline{
K^{\alpha,+}_T}$, and  let $w\in \overline{ K^{\alpha,+}_T}$ be a
complex number with  a sufficiently large module.  Suppose that
$\zeta\in \Bbb C$, $\varphi\in(0,\alpha]$ and $\varphi<\vartheta$;
here $\vartheta$ is the angle from the condition~\eqref{cndtn}. If
the cone $\mathcal K^\varphi_\zeta$ is free from the spectrum of
the operator pencil $\Bbb C\ni\lambda\mapsto \mathfrak A(\lambda)$
then the operator
$$
\mathfrak A(D_{(\phi)})-{\mathfrak Q^+_w}(z, D_{(\phi)}):\mathsf
D^{\ell}_\zeta( K^\varphi_w)\to\mathsf H^{\ell-m}_\zeta(
K^\varphi_w;X_0)
$$
of the subsidiary equation~\eqref{MP} yields an isomorphism.
\end{thm}
\begin{pf}  The proof is similar to the proof of Theorem~\ref{thm2}. Indeed,
by Theorem~\ref{t1} the operator $\mathfrak A(D_{(\phi)})$
implements an isomorphism~\eqref{pencilD'} and its inverse
operator is bounded uniformly in $w\in\Bbb C$. It remains to note
that under the assumptions of the theorem the norm of the
operator~\eqref{pert'} is sufficiently small. \qed\end{pf}

Now we are in position to extend the results of
Theorem~\ref{comsc} to the case of equations with variable
coefficients. The following theorem presents a more general result
than Theorem~\ref{T2} cited in the introductory part.
\begin{thm} \label{---}Let the conditions i, ii on the coefficients of the equation~\eqref{VC}
be fulfilled in some cone $\overline{ K^{\alpha,+}_T}$. Suppose
that the polynomial operator pencil~$\Bbb C\ni\lambda\mapsto
\mathfrak A(\lambda)$ satisfies the condition~\eqref{cndtn}  for
some $\vartheta\in(0,\pi/2)$ and $R>0$. Let the angle
$\varphi\in(0,\alpha]$, $\varphi<\vartheta$, and the vertex
$\zeta\in\Bbb C$ of the cone $\mathcal K^\varphi_\zeta$ be such
that the closed cone $\overline{\mathcal K^\varphi_\zeta}$ is free
from the spectrum of the pencil $\Bbb C\ni\lambda\mapsto \mathfrak
A(\lambda)$.   Assume that $\ell\geq m$, $\ell\in\Bbb Z$, and the
right hand side $F\in\mathsf W^{\ell-m}_\zeta(\Bbb R;X_0)$ of the
equation~\eqref{VC} is representable in the form
$F=J+G\!\!\upharpoonright_{\Bbb R}$, where $J \in\mathsf
W^{\ell-m}_\zeta(\Bbb R^-+T;X_0)$, $G\in \mathsf
H^{\ell-m}_\zeta(K^\varphi_T;X_0)$, and $T$ is a sufficiently
large positive number. In addition we assume that $\chi F\in
\mathsf W^\ell_0(\Bbb R;X_0)$, where $\chi\in C^\infty(\Bbb R)$ is
a compactly supported cutoff function, $\chi=1$ in a neighbourhood
of the point $t=R$.
 Then a solution $u\in \mathsf D^{\ell}_\zeta(\Bbb
R)$ of the equation~\eqref{VC} can be extended to an analytic
function $K^{\varphi,+}_R\ni z\mapsto u(z)\in X_{m}$ such that:
(i)  the extension meets the inclusions
 $u\in C^{\ell-j-1}\bigl(\overline{K^{\varphi,+}_R}\cup\Bbb R; X_j
 \bigr)$,
$ j=0,\dots,\max\{m,\ell-1\}$; (ii) for  almost all $t>R$
 the value $\|D_z^{\ell-m}u(z)-D_t^{\ell-m}u(t)\|_m$ tends to zero as $z$
goes to $t$
 by a non-tangential to $\Bbb R$ path in $
 K^{\varphi,+}_R$. Moreover, the extension satisfies the following equation and  estimate
\begin{equation}\label{inK} \mathfrak A(D_z)u(z)- \sum_{j=0}^m Q_j(z)
D^{m-j}_zu(z)=G(z),\quad z\in K^{\varphi,+}_T,
\end{equation}
\begin{equation}\label{estimateX`}
\begin{aligned}
\sum_{j=0}^\ell\int_{e^{-i\psi}\Bbb R^+}& \|e^{-i\zeta
z}D^j_zu(z+T)\|^2_{\ell-j}\,|dz|\leq C \bigl( \|\chi u;\mathsf
D^{\ell}_0(\Bbb R)\|^2\\&+\|G;\mathsf
H^{\ell-m}_\zeta(K^\varphi_T;X_0)\|^2+\|\chi (J+G);\mathsf
W^{\ell}_0(\Bbb R;X_0)\|^2\bigr),
\end{aligned}
\end{equation}
 where   the constant $C$ is independent of
$\psi\in(0,\varphi)$, $J$,  $G$, and $\chi$.
\end{thm}
\begin{pf}
 Let $\rho\in C^\infty
(\Bbb R)$ be a cutoff function such that $\rho=1$ in a
neighbourhood of the point $t=T$, $\rho(t)=0$ for $t<T-1/2$, and
$\rho\chi=\rho$. Then
$$\bigl(\mathfrak A(D_t)-\mathfrak Q^+_{T-1}(t, D_t)\bigr )\rho u=\rho\chi
F+[\mathfrak A(D_t)-\mathfrak Q^+_{T-1}(t, D_t), \rho]\chi u, $$
 where $\rho\chi F=\rho F\in\mathsf
W^{\ell}_\zeta (\Bbb R; X_0)$ and $\ord [\mathfrak
A(D_t)-\mathfrak Q^+_{T-1}(t, D_t),\rho]\leq m-1$. The operator
$\mathfrak A(D_t)-\mathfrak Q^+_{T-1}(t, D_t):\mathsf
D^s_\zeta(\Bbb R)\to \mathsf W^{s-m}_\zeta (\Bbb R; X_0)$ yields
an isomorphism for all $s\geq m$, $s\in\Bbb Z$
(Theorem~\ref{thm2}). Hence we have $\rho u\in\mathsf
D^{\ell+m}_\zeta(\Bbb R)$ and
\begin{equation}\label{n1`}
\|\rho u; \mathsf D^{\ell+m}_\zeta(\Bbb R)\|\leq C\bigl( \|\chi
F;\mathsf W^\ell_0(\Bbb R; X_0)\|+\|\chi u;\mathsf D^{\ell}_0(\Bbb
R)\|\bigr).
\end{equation}
 The function $\rho u\in\bigcap_{j=0}^m \mathsf
W^{\ell+m-j}_\zeta (\Bbb R; X_j)$ 
 has traces satisfying
 the estimates
\begin{equation}\label{n2`}
\|\bigl (D^j_t u \bigr)(0); X_m\|\leq \|\rho u; \mathsf
D^{\ell+m}_\zeta(\Bbb R)\|,\quad j=0,\dots,\ell-1.
\end{equation}
Let $\Phi$ denote an analytic function $\overline{K^\varphi_T}\ni
z\mapsto \Phi(z)\in X_m$ such that
\begin{equation}\label{4.8`}
\bigl(D^j_z\Phi\bigr)(T)=\bigl(D_t^j u \bigr)(T)\in X_m,\quad
j=0,\dots,\ell-1,
\end{equation}
\begin{equation}\label{sc`}
\begin{aligned}
\sum_{j=0}^{\ell+m} \int_{e^{-i\psi}\Bbb R^+}\|&e^{-i\zeta
z}\bigl(D^j_z\Phi \bigr)(z+T)\|^2_m \,|dz|\\&\leq C\bigl( \|\chi
F;\mathsf W^\ell_0(\Bbb R; X_0)\|^2+\|\chi u;\mathsf
D^{\ell}_0(\Bbb R)\|^2\bigr), \quad\psi\in[0,\varphi].
\end{aligned}
\end{equation}
By Proposition~\ref{pr1} we have
 $\theta(\cdot-T)(u-\Phi)\in \cap_{j=0}^m \mathsf W^{\ell-j}_\zeta(\Bbb
 R^++T;X_j)$.
 The equalities
 $$ [\mathfrak A(D_t)-\mathfrak Q^+_{T}(t, D_t),\theta(\cdot-T)] (u-\Phi)=0,$$
$$\bigl(\mathfrak A(D_t)-\mathfrak Q^+_{T}(t, D_t)\bigr)\theta(\cdot-T) (u-\Phi)=\theta(\cdot-T) (G-\mathfrak
A(D_t)\Phi+\mathfrak Q^+_{T}(t, D_t)\Phi)$$ hold. Moreover, one
can easily see that $D_t^j(G-\mathfrak A(D_t)\Phi+\mathfrak
Q^+_{T}(t, D_t)\Phi)(T)=0$ for $j=0,\dots,\ell-m-1$. By
Theorem~\ref{THMemb} we have  $G\in
C^{\ell-m-1}(\overline{K^\varphi_T};X_0)$, thus for
$\psi\in[0,\varphi]$ and $j=0,\dots,\ell-m-1$ we get
$$D_{(\psi)}^j(G-\mathfrak A(D_t)\Phi+\mathfrak Q^+_{T}(t,
D_t)\Phi)|_{t=T}=0.$$ This together with Proposition~\ref{pr1} and
the estimate~\eqref{sc`} gives
\begin{equation}\label{ema1`}
\begin{aligned}
\|\theta(\cdot-T) \bigl(G-\mathfrak A(D_t)\Phi+\mathfrak
Q^+_{T}(t, D_t)\Phi\bigr); \mathsf H^{\ell-m}_\zeta
(K_T^\varphi;X_0)\|^2\\\leq C \bigl( \|\chi u;\mathsf
D^{\ell}_0(\Bbb R)\|^2+\|G;\mathsf
H^{\ell-m}_\zeta(K^\varphi_T;X_0)\|^2+\|\chi F;\mathsf
W^{\ell}_0(\Bbb R;X_0)\|^2\bigr).
\end{aligned}
\end{equation}
As a consequence of Theorem~\ref{thm2} we get the estimate
\begin{equation}\label{ema2`}
\begin{aligned}
\|\theta&(\cdot-T) (u-\Phi); \mathsf D^\ell_\zeta
(K^\varphi_T)\|^2\\&\leq c \|\theta(\cdot-T) \bigl(G-\mathfrak
A(D_t)\Phi+\mathfrak Q^+_{T}(t, D_t)\Phi\bigr); \mathsf
H^{\ell-m}_\zeta (K_T^\varphi;X_0)\|^2.
\end{aligned}
\end{equation}
By Theorem~\ref{THMemb} the inclusion $\theta(\cdot-T) (u-\Phi)\in
C^{\ell-j-1}(\overline {K^\varphi_T}; X_j)$ is fulfilled for
$j=0,\dots,\max\{m,\ell-1\}$; by the Sobolev theorem for $u\in
\mathsf D^\ell_\zeta (\Bbb R)$ we have $u\in C^{\ell-j-1}(\Bbb R;
X_j)$ with the same restrictions on $j$. This together with the
relations~\eqref{4.8`} proves the property (i). The property (ii)
follows from Theorem~\ref{p3}, (i). It is clear that $u$ satisfies
the equation~\eqref{inK}. The estimate~\eqref{estimateX`} is a
consequence of the estimates~\eqref{ema1`},~\eqref{ema2`}
and~\eqref{sc`}. \qed\end{pf}
\begin{rem}\label{r---} In Theorem~\ref{---} the assumptions on the right hand
side $F\in \mathsf W^{\ell-m}_\zeta(\Bbb R; X_0) $ of the
equation~\eqref{VC} are a priori fulfilled for all $T>T_0$ if
 the function  $F$ can be extended to an analytic
function $\overline{K^{\varphi,+}_{T_0}}\ni z\mapsto F(z)\in X_0$
satisfying the uniform in $\psi\in[0,\varphi]$ estimate
\begin{equation}\label{class}
\sum_{j=0}^{\ell-m}\int_{e^{-i\psi}\Bbb R^+}\|e^{-i\zeta z}D_z^j
F(z+T_0)\|_0^2\,|dz|\leq Const.
\end{equation}
Indeed, from the estimate~\eqref{class} it follows that the
functions $K^{\varphi,+}_{T_0}\ni z\mapsto D_z^j F(z)\in X_0$,
$j=0,\dots,\ell-m$, extended to $K^{\varphi,-}_{T_0}$ by zero are
in the class $\mathscr H^0_{-\zeta}(\mathcal
K^{-\varphi}_{T_0};X_0)$; here $\mathcal
K^{-\varphi}_{T_0}=K^\varphi_{T_0}$. Then by Proposition~\ref{pr2}
we have
\begin{equation*}
\sum_{j=0}^{\ell-m}\int_{e^{-i\psi}\Bbb R^+}\|e^{-i\zeta z}D_z^j
F(z+T)\|_0^2\,|dz|\leq C,\quad T>T_0,\  \psi\in [0,\varphi].
\end{equation*}
We find  an entire function $\Bbb C\ni z\mapsto \Phi(z)\in X_0$
such that $D_z^j\Phi(T)=D_z^j F(T)$ for all $j=0,\dots,\ell-m-1$,
and
\begin{equation*}
\sum_{j=0}^{\ell-m}\int_{e^{-i\psi}\Bbb R^-}\|e^{-i\zeta z}D_z^j
\Phi(z+T)\|_0^2\,|dz|\leq C,\quad \psi\in [0,\varphi].
\end{equation*}
By setting $G=\Phi$ on $K^{\varphi,-}_T$ and $G=F$ on
$K^{\varphi,+}_T$ we define $G\in \mathsf
H^{\ell-m}_\zeta(K^\varphi_{T};X_0)$. It is clear that $J=F-G\in
\mathsf W^{\ell-m}_\zeta(\Bbb R^-+ T; X_0)$.

\end{rem}

\section{Appendix. Proof of Proposition~\ref{p2.6}, Lemma~~\ref{ApLemma2}}

\begin{pf*}{\bf PROOF of Proposition~\ref{p2.6}.}\label{proofpage} (i) Here we prove
the representation
\begin{equation}\label{Acau}
 \mathcal F(\lambda)=\int_{\partial \mathcal K^{\varphi,+}_\zeta}\frac
{e^{iw(\mu-\lambda)}(\mu-\eta)^s \mathcal F(\mu)}{2\pi i
(\lambda-\eta)^s(\mu-\lambda)}\, d\mu,\quad \lambda\in \mathcal
K^{\varphi,+}_\zeta,\ \eta\notin \overline{\mathcal
K^{\varphi,+}_\zeta},\ s\leq\ell,
\end{equation}
for a function $\mathcal F\in \mathscr H^\ell_w(\mathcal
K^\varphi_\zeta;X)$, where $\mathcal K^{\varphi,+}_\zeta
=\{\lambda\in\mathcal K^\varphi_\zeta:\Im\lambda>\Im\zeta\}$ and
 $\partial \mathcal K^{\varphi,+}_\zeta$ is the boundary
of $\mathcal K^{\varphi,+}_\zeta$, the function $(\cdot-\eta)^s$
is analytic in $\overline{\mathcal K^{\varphi,+}_\zeta}$;
cf.~\eqref{cau}. The proof of the second
representation~\eqref{1cau} in Proposition~\ref{p2.6},(i) is
similar.

Let us define the function $\mathcal G\in \mathscr H^0_0(\mathcal
K^\varphi_0;X)$ by the equality
\begin{equation}\label{G}
\mathcal
G(\lambda-\zeta)=\left\{%
\begin{array}{ll}
       \exp\{iw\lambda\}(\lambda-\eta)^s \mathcal
F(\lambda), \quad\lambda\in \mathcal K^{\varphi,+}_\zeta ; \\
       \exp\{iw\lambda\}(\lambda-\tau)^s \mathcal
F(\lambda), \quad\lambda\in \mathcal K^{\varphi,-}_\zeta. \\
\end{array}%
\right.
\end{equation}
 To prove the representation \eqref{Acau} it suffices to
establish that
\begin{equation}\label{cau1}
\mathcal G(\lambda)=\frac 1 {2\pi
i}\left(\int_0^{+\infty}\frac{\mathcal G(\xi)}{\xi-\lambda} \,d
\xi-\int_0^{+\infty}\frac{\mathcal
G(e^{i\varphi}\xi)}{e^{i\varphi}\xi-\lambda} e^{i\varphi}\,d
\xi\right)\quad\lambda\in \mathcal K^{\varphi,+}_0.
\end{equation}

By the Cauchy integral theorem we have
\begin{equation}\label{cau2}
\mathcal G(\lambda)=\frac 1 {2\pi i}\int_R^{R+1} da\oint_{\mathscr
C(a, \psi,\phi)}\frac{\mathcal G(\mu)}{\mu-\lambda}\,d \mu,
\end{equation}
where  $R>0$, the contour integration runs anticlockwise along the
closed path
\begin{equation*}
\begin{aligned}
&\mathscr C(a,\psi,\phi)=\{\mu :\mu=a
e^{i\vartheta},\vartheta\in[\psi,\phi]\}\cup\{\mu:\mu=\xi
e^{i\phi},1/a\leq \xi\leq a\}\\&\cup  \{\mu:\mu=
e^{i\vartheta}/a,\vartheta\in[\psi,\phi]\}\cup\{\mu:\mu=\xi
e^{i\psi}, 1/a\leq \xi\leq a\},\ 0<\psi<\phi<\varphi,
\end{aligned}
\end{equation*}
 and $\lambda$ is inside of the contour
$\mathscr C(R,\psi,\phi)(\subset \mathcal K^{\varphi,+}_0)$. We
will show that for all $\lambda\in \mathcal K^{\varphi,+}_0$ the
equality
\begin{equation}\label{apL1}
\int_0^{+\infty}\frac{\mathcal G(\xi)}{\xi-\lambda} \,d
\xi-\int_0^{+\infty}\frac{\mathcal
G(e^{i\varphi}\xi)}{e^{i\varphi}\xi-\lambda} e^{i\varphi}\,d
\xi=\lim\limits_{\substack{\psi\to 0+\\
\phi\to\varphi-}}\lim\limits_{R\to+\infty}\int_R^{R+1}
da\oint_{\mathscr C(a, \psi,\phi)}\frac{\mathcal
G(\mu)}{\mu-\lambda}\,d \mu
\end{equation}
is fulfilled, where the limits are taken in the space $X$. The
equalities~\eqref{cau2}, \eqref{apL1} prove the
representation~\eqref{cau1} (and consequently the
representation~\eqref{Acau}).

Let us demonstrate~\eqref{apL1}.  At first we estimate the norm in
$X$ of the integral
$$
\int_R^{R+1}\int_{\{\mu :\mu=a
e^{i\vartheta},\vartheta\in[\psi,\phi]\}}\frac{\mathcal
G(\mu)}{\mu-\lambda}\,d \mu\, da=\int_R^{R+1}\int_\psi^\phi
\frac{\mathcal G( ae^{i\vartheta})}{ae^{i\vartheta}-\lambda} i a
e^{i\vartheta}\,d \vartheta\, da.
$$ Interchanging the order of integration and applying the Cauchy-Schwarz inequality, we
get
\begin{equation}\label{cau3}
\begin{aligned}
&\left\|\int_R^{R+1}\int_\psi^\phi \frac{\mathcal G(
ae^{i\vartheta})a}{ae^{i\vartheta}-\lambda}   \,d \vartheta\,
da\right\|^2
\leq (\phi-\psi)\int_\psi^\phi
\left(\int_R^{R+1}\left\|\frac{\mathcal G( ae^{i\vartheta})a
}{ae^{i\vartheta}-\lambda}\right\|\, da\right)^2  \,d \vartheta
\\
&\leq (\phi-\psi)\int_\psi^\phi \left(\int_R^{R+1}\left|\frac{a
}{ae^{i\vartheta}-\lambda}\right|^2\, da\right)
\left(\int_R^{R+1}\|\mathcal G( ae^{i\vartheta})\|^2\, da\right)
\,d \vartheta
\\
&\leq C \int_\psi^\phi \int_R^{R+1}\|\mathcal G(
ae^{i\vartheta})\|^2\, da \,d \vartheta.
\end{aligned}
\end{equation}
In the same way we derive
\begin{equation}\label{cau4}
\begin{aligned}
&\left\|\int_R^{R+1}\int_{\{\mu :\mu=
e^{i\vartheta}/a,\vartheta\in[\psi,\phi]\}}\frac{\mathcal
G(\mu)}{\mu-\lambda}\,d \mu\,
da\right\|^2\leq\left(\int_R^{R+1}\int_\psi^\phi
\left\|\frac{\mathcal G(
e^{i\vartheta}/a)/a}{e^{i\vartheta}/a-\lambda}\right\| \,d
\vartheta\, da\right)^2
\\
&\leq (\phi-\psi)\int_\psi^\phi
\left(\int_R^{R+1}|{e^{i\vartheta}/a-\lambda}|^{-2}\, da\right)
\left(\int_R^{R+1}\|\mathcal G( e^{i\vartheta}/a)\|^2 a^{-2}\,
da\right) \,d \vartheta
\\
&\leq C \int_\psi^\phi \int_{1/(R+1)}^{1/R}\|\mathcal G(
ae^{i\vartheta})\|^2\, da \,d \vartheta.
\end{aligned}
\end{equation}
Since $\mathcal G\in\mathscr H^0_0(\mathcal K^\varphi_0;X)$, the
integrals $\int_R^{R+1}\|\mathcal G( ae^{i\vartheta})\|^2\,da$,
$\int_{1/(R+1)}^{1/R}\|\mathcal G( ae^{i\vartheta})\|^2\,da$ are
uniformly bounded and tend to zero as $R\to+\infty$. Consequently,
the right hand sides of the inequalities \eqref{cau3} and
\eqref{cau4} tend to zero as $R\to +\infty$.

Let us consider the integral
$$
I(\psi,R)=\int_R^{R+1}\,da\int_{\{\mu:\mu=\xi e^{i\psi}, 1/a\leq
\xi\leq a\}}\frac{\mathcal G(\mu)}{\mu-\lambda} \,d\mu.
$$
For brevity, in the next formula we write $[\dots]$ instead of the
integrand
$$
{\mathcal G(\xi e^{i\psi})}/(\xi e^{i\psi}-\lambda).
$$
Interchanging the order of integration, we arrive at the formula
$$
I(\psi,R)=\int_R^{R+1}\,da\int_{1/a}^{ a}[\dots] \,d \xi
$$
$$
=\int_{1/R}^R\,d\xi\int_R^{R+1}\,da[\dots]+\int_{1/(R+1)}^{1/R}\,d\xi\int_{1/\xi}^{R+1}\,da[\dots]+\int_R^{R+1}\,d\xi\int_R^\xi\,da[\dots]
$$
$$
=\int_{1/R}^R\,d\xi[\dots]+\int_{1/(R+1)}^{1/R}\,d\xi(R+1-1/\xi)[\dots]+\int_R^{R+1}\,dt
(t-R)[\dots].
$$
Note that the last two integrals tend in $X$ to zero as
$R\to+\infty$. Indeed, the estimates
$$
\left\|\int_{1/(R+1)}^{1/R}(R+1-1/\xi){\mathcal G(\xi
e^{i\psi})}(\xi e^{i\psi}-\lambda)^{-1}\,d\xi\right\|^2\leq
C\int_{1/(R+1)}^{1/R}\|\mathcal G(\xi e^{i\psi})\|^2\,d\xi,
$$
$$
\left\|\int_{R}^{R+1}(\xi-R){\mathcal G(\xi e^{i\psi})}(\xi
e^{i\psi}-\lambda)^{-1}\,d\xi\right\|^2\leq
C\int_{R}^{R+1}\|\mathcal G(\xi e^{i\psi})\|^2\,d\xi,
$$
are valid, where the right hand sides tend to zero as
$R\to+\infty$ whenever $\mathcal G\in\mathscr H^0_0(\mathcal
K^\varphi_0;X)$ and $\psi\in(0,\varphi)$. We have proved the
equality
\begin{equation}\label{2cau}
\mathcal G(\lambda)=\frac 1 {2\pi
i}\left(\int_0^{+\infty}\frac{\mathcal
G(e^{i\psi}\xi)}{e^{i\psi}\xi-\lambda} e^{i\psi}\,d
\xi-\int_0^{+\infty}\frac{\mathcal
G(e^{i\phi}\xi)}{e^{i\phi}\xi-\lambda} e^{i\phi}\,d \xi\right),
\end{equation}
 where   $\arg\lambda\in(\psi,\phi)$, $\lambda\neq 0$, and
 $0<\psi<\phi<\varphi$.
Due to the assertion 1.(i) of Proposition~\ref{p2}, we can pass in
\eqref{2cau} to the limit as $\psi\to 0+$ and $\phi\to\varphi-$.
As a result we obtain the equality \eqref{apL1}. The
representation~\eqref{cau1} is proved.

(ii) To prove the assertion (ii) of Proposition~\ref{p2.6} it
suffices to show that the value $\|\mathcal G(\lambda)\|$
uniformly tends to zero as $\lambda$ goes to infinity in the set
$\{\lambda\in\mathcal K^\varphi_0:\dist\{\lambda,\partial \mathcal
K^\varphi_0\}\geq\epsilon>0\}$; here $\mathcal G\in\mathscr
H^0_0(\mathcal K^\varphi_0;X)$ is identified with $\mathcal
F\in\mathscr H^\ell_w(\mathcal K^\varphi_\zeta;X)$ by the
rule~\eqref{G}. For all $0<\varrho<\epsilon$ we have
$$
\mathcal G(\lambda)=\frac  1 {2\pi i}\int_{|\mu-\lambda|=\varrho}
\frac {\mathcal G(\mu)} {\mu-\lambda}\,d\mu=\frac 1
{2\pi}\int_0^{2\pi}\mathcal G(\lambda+\varrho
e^{i\alpha})\,d\alpha,
$$
\begin{equation}\label{for}
\pi\epsilon^2\mathcal
G(\lambda)=\int_0^\epsilon\int_0^{2\pi}\mathcal G(\lambda+\varrho
e^{i\alpha})\,d\alpha\,\varrho\,d\varrho.
\end{equation}
Using the Cauchy-Schwarz inequality, we deduce from \eqref{for}
that
\begin{equation}\label{cau6}
\pi^2\epsilon^4\|\mathcal G(\lambda)\|^2\leq 2\pi \int_0^\epsilon
\varrho\,d\varrho \cdot \int_0^\epsilon \int_0^{2\pi}\|\mathcal
G(\lambda+\varrho e^{i\alpha})\|^2\,d\alpha\,\varrho\,d\varrho.
\end{equation}
Since the circle $\{\mu\in\Bbb C:|\mu-\lambda|=\epsilon\}$ is
contained in the set
$$
\{\mu\in C: \arg \mu\in [\arg\lambda-
2\epsilon/|\lambda|,\arg\lambda+2\epsilon/|\lambda|],
|\lambda|-\epsilon\leq|\mu|\leq|\lambda|+\epsilon\},
$$
 we have
 \begin{equation}\label{cau7}
 \begin{aligned}
  \int_0^\epsilon \int_0^{2\pi}\|\mathcal
G(\lambda+\varrho e^{i\alpha})\|^2\,d\alpha\,\varrho\,d\varrho
\leq\int_{-2\epsilon/|\lambda|}^{2\epsilon/|\lambda|}\,d\psi\int_{|\lambda|-\epsilon}^{|\lambda|+\epsilon}
\|\mathcal G(\xi e^{i(\arg\lambda+\psi)})\|^2\xi\,d\xi
 \\
 =\int_{-2\epsilon}^{2\epsilon}\,d\psi \int_{|\lambda|-\epsilon}^{|\lambda|+\epsilon}
\|\mathcal G(\xi
e^{i(\arg\lambda+\psi/|\lambda|)})\|^2(\xi/|\lambda|)\,d\xi
\\
\leq C \int_{-2\epsilon}^{2\epsilon}\,d\psi
\int_{|\lambda|-\epsilon}^{|\lambda|+\epsilon} \|\mathcal G(\xi
e^{i(\arg\lambda+\psi/|\lambda|)})\|^2\,d\xi.
 \end{aligned}
 \end{equation}
 Due to the inclusion $\mathcal G\in\mathscr H^0_0(\mathcal
 K^\varphi_0;X)$ the integral $\int_{|\lambda|-\epsilon}^{|\lambda|+\epsilon} \|\mathcal G(\xi
e^{i(\arg\lambda+\psi/|\lambda|)})\|^2\,d\xi$ is uniformly bounded
in $\psi\in[-2\epsilon,2\epsilon]$ for all sufficiently large
$|\lambda|$ and tends to zero as $|\lambda|\to +\infty$. This
together with \eqref{cau6} and \eqref{cau7} finishes the proof.
\qed\end{pf*}

\begin{lem}\label{ApLemma2} For any analytic function $\Bbb
C^+\ni\lambda\mapsto \mathcal F(\lambda)\in X$ satisfying the
uniform in $\eta\in\Bbb C^+$ estimate
\begin{equation}\label{apL3} \|\mathcal
F(\cdot+\eta); L_2(\Bbb R;X)\|\leq C(\mathcal F)
\end{equation}
 there exists a function $\mathcal G\in
L_2(\Bbb R;X)$ such that for all $\lambda\in\Bbb C^+$ the equality
\begin{equation}\label{apL2}
\sqrt{\lambda}\mathcal F(\lambda)=\frac 1{2\pi
i}\int_{-\infty}^{+\infty} \frac {\mathcal G(\xi)}{(\log
\xi-\log\lambda)\sqrt{\xi}}\,d\xi
\end{equation}
is valid; here we use the analytic in $\overline{\Bbb
C^+}\setminus\{0\}$ principal branches of logarithm and square
root, the integral is absolutely convergent in $X$.
\end{lem}

\begin{pf}
Let $R>0$ and $0<\eta_1<\eta_2$. From the Cauchy integral theorem
it follows the equality
\begin{equation}\label{re}
\sqrt{\lambda}\mathcal F(\lambda)=\frac 1{2\pi i}\int_R^{R+1}\,da
\oint_{\mathscr C(a,\eta_1,\eta_2)} \frac {\mathcal F(\mu)}{(\log
\mu-\log\lambda)\sqrt{\mu}}\,d\mu,
\end{equation}
where $|\Re\lambda|<R$, $\Im\lambda\in (\eta_1,\eta_2)$, and the
closed contour $\mathscr C(a,\eta_1,\eta_2)$ is given by the
formula
$$
\begin{aligned}
\mathscr C(a,\eta_1,\eta_2)=\{\mu\in\Bbb
C:\mu=\xi+i\eta_1,\xi&\in[-a,a]\}\\ \cup\{\mu\in\Bbb
C:\mu=a+i\eta, \eta\in[\eta_1,\eta_2]\} \cup&\{\mu\in\Bbb
C:\mu=-\xi+i\eta_2,\xi\in[-a,a]\}\\\cup&\{\mu\in\Bbb
C:\mu=-a-i\eta, \eta\in[-\eta_2,-\eta_1]\}.
\end{aligned}
$$
 Let us check that   passing to the limit in \eqref{re} as
$R\to+\infty$ leads to the equality
\begin{equation}\label{re1}
\sqrt{\lambda}\mathcal F(\lambda)=\frac 1{2\pi i}\left(\int_{\Bbb
R+i\eta_1} \frac {\mathcal F(\mu)}{(\log
\mu-\log\lambda)\sqrt{\mu}}\,d\mu-\int_{\Bbb R+i\eta_2} \frac
{\mathcal F(\mu)}{(\log \mu-\log\lambda)\sqrt{\mu}}\,d\mu\right).
\end{equation}
Indeed, by Cauchy-Schwarz inequality we get
\begin{equation}\label{re3}
\begin{aligned}
\left(\int_{R}^{R+1}\,da\int_{\eta_1}^{\eta_2}\left\|\frac
{\mathcal F(\pm a+i\eta)}{(\log(\pm a+i\eta)-\log\lambda)\sqrt{\pm
a+i\eta}}\right\|\,d\eta \right)^2
\\
\leq (\eta_2-\eta_1)\int_{\eta_1}^{\eta_2} \left(\int_R^{R+1}
{|\log(\pm a+i\eta)-\log\lambda|^{-2}|\pm
a+i\eta|^{-1}}\,{da}\right)
\\
\times\left(\int_R^{R+1}\|\mathcal F(\pm a+i\eta)\|^2\,d a\right)
\,d\eta.
\end{aligned}
\end{equation}
Here the integral $\int_R^{R+1}\|\mathcal F(\pm a+i\eta)\|^2\,d a$
is uniformly bounded in $\eta\in\Bbb R^+$ and tends to zero as
$R\to+\infty$ due to the estimate \eqref{apL3}. It is easy to see
that
\begin{equation}\label{re5}
\begin{aligned}& \int_R^{R+1}{|\log(\pm a+i\eta)-\log\lambda|^{-2}|\pm
a+i\eta|^{-1}}\,{da}
\\
&\leq \int_R^{+\infty}(\log a-\log|\lambda|)^{-2} a^{-1}
\,da=1/\log (R/|\lambda|)\to 0, \ R\to+\infty.
\end{aligned}
\end{equation}
Thus the right hand side of the inequality \eqref{re3} tends to
zero as $R\to+\infty$.

 We now consider the integral
\begin{equation}\label{re4}
\begin{aligned}
\int_R^{R+1}&\int_{-a}^a\frac{\mathcal
F(\xi+i\eta)}{(\log(\xi+i\eta)-\log\lambda)\sqrt{\xi+i\eta}}\,d\xi
\,d a
=\int_{-R}^R[\dots]\,d\xi\\&+\int_{-R-1}^{-R}(R+1+\xi)[\dots]\,d\xi+\int_R^{R+1}(R+1-\xi)[\dots]\,d\xi,
\end{aligned}
\end{equation}
where, for brevity, $[\dots]$ denotes the same integrand as in the
left hand side of the equality \eqref{re4}, and $\eta=\eta_1$ or
$\eta=\eta_2$. The interchange of integrations in \eqref{re4} is
legal because of \eqref{apL3} and the estimate
$$
\int_{-\infty}^{+\infty}{|\log(\xi+i\eta)-\log\lambda|^{-2}|\xi+i\eta|^{-1}}\,{d\xi}
\leq Const(\eta,\lambda),\quad\lambda\notin\Bbb R+i\eta.
 $$
  Note
that the last two integrals in \eqref{re4} go to zero in $X$ as
$R\to+\infty$. Indeed, for the first of these integrals it follows
from \eqref{re5}, \eqref{apL3}, and the estimate
$$
\left\|\int_{-R-1}^{-R}(R+1+\xi)\frac{\mathcal
F(\xi+i\eta)}{(\log(\xi+i\eta)-\log\lambda)\sqrt{\xi+i\eta}}\,d\xi\right\|^2
$$
$$\leq
\left(\int_{-R-1}^{-R}{|\log(\xi+i\eta)-\log\lambda|^{-2}|\xi+i\eta|^{-1}}\,d\xi\right)\left(\int_{-R-1}^{-R}\|\mathcal
F(\xi+i\eta)\|^2\,d\xi\right);
$$
a similar estimate is valid for the last integral in \eqref{re4}.

It is easily seen that the first integral in the right hand side
of the equality \eqref{re4} tends to the first integral in the
formula \eqref{re1} (as $R\to+\infty$) if $\eta=\eta_1$ and to the
second integral if $\eta=\eta_2$. We have shown that the right
hand side of the equality \eqref{re} tends to the right hand side
of \eqref{re1} in the space $X$ as $R\to+\infty$. The equality
\eqref{re1} is proved.

 On the next step we pass in \eqref{re1} to the limit as
$\eta_1\to 0+$ and then as $\eta_2\to+\infty$. As a result we get
the equality \eqref{apL2} and complete the proof.

Let us establish the uniform in $\eta\in[0,\epsilon]$
 estimate
\begin{equation}\label{star}
\int_{-\infty}^{+\infty}{|\log(\xi+i\eta)-\log\lambda|^{-2}|\xi+i\eta|^{-1}}\,{d\xi}
\leq Const(\epsilon,\lambda),
\end{equation}
where $\lambda\in\Bbb C^+$ and $\epsilon$ is a sufficiently small
positive number,  $\epsilon<<\Im\lambda$. If $|\xi|<\epsilon$,
$\eta\in[0,\epsilon]$, and $\epsilon$ is sufficiently small then
$|\xi+i\eta|<\exp(-2\max\{1,|\log\lambda|\})$ and the inequalities
$$
|\log(\xi+i\eta)-\log\lambda|\geq
|\log|\xi+i\eta||-|\log\lambda|\geq |\log|\xi+i\eta||/2,
$$
$$
|\log(\xi+i\eta)-\log\lambda|^2|\xi+i\eta|\geq
(\log|\xi+i\eta|)^2|\xi+i\eta|/4\geq(\log|\xi|)^2|\xi|/4
$$
hold. Let $R$ be a fixed number, $R>\max\{1,|\lambda|\}$. For
$\xi\in[-R,-\epsilon]\cup[\epsilon, R]$ we have
$$
|\log(\xi+i\eta)-\log\lambda|^{-2}|\xi+i\eta|^{-1}\leq
C(\epsilon), \quad\eta\in[0,\epsilon].
$$
Now the estimate \eqref{star} follows from the inequalities
$$
\int_{-\infty}^{+\infty}{|\log(\xi+i\eta)-\log\lambda|^{-2}|\xi+i\eta|^{-1}}\,{d\xi}\leq
8\int_0^{\epsilon} (\log|\xi|)^{-2}|\xi|^{-1}\,d\xi$$
$$+2\int_{\epsilon}^RC(\epsilon)\,d\xi
+2\int_R^{+\infty}(\log\xi-\log|\lambda|)^{-2}\xi^{-1}\,d\xi\leq
Const(\epsilon,\lambda).
$$
Since  $(\log(\xi+i\eta)-\log\lambda)^{-1}(\xi+i\eta)^{-1/2}$
tends to $(\log\xi-\log\lambda)^{-1}\xi^{-1/2}$ as $\eta\to 0+$
for $\xi\in\Bbb R\setminus\{0\}$, and the estimate \eqref{star}
holds, we conclude that
\begin{equation}\label{rel1}
\|(\log(\cdot+i\eta)-\log\lambda)^{-1}(\cdot+i\eta)^{-1/2}-(\log\cdot-\log\lambda)^{-1}(\cdot)^{-1/2};L_2(\Bbb
R)\|\to 0
\end{equation}
as $\eta\to 0+$.
 Proposition~\ref{p1} asserts that
every analytic function $\Bbb C^+\ni\lambda\mapsto\mathcal
F(\lambda)\in X$ satisfying the uniform estimate \eqref{apL3} has
non-tangential boundary limits $\mathcal F_0\in L_2(\Bbb R;X)$ and
\begin{equation}\label{rel2}
\|\mathcal F(\cdot+i\eta)-\mathcal F_0(\cdot);L_2(\Bbb R;X)\|\to
0,\quad\eta\to 0+.
\end{equation}
 Having \eqref{rel1} and \eqref{rel2} we can pass in
\eqref{re1} to the limit as $\eta_1\to 0+$. We have
\begin{equation}\label{rel3}
\sqrt{\lambda}\mathcal F(\lambda)=\frac 1{2\pi i}\left(\int_{\Bbb
R} \frac {\mathcal F_0(\mu)}{(\log
\mu-\log\lambda)\sqrt{\mu}}\,d\mu-\int_{\Bbb R+i\eta_2} \frac
{\mathcal F(\mu)}{(\log \mu-\log\lambda)\sqrt{\mu}}\,d\mu\right).
\end{equation}
It remains to show that the second integral in \eqref{rel3} goes
to zero in $X$ as $\eta_2\to+\infty$. Let $R>\max\{1,|\lambda|\}$.
For all $\eta>R$ we have
$$
\begin{aligned}
\int_{-\infty}^{+\infty}{|\log(\xi+i\eta)-\log\lambda|^{-2}|\xi+i\eta|^{-1}}\,{d\xi}
\leq \int_{-R}^R(\log|\eta|-\log|\lambda|)^{-2}\,
d\xi\\+2\int_R^{+\infty}(\log\xi-\log|\lambda|)^{-2}\xi^{-1}\,d\xi\leq
Const(\lambda,R).
\end{aligned}
$$
As a consequence  we get
\begin{equation*}
\|(\log(\cdot+i\eta)-\log\lambda)^{-1}(\cdot+i\eta)^{-1/2};L_2(\Bbb
R)\|\to 0, \quad \eta\to +\infty.
\end{equation*}
This together with \eqref{apL3} allows us to pass in \eqref{rel3}
to the limit as $\eta_2\to+\infty$. The equality \eqref{apL2} is
valid for $\mathcal G\equiv\mathcal F_0$. \qed\end{pf}


\end{document}